\documentclass{article}

\usepackage{comment}
\usepackage{amsfonts}
\usepackage{graphics}
\usepackage{color}
\usepackage{amssymb}
\usepackage{verbatim}
\usepackage{amsmath}
\usepackage{esint}
\usepackage{amssymb}
\usepackage{amsthm}
\usepackage{amscd}
\usepackage{euscript}

\usepackage{mathtools}

\usepackage[all]{xy}

\usepackage{hyperref}

\newtheorem{Theorem}{Theorem}

\newtheorem{Corollary}[Theorem]{Corollary}

\newtheorem{Definition}[Theorem]{Definition}
\newtheorem{Remark}[Theorem]{Remark}

\newtheorem{Lemma}[Theorem]{Lemma}

\newtheorem{Proposition}[Theorem]{Proposition}

\newtheorem{Fundamental Theorem}{Fundamental Theorem}

\newenvironment{Proof}[1][Proof:\,\,]{\hskip-0.5cm\textbf{#1}}{$\qed$}

\usepackage[all]{xy}
\def \ch {\mathfrak{ch}}

\def \String {{\rm String}}

\def \aut {{\rm Aut}}
\def \Aut {{\mathcal{AUT}}}
\def \M {\EuScript{M}}
\def \hom {{\mathrm{hom}}}
\xyoption{v2}
\xyoption{all}
\xyoption{2cell}

\def \D {\Delta}
\def \tn {\otimes}
\def \BA {{\mathcal{BA}}}

\def \Prim {{\mathrm{Prim}}}
\def \t {\triangleright}

\def \ad {\t_{\rm ad}}
\def \tr {\t_{\rm \r}}
\def \ta {\t_{\rm ad}}

\def \tla {\lessdot}
\def \tra {\gtrdot}

\def \U {\EuScript{U}}

\def \b {\beta}
\def \dis  {\displaystyle}
\def \e {\epsilon}
\def \l {\lambda}

\def \ra {\xrightarrow}
\def \V {{\mathcal{V}}} 
\def \d {\partial}
\def \H {\mathcal{H}}

\def \r {\rho}
\def \T {\EuScript{T}}
\def \Tc {\mathcal{T}}
\def \P {\EuScript{P}}
\def \h {[[h]]}

\def \id {\mathrm{id}}

\def \w {\omega}
\def \A {\EuScript{A}}
\def \F {\EuScript{F}}
\def \Fc {\mathcal{F}}

\def \g {\gamma}
\def \N {\mathbb{N}}
\def \G {\Gamma}
\def \dt {\frac{\d}{\d t}}
\def \ds {\frac{\d}{\d s}}
\def \du {\frac{\d}{\d u}}

\def \dx {\frac{\d}{\d x}}

\def \bw {{\overline{\w}}}
\def \bm {\overline{m}}

\def \C {\mathbb{C}}
\def \Hh {H[[h]]}
\def \Ih {I[[h]]}
\def \Ch {\C[[h]]}
\def \sn {\sum_{n=0}^{+ \infty}}
\def \X  {{\cal X}}
\def \fo {\textrm{ for each }}
\def \an {\textrm{ and }}
\def \Der {{\mathrm{Der}}}
\def \lg {{\mathfrak{g}}}
\def \Ad {A_\bullet}

\def \Cc {\mathcal{C}}

\def \der {{\rm der}}

\def \nss {\hspace*{-.6cm}}
\def \ae {\mathfrak{e}}

\def \Lie {\mathrm{Lie}}
\def \i {\iota}
\def \ve {\xi}
\def \gl {{\mathfrak{gl}}}
\def \lgl {{\mathfrak{GL}}}
\def \GL {{\mathrm{GL}}}
\def \we {{\textrm{ where }}}
\def \Fz {{{\cal Z}}}

\makeatother
\newcommand{\mortimes}[3]{{   {\xymatrix{ & {*} \ar@/^1pc/[rr]^{{#3} }\ar@/_1pc/[rr]_{{#1}} & \boxed{\times}\Uparrow {#2} & {*}}}}}
\newcommand{\bmortimes}[3]{{   {\xymatrix{ & {*} \ar@/^2pc/[rrrr]^{{#3} }\ar@/_2pc/[rrrr]_{{#1}} && \boxed{\times}\Uparrow {#2} && {*}}}}}

\newcommand{\bbmortimes}[3]{{   {\xymatrix{ & {*} \ar@/^1.1pc/[rrrr]^{{#3} }\ar@/_1.1pc/[rrrr]_{{#1}} && \boxed{\times}\Uparrow {#2} && {*}}}}}

\newcommand{\mort}[5]{{   {\xymatrix{ & {#1} \ar@/^1pc/[rr]^{{#4} }\ar@/_1pc/[rr]_{{#2}} & \Uparrow {#3} & {#5}}}}}
\newcommand{\morts}[5]{{  \hskip-1cm {\xymatrix{ & {#1} \ar@/^0.5pc/[r]^{{#4} }\ar@/_0.5pc/[r]_{{#2}} \ar@{{}{ }{}} [r]|{\Uparrow#3}  & {#5}}}}}

\newcommand{\mortimesvcomp}[5] {{ \xymatrix{ & {*}\ar@/^2pc/[rr]^{{#5}}\ar[rr]|{{#3}} \ar@/_2pc/[rr]_{{#1}} \ar @/^/ @{{}{ }{}} [rr]^{\boxed{\times}\Uparrow {#4}} \ar @/_/ @{{}{ }{}} [rr]_{\boxed{\times}\Uparrow{#2}}
 && {*}}}}

\newcommand{\mortimesvvcomp}[5] {{ \xymatrix{ & {*}\ar@/^3pc/[rrrr]^{{#5}}\ar[rrrr]|{{#3}} \ar@/_3pc/[rrrr]_{{#1}} \ar @/^/ @{{}{ }{}} [rrrr]^{\boxed{\times}\Uparrow {#4}} \ar @/_/ @{{}{ }{}} [rrrr]_{\boxed{\times}\Uparrow{#2}}
 &&&& {*}}}}
\newcommand{\mortimesvhcomp}[9] {{
 \xymatrix{ 
& {*}\ar@/^2pc/[rr]^{{#5}}\ar[rr]|{{#3}} \ar@/_2pc/[rr]_{{#1}} \ar @/^/ @{{}{ }{}} [rr]^{\boxed{\times}\Uparrow {#4}} \ar @/_/ @{{}{ }{}} [rr]_{\boxed{\times}\Uparrow{#2}}
 && {*} \ar@/^2pc/[rr]^{{ \d({f'}^{-1})\d({e'}^{-1}) g'   }}\ar[rr]|{{#8}} \ar@/_2pc/[rr]_{{#6}} \ar @/^/ @{{}{ }{}} [rr]^{\boxed{\times}\Uparrow {#9}} \ar @/_/ @{{}{ }{}} [rr]_{\boxed{\times}\Uparrow{#7}}
 && {*} }
}}

\newcommand{\rwtimes}[4] {{   {\xymatrix{ & {*} \ar@/^1pc/[rr]^{{#3} }\ar@/_1pc/[rr]_{{#1}} &\boxed{\times} \Uparrow {#2} & {*} \ar[r]^{{#4}} &{*} }}}}

\newcommand{\lwtimes}[4] {{   {\xymatrix{ & {*} \ar[r]^{{#1}}  & {*} \ar@/^1pc/[rr]^{{#4}} \ar@/_1pc/[rr]_{{#2}} & \boxed{\times}\Uparrow {#3} & {*}  } } 
                    }}

\newcommand{\mortimeshcomp}[6]{{   {\xymatrix{ & {*} \ar@/^1pc/[rr]^{{#3} }\ar@/_1pc/[rr]_{{#1}} & \boxed{\times}\Uparrow {#2} & {*}  \ar@/^1pc/[rr]^{{#6} }\ar@/_1pc/[rr]_{{#4}} & \boxed{\times}\Uparrow {#5} & {*} }}}}

\newcommand{\morplus}[3]{{   {\xymatrix{ & {*} \ar@/^1pc/[rr]^{{#3} }\ar@/_1pc/[rr]_{{#1}} & \boxed{+}\Uparrow {#2} & {*}}}}}

\newcommand{\morplush}[3]{{   {\xymatrix{ & {*} \ar@/^1.5pc/[rr]^{{#3} }\ar@/_1.5pc/[rr]_{{#1}} & \boxed{+}\Uparrow {#2} & {*}}}}}

\newcommand{\morplusvcomp}[5] {{ \xymatrix{ & {*}\ar@/^2pc/[rr]^{{#5}}\ar[rr]|{{#3}} \ar@/_2pc/[rr]_{{#1}} \ar @/^/ @{{}{ }{}} [rr]^{\boxed{+}\Uparrow {#4}} \ar @/_/ @{{}{ }{}} [rr]_{\boxed{+}\Uparrow{#2}}
 && {*}}}}

\newcommand{\rwplus}[4] {{   {\xymatrix{ & {*} \ar@/^1pc/[rr]^{{#3} }\ar@/_1pc/[rr]_{{#1}} &\boxed{+} \Uparrow {#2} & {*} \ar[r]^{{#4}} &{*} }}}}

\newcommand{\lwplus}[4] {{   {\xymatrix{ & {*} \ar[r]^{{#1}}  & {*} \ar@/^1pc/[rr]^{{#4}} \ar@/_1pc/[rr]_{{#2}} & \boxed{+}\Uparrow {#3} & {*}  } } 
                    }}

\newcommand{\morplushcomp}[6]{{   {\xymatrix{ & {*} \ar@/^1pc/[rr]^{{#3} }\ar@/_1pc/[rr]_{{#1}} & \boxed{+}\Uparrow {#2} & {*}  \ar@/^1pc/[rr]^{{#6} }\ar@/_1pc/[rr]_{{#4}} & \boxed{+}\Uparrow {#5} & {*} }}}}
\newcommand{\bmorplus}[3]{{   {\xymatrix{ & {*} \ar@/^2pc/[rrrr]^{{#3} }\ar@/_2pc/[rrrr]_{{#1}} && \boxed{+}\Uparrow {#2} && {*}}}}}

\newcommand{\bmorplush}[3]{{   {\xymatrix{ & {*} \ar@/^3pc/[rrrr]^{{#3} }\ar@/_3pc/[rrrr]_{{#1}} && \boxed{+}\Uparrow {#2} && {*}}}}}
\makeatletter

\allowdisplaybreaks
\begin{document}
\title{Crossed modules of Hopf algebras and of associative algebras and  two-dimensional holonomy}

\author{
Jo\~{a}o  Faria Martins\thanks{ This work was supported by CMA/FCT/UNL, under the project UID/MAT/00297/2013  and by FCT(Portugal) through the  ``Geometry and Mathematical Physics Project'',
FCT EXCL/MAT-GEO/0222/2012.} \thanks{I would like to thank Camilo Arias-Abad and Kadir Emir for useful discussions. {I would like to express my gratitude to the referee for thorough comments.}   
} \\ 
{\footnotesize Departamento de Matem\'{a}tica  and Centro de Matem\'{a}tica e Aplica\c{c}\~{o}es,} \\ 
{\footnotesize Faculdade de Ci\^{e}ncias e Tecnologia,   Universidade NOVA de Lisboa} \\
{\footnotesize  Quinta da Torre, 2829-516 Caparica, Portugal } \\
{\it \small jn.martins@fct.unl.pt } }
\date{\today}
\maketitle
\begin{abstract}
After a thorough treatment of all algebraic structures involved, we address two dimensional holonomy operators with values in crossed modules of Hopf algebras and in crossed modules of associative algebras (called here crossed modules of bare algebras.) In particular, we will consider two general formulations of the {two-dimensional} holonomy of a (fully primitive) Hopf 2-connection (exact and blur), the first being multiplicative the second being additive, proving that they coincide in a certain natural quotient (defining what we called the fuzzy holonomy of a fully primitive Hopf 2-connection). 
\end{abstract}
\maketitle
\tableofcontents
\section{Introduction}
A crossed module of groups $\X=(\d\colon E \to G, \t)$ is given by a map $\d\colon E \to G$ of groups, together with an action $\t$ of $G$ on $E$ by automorphisms. This action must satisfy two natural properties (the Peiffer relations):
\begin{align*}
   \d(g \t e) &=g\, \d(e)\, g^{-1}, \quad \fo g \in G, \an e \in E &&\textrm{(first Peiffer relation)},\\
   \d(e) \t f &= e\, f\, e^{-1},\quad  \fo e,f \in E &&\textrm{(second Peiffer relation)}.
  \end{align*}
Note that the second Peiffer law implies that $\ker(\d)$ is an abelian subgroup of $E$.  These are very flexible axioms. For example, given crossed modules $\X=(\d\colon E \to G, \t)$ and  $\X'=(\d\colon E' \to G', \t)$, then $\X \times \X'=\big((\d\times \d')\colon E\times E' \to G \times G', \t\times \t')$, with the obvious product action $\t \times \t'$, is a crossed module. On the other hand, if $L$ is a subgroup of $\ker(\d) \subset E$, such that $G \t L \subset L$, which implies that $L$ is normal (because of the second Peiffer relation), then 
  $\X/L=(\d\colon E/L  \to G, \t)$ with the obvious quotient action is also a crossed module. We note that morphisms of crossed module are defined in the obvious way. 

Crossed modules of groups were invented by Whitehead in \cite{WJHC}, naturally appearing in the context of homotopy theory, being algebraic models for homotopy 2-types \cite{MLW}  (for a modern account of this see \cite{BHS,Baues}). Given a pointed fibration $F \to E \to B$, then the inclusion of the fibre $F$ in the total space $E$ induces a crossed module $\big(\pi_1(F) \to \pi_1(E)\big)$. A pointed pair of spaces $(M,N)$ also has a fundamental crossed module $\big(\pi_2(M,N) \to \pi_1(N)\big)$, with the obvious boundary map and action of $\pi_1$ on $\pi_2$. The relation between group crossed modules and strict 2-groups (small categories whose sets of objects and of morphisms are groups, with all structure maps, including the  composition, being group morphisms) was elucidated in \cite{BrS}, where we can find the first reference to the fact that the categories of crossed modules and of strict 2-groups are equivalent. This relation was fully generalized for crossed complexes and $\omega$-groupoids in \cite{BHiggins}. 

We can see a (strict) 2-group as being a small 2-category with a single object where all morphisms are invertible; \cite{BL}. Given a crossed module $\X=(\d\colon E \to G, \t)$, the 1-morphisms of the associated 2-group, denoted by $\Cc^\times(\X)$, have the form $*\ra{g} *$, where $g \in G$, with the obvious composition. The 2-morphisms of  $\Cc^\times(\X)$ have the form below, composing vertically and horizontally (for conventions see \ref{gcmc}): $$\mortimes{g}{e}{\d(e)^{-1} g},\we g \in G, \an e \in E .$$

A general theory of (Lie) 2-groups (including non-strict ones) appears in \cite{BL}. This theory parallels the theory of Lie-2-algebras (strict and non-strict), which was developed in  \cite{BC}. We mention that a strict Lie 2-algebra is uniquely represented by a crossed module $\mathfrak{X}=(\d\colon \ae \to \lg, \t)$ of Lie algebras (also called a differential 
crossed module), where  $\t$ is a left action of $\lg$ on $\ae$ by derivations, and $\d\colon \ae \to \lg$ is a Lie 
algebra map, furthermore satisfying the (obvious) differential Peiffer relations. Given a crossed module $\X=(\d\colon E \to G, \t)$, of Lie groups, the induced Lie algebra map $\d\colon \ae \to \lg$, and action by derivations of $\lg$ on $\ae$, defines a crossed module  $\mathfrak{X}=(\d\colon \ae \to \lg, \t)$, of Lie algebras.

On the purely algebraic side, crossed modules of (Lie) groups and of Lie algebras arise in a variety of ways. Given a Lie group $G$, we have the actor crossed module of $G$, being $\Aut(G)=(G \ra{\rm ad} \aut(G),\t)$, where $\aut(G)$ is the Lie group of automorphisms of $G$, acting in $G$ in the obvious way  and ${\rm ad}$ is the map that sends every $g\in G$ to the map ${\rm ad}_g\colon G \to G$ which is conjugation by $g$.
If we have a central extension $\{1\}\to A \ra{i} B \ra{\d} K\to \{1\}$ of groups, given any section $s\colon K\to B$ of $\d\colon B \to K$, then $k \t b\doteq s(k)\, b\, s(k)^{-1}$, where $k\in K$ and $b \in B$, is an action of $K$ on $B$, by automorphisms independent of the chosen section, defining a crossed module with underlying group map being $\d\colon B \to K$.  

To a chain-complex $\V=(\dots \ra{\b} V_n \ra{\b} V_{n-1}\ra{\b}\dots)$ of vector spaces we can associate crossed modules of Lie algebras and (in the finite dimensional case) of Lie groups, denoted respectively by $\lgl(\V)=\big(\d\colon \gl_1(\V) \to \gl_0(\V),\t)$ and by $\GL(\V)$; see \cite{BL,FMM}. The Lie algebra $\gl_0(\V)$ is made out of chain-maps $\V\to \V$, with the usual bracket. On the other hand, the Lie algebra $\gl_1(\V)$ is given by homotopies (degree one maps) $\V \to \V$, up to 2-fold homotopies.  (The bracket is not the commutator of the underlying linear maps, the latter having degree two.) The 2-group associated to the crossed module $\GL(\V)$ has a single object, has morphisms being the invertible chain maps $\V \to \V$ and 2-morphisms being the chain homotopies (up to 2-fold homotopy) between them. The vertical composition of 2-morphisms is induced by the sum of homotopies. 

We can also consider pre-crossed modules of groups and of Lie algebras. These are defined similarly to crossed modules, however not imposing the second Peiffer condition. Crossed modules form a full subcategory of the category of pre-crossed modules (both in the group and Lie algebra cases). Moreover, any pre-crossed module can  naturally be converted to a crossed module by dividing out the second Peiffer relations. This defines reflection functors $\{\textrm{pre-crossed modules}\} \to \{\textrm{crossed modules}\}.$

Crossed modules of groups and of Lie algebra  naturally have \textit{free} objects. Let $\lg$ be a Lie algebra. Let also $E$ be a set, with a map $\d_0\colon E \to \lg$. We have a free differential crossed module $\mathfrak{F}=(\d\colon \mathfrak{f} \to \lg)$ on the map $\d_0\colon E \to \lg$. Here $\d\colon \mathfrak{f} \to \lg$ is a Lie algebra map, and we have a (set) inclusion $i\colon E \to \mathfrak{f}$, satisfying $\d \circ i =\d_0$. This free crossed module $\mathfrak{F}$ satisfies (and is defined by) the following universal property: if we have a differential crossed module $\mathfrak{X}'=(\d'\colon \ae' \to \lg',\t)$,  a Lie algebra map $f\colon \lg \to \lg'$ and a set map $g_0\colon E \to \ae'$, satisfying $\d'\circ g_0 = f \circ \d_0$, then there exists a unique Lie algebra map $g\colon \mathfrak{f} \to \ae'$, extending $g_0$, and such that, furthermore, the pair $\big(g\colon  \mathfrak{f} \to \ae', f\colon \lg \to \lg')$ is a map of differential crossed modules $\mathfrak{F} \to \mathfrak{X}'$. 

Models for the free differential crossed module on a set map $\d_0\colon E \to \lg$ appear in \cite{El,CFM1}.
 Consider the action $\t$ of $\lg$ on $\U(\lg)$, the universal enveloping algebra of $\lg$ by left multiplication. Consider also the vector space $\oplus_{e \in E} \U(\lg).e$, with the obvious action $\t$ of $\lg$. Consider the linear map $\d_1\colon \oplus_{e \in E} \U(\lg).e \to \lg$ such that {$\d_1(a.e)=a\ad \d_0(e)$} (here $a \in \U(\lg)$ and $e \in E$), where $\ad$ is the action of $\U(\lg)$ in $\lg$ induced by the adjoint action of $\lg$ on $\lg$. Clearly $\d_1(X \t (a.e))=[X,\d_1(a.e)]$, where $X \in \lg$, $a \in \U(\lg)$ and  $ e \in E$. Let $\mathfrak{f}_1$ be the free Lie algebra on the vector space $ \oplus_{e \in E} \U(\lg).e$. We therefore have a Lie algebra map $\d_2\colon \mathfrak{f}_1 \to \lg$, extending $\d_1$. Since any linear map $V \to V$, $V$ a vector space, induces a derivation at the level of the free Lie algebra on $V$, we have a Lie algebra action $\t$ of $\lg $ on $\mathfrak{f}_1$, by derivations. Clearly this defines a pre-crossed module of Lie algebras. The free crossed module on 
the map $\d_0\colon E \to \lg$ is obtained by converting the latter differential pre-crossed module into a crossed module of Lie algebras.

Crossed modules (of groups and of Lie algebras) are closely related to cohomology; see \cite{WoW,KB,Ger}. Namely, in the group case (not considering any topology),  given a group $K$ and a $K$-module $A$, there is a one-to-one correspondence between group cohomology classes $\w\in H^3(K,A)$ and weak homotopy classes of  crossed modules $(\d\colon E \to G, \t)$, fitting inside the exact sequence $\{0\} \to A \to E \ra{\d} G \to K\cong G/\d(E) \to \{1\}$. We actually have a group isomorphism considering   Baer sums of  crossed modules; \cite{KB}.

The same correspondence holds for Lie algebras, where we have a one-to-one correspondence between cohomology classes $\w \in H^3(\mathfrak{k},\mathfrak{a})$, and weak equivalence classes of crossed modules of Lie algebras  $\mathfrak{X}=(\d\colon \ae \to \lg, \t)$, fitting inside the exact sequence $\{0\} \to \mathfrak{a} \to \ae \ra{\d} \lg \to \lg/\d(\ae)\cong \mathfrak{k}$; see \cite{W,Ger,BC}. 

One has to be careful that, in the Lie group case, there are several different definitions of the Lie group cohomology $H^3_*(K,A)$, for a Lie group $K$ acting on a vector space $A$. For a very nice review on these topics see \cite{WoW,BL}. A version of Lie-group cohomology $H_g(K,A)$ whose obvious linearization yields an isomorphism  $H^3_g(K,A) \cong 
H^3(\mathfrak{k},\mathfrak{a})$  is obtained by considering germs of group cocycles $\w\colon K^n \to A$, where a cocycle germ satisfies the same laws of the usual group cohomology, except that it is only defined in a neighborhood of the identity (in addition to being smooth), and two of these considered to be equivalent if they coincide in a neighborhood of the identity. There is a map $H_g^n(K,A) \to H^n(\mathfrak{k},\mathfrak{a})$ induced by derivation which is an isomorphism. If $K$ is compact, these local group cocycles never extend to  everywhere defined,  cohomologically non-trivial, continuous cocycles. The reason is that if $\w\colon K^{n+1} \to K$ is a continuous group cocycle then $(a_1,\dots a_n) \in G^n  \stackrel{\alpha}{\mapsto}\int_G \w(a_1,\dots,a_n,g)dg$ is such that $d\alpha=\pm\w$. Nevertheless, for integral cohomology classes, one can find extensions to globally defined cocycles, however not continuous outside of a neighborhood of the identity. These integral cohomology classes $\w$ 
can frequently be 
realized by  crossed modules $(P \to B,\t)$, of (infinite dimensional) Lie groups, fitting inside the exact sequence: $\{a\} \to A \to P \to B \to K \to \{1\}$.

Consider a compact semisimple Lie group $G$ with Lie algebra $\lg$. We have a cohomology class $\w\colon \lg^{\tn 3} \to \mathbb{C}$, given by $\w(X,Y,Z)=\langle [X,Y],Z \rangle$, where $\langle,\rangle$ is the Cartan-Killing form. A  geometric realization of the class $v=2\pi i \w$ as being the differential crossed module associated to a crossed module of infinite dimensional Lie groups, denoted by  $\String$, appears in \cite{BSAS}. In the ${\rm SU}(2)\cong S^3$ case, this crossed module of Lie groups fits inside the exact sequence $\{0\} \to U(1) \ra{i} \widehat{\Omega(S^3)} \ra{\d} \mathcal{P}(S^3) \ra{p} {\rm SU}(2)$. Here
$ \mathcal{P}(S^3)$ denotes the group of smooth paths $[0,1] \to S^3$, starting at the identity, with point-wise product, $p$ being the map choosing the final point of the path; $\widehat{\Omega(S^3)}$ is the  
Kac-Moody group, defining (see \cite{Mu,PS}) a universal central extension $U(1) \to \widehat{\Omega(S^3)} \ra{p}  \Omega(S^3)$ of the loop group $\Omega(S^3)$,  the group of smooth paths $[0,1] \to S^3$, starting and ending at the identity. The boundary map $\d\colon  \widehat{\Omega(S^3)}  \to \mathcal{P}(S^3)$ is  the composition of the projection map $p\colon \widehat{\Omega(S^3)} \to \Omega(S^3)$ and the inclusion map $\Omega(S^3) \to \mathcal{P}(S^3)$.

A general theory of 2-bundles was initiated in \cite{BS}, with very strong emphasis on their two dimensional holonomy. These 2-bundles have a structure Lie 2-group (equivalent to a Lie  crossed module), which (as far as two dimensional holonomy is concerned) is taken to be strict, however not necessarily finite dimensional. We note that 2-bundles with connection are closely related with gerbes with connection; see \cite{BreenMessing}. Further work on the holonomy of a 2-bundle with a 2-connection  appears in \cite{SW1,SW2,FMP1,FMP2,FMP3}.

Let $\X=(\d\colon E \to G,\t)$ be a crossed module of  Lie groups, with associated differential crossed module 
$\mathfrak{X}=(\d\colon \ae \to \lg, \t)$. Given a 2-bundle $P$ over a manifold $M$, with structure 2-group $\X$, a 2-connection in $P$ is locally given by a 1-form $\w$ in $M$ with values in {$\lg$}, and an $\ae$-valued 2-form $m$, in $M$, such that $\d(m)= \F_\w$, the curvature of $\w$, in our conventions:
 $\F_\w=d \w -\frac{1}{2} \w \wedge \w .$ The 2-curvature 3-form of the pair $(\w,m)$ is, by definition, $\M=d m - \w \wedge^\t m.$ A 2-connection is said to be flat if its 2-curvature 3-form vanishes. 
 
 Let us roughly explain what we mean by the two dimensional holonomy of a local 2-connection $(\w,m)$, referring details to the text and the references mentioned above. 
Let $M$ be a manifold. We denote by $\P^1(M)$ the space of continuous and piecewise smooth paths $\g\colon [0,1] \to M$ ({1-path}s), represented as: 
$$x=\g(0) \ra{\g} \g(1)=y. $$
If $\g_1$ and $\g_2$ are 1-paths with $\g_1(1)=\g_2(0)$, we define their concatenation, written as $\left (x \ra{\g_1} y \ra{\g_2} z\right)=x\ra{\g_1 \g_2 } z,$  in the obvious way. We  define a 2-path as being  a map $\G\colon [0,1]^2 \to M$, which is to  be continuous and piecewise smooth, for some paving of the square by polygons, transverse to the boundary of the square. We impose, moreover, that $\G(\{0,1\} \times [0,1])$ has at most two elements. Putting $\{x\}=\G(\{0\} \times [0,1])$, $\{y\}=\G(\{1\} \times [0,1])$, and also $\g_1(t)=\G(t,1)$, $\g_0(t)=\G(t,0)$, where $t \in [0,1]$, defining {1-path}s $\g_0,\g_1\colon [0,1] \to M$, we denote the  2-path $\G$ in the (very suggestive) form: $$\G=\hskip-1cm\mort{x}{\g_0}{\G}{\g_1}{y} .$$
 The set of 2-paths in $M$ is denoted by $\P^2(M)$. {Clearly 2-paths compose horizontally and vertically, having obvious horizontal and vertical reverses. We also have whiskerings of 2-paths by 1-paths.}

 { Let $\P(M)$ be the triple $(\P^2(M),\P^1(M),M)$, with the natural compositions, boundaries and reverses.  If  $\X=(\d\colon E \to G,\t)$ is a crossed module of  Lie groups, with tangent differential crossed module
$\mathfrak{X}=(\d\colon \ae \to \lg, \t)$, then a  pair of differential forms $\left(\w\in\Omega^1(M,\lg), m\in \Omega^2(M,\ae)\right)$, in $M$, with $\d(m)=\F_\w$, defines a  local holonomy operator ${\Fc^\times_{(\w,m)}}$. This operator assign to every 2-path $\G$ a 2-morphism of the 2-group $\Cc^\times(\X)$ associated with $\X$. Denote it by:} \begin{equation}
 {\Fc^\times_{(\w,m)}}\left(
{\hskip-1cm
\xymatrix
{ 
& {x} \ar@/^1pc/[rr]^{{\g_2} }\ar@/_1pc/[rr]_{\g_1} & \Uparrow {\G} & {y}
}
} \right)=\hskip-1cm\bbmortimes{\left(P_\w(\g_1)\right)^{-1}}{ \left( Q_{(\w,m)}(\G)\right)^{-1}}{\left(P_\w(\g_2)\right)^{-1}}  .
\end{equation}
{The operator ${\Fc^\times_{(\w,m)}}$ preserves all compositions, boundaries and reverses, in $\P(M)$ and in $\Cc^\times(\X)$.} 

In the context of graded vector bundles, a different type of higher dimensional holonomy (of representations up to homotopy) was developed in \cite{AC,AF1}. In a natural (non-trivial) quotient (cotruncation in dimension two) the two-dimensional part of the holonomy of a representation up to homotopy coincides with the two-dimensional  holonomy of a naturally associated 2-connection, as recently proved in \cite{AF2}. A similar result is mentioned in \cite{CFM2}. Later on, will also give a new insight to this non-trivial fact.

In this paper we propose a fully algebraic construction of the two dimensional holonomy operators.
The main reason for doing this is the following:
(local) 2-connections in a smooth manifold $M$ may arise simply as a pair $\big(\w\in \Omega^1(M,\lg),m \in \Omega^2(M,\ae)\big)$, with $\d(m)=\F_\w$, the curvature of $\w$, where $\mathfrak{X}=(\d\colon \ae \to \lg,\t)$ is a differential crossed module, with no obvious crossed module of Lie groups associated to it. The Lie algebras $\lg$ and $\ae$ may very well be infinite dimensional, and without any assigned or natural topology, which makes it impossible to apply Lie third theorem, in order to find a Lie crossed module $\cal{X}$ integrating $\mathfrak{X}$.

  Examples of such infinite dimensional Lie algebras and differential crossed module arise as quotients of free Lie algebras and of free differential crossed  modules.  An important case is the Lie algebra  $\ch_n$ (here $n$ is a positive integer), the Lie algebra of horizontal chord diagrams in $n$ strands, which is formally  generated by the symbols $r_{ab}$, where $1\leq a,b\leq n$, satisfying the infinitesimal braid group relations; \cite{Kassel,BN}:
\begin{align*}
& r_{ab} = r_{ba} ,  \\ & [r_{ab},r_{cd}] = 0 \; \mbox{ for } \{a,b\}\cap\{c,d\}=\emptyset \, , \\
& [r_{ab}+r_{ac},r_{bc}]  = 0 = [r_{ab},r_{ac}+r_{bc}]   \,\quad \quad \quad   \textrm{(called the ``four term relations'')}.&
\end{align*}

{The configuration space $\C(n)$} of $n$ (distinguishable) particles in the complex plane is, by definition, $\C(n)=\{(z_1,\dots,z_n)\in \C^n: z_i \neq z_j \textrm{ if } i \neq j\}.$  Consider the following $\ch_n$-valued connection form in $\C(n)$, called the universal Knizhnik-Zamolodchikov (KZ) connection, which is a flat connection, $\F_A=0$:
\begin{align}\label{defA}
A=\sum_{1 \leq a<b \leq n} \w_{ab}r_{ab}, \quad \quad \textrm{ where } \w_{ab}=\frac{d z_a-dz_b}{z_a-z_b}.
\end{align}
(We mention that the Kontsevich Integral \cite{Ko,BN} is essentially given by  the holonomy of the connection form $A$.)  Since there is no Lie group naturally assigned to $\ch_n$, we must adress the holonomy of $A$ in an algebraic fashion. We can do this is by introducing a formal parameter $h$, passing to the connection form  $A_h=\sum_{1 \leq a<b \leq n} h\w_{ab}r_{ab}.$ This connection is flat in a strong (graded) sense: $d A_h+ \frac{1}{2}(A_h \wedge A_h)=0$, since $dA=0$ and $A\wedge A=0$.  Then given a path $\gamma \colon [0,1] \to \C(n)$, the holonomy of $A_h$, over $\g$, it is given by the solution of the differential equation $\frac{d}{dt} a(t)=hA(\g'(t))\, a(t)$, $a(0)=1$, in the space $\U(\ch_n))[[h]]$ of formal power series over the universal enveloping algebra $\U(\ch_n)$ of $\ch_n$.   This has a  unique solution, which can be expressed by using Chen's iterated integrals; \cite{Chen} (we fully review the construction here). The map  $t \mapsto a(t)$ takes values in the group of group like 
elements of $\U(\ch_n))[[h]]$, namely $\D(a(t))=a(t) \tn a(t)$, where $\D$ is the coproduct map in the usual Hopf algebra 
structure given to $\U(\ch_n)$, extended to $\U(\ch_n)[[h]]$. 

Categorifications, in the realm of differential crossed modules, of the Lie algebra of horizontal chord diagrams were addressed in \cite{CFM1,CFM2}. {The calculations in \cite{CFM3} indicate that  the differential crossed module $2\mathfrak{ch}_n=(\d \colon 2{\rm ch}_n \to {\rm ch}_n^+,\t) $ associated to the Lie 2-algebra of horizontal {categorified} chord diagrams in $n$ strands should be defined in the following way. Firstly, ${\rm ch}_n^+$ is the quotient of the free Lie algebra on the symbols  $r_{ab}$ ($1 \leq a < b \leq n$) by the relations: 
$$ r_{ab} = r_{ba} ,  \quad \quad [r_{ab},r_{cd}] = 0 \; \mbox{ for } \{a,b\}\cap\{c,d\}=\emptyset  .$$
(Note that the 4-term relations are not imposed.) In order to construct  $2{\rm ch}_n$, consider the  free differential crossed module $(\d \colon 2{\rm ch}_n^+ \to {\rm ch}_n^+,\t)$  on the symbols $P_{abc}$ and $Q_{abc}$, with $1 \leq a < b < c \leq n$, such that:}
\begin{align*}
\d(P_{abc})&=[r_{bc},r_{ab}+r_{ac}],  && 
\d(Q_{abc})=[r_{ab},r_{ac}+r_{bc}]. 
\end{align*}
And then we form the quotient $2{\rm ch}_n$ of $2{\rm ch}_n^+$ by the ${\rm ch}_n^+$ module generated by the relations below: 
\begin{align*}
r_{ab} \t P_{ijk}&=0=r_{ab} \t Q_{ijk}, \textrm{ if } \{a,b\} \cap \{i,j,k\}=\emptyset\\
 (r_{ad}+r_{bd}+r_{cd}) \t { P_{abc}} -
 (r_{ab} + r_{ac})\t {Q_{bcd}}  + r_{bc} \t ({Q_{abd}}+ {Q_{acd}})&=0 \,, \\ 
 \left( r_{ab}+r_{ac}+r_{ad}\right) \t {P_{bcd}} + r_{cd} \t \big({P_{abc}} + {P_{abd}} \big) - \left( r_{bc}+r_{bd}\right) \t {P_{acd}} &=0\, ,\\ 
 (r_{ad} + r_{bd}+ r_{cd}) \t {Q_{abc}} + r_{ab} \t \big( {Q_{acd}} + {Q_{bcd}} \big)  - (r_{ac} + r_{bc}) \t {Q_{abd}} &= 0 \, ,\\ 
 (r_{ab} + r_{ac} + r_{ad}) \t {Q_{bcd}} + r_{bc} \t \big( {P_{abd}} + {P_{acd}} \big) - (r_{bd} + r_{cd}) \t {P_{abc}} &= 0 \, ,\\
  r_{ab} \t \big ( {P_{acd}}+{P_{bcd}}\big) - r_{cd} \t \big({Q_{abc}}+{Q_{abd}}\big) &= 0\,,\\
 r_{ac} \t \big ( {P_{abd}}-{P_{bcd}}-{Q_{bcd}}\big) + r_{bd} \t \big({Q_{abc}}+{P_{abc}}-{Q_{acd}}\big) & = 0\,.
 \end{align*}
 Note that the left-hand-side of each of these relations maps to zero in ${\rm ch}_n^+$, so that the quotient is well defined (for a proof see \cite{CFM1}.) Therefore we have a crossed module $2\mathfrak{ch}_n=(\d \colon 2{\rm ch}_n \to {\rm ch}_n^+,\t)$.

It is proven in \cite{CFM1,CFM3} that the pair $(A,B)$, where $A$ is as in \eqref{defA}, and $B$ is as below, is a flat 2-connection in the configuration space of $n$ particles in $\C$, called the universal KZ 2-connection:
 $$B= 2\sum_{1\leq a<b<c \leq n}\w_{bc}\wedge \w_{ca}\,P_{abc} -2\sum_{1\leq a<b<c\leq n}\w_{ca}\wedge \w_{ab}\,Q_{abc} \,.$$
 
 Our previous work on the categorification of the Lie algebra of horizontal chord diagrams, and the strongly related construction of the KZ-2-connection, was motivaded by the possible applications to the categorification of the Kontsevich Integral, in order to extended it, at least partially, to invariants of knotted surfaces in the 4-sphere. The project was also motivated by the issue of constructing categorified Drinfeld associators, presumably important for constructing explicit examples of braided monoidal bicategories. (In the uncategorified setting, Drinfeld associators can also be constructed from the holonomy of the KZ-connection).   Both of these applications require the definition of the two dimensional holonomy of the KZ-2-connection. Given that there is no crossed module of Lie groups associated to the differential crossed module $2\mathfrak{ch}_n$ {of categorified horizontal chord diagrams} a fully algebraic definition of this surface holonomy is very much needed. 
 
 \subsection{Plan of the paper}
  
 The algebraic and categorical foundations of this paper will be quite heavy, occupying the biggest part of it. We note that the word ``counit'' will be used both for designating categorical counits (in the context of adjoint functors) and for designating coalgebraic counits (in the context of Hopf algebras.) {Compositions in the text are intended from right to left.}

 \subsubsection{Algebraic setting: {one dimensional} holonomy}
 We will consider three types of associative algebras: unital algebras, where morphisms of unital algebras will, by imposition, always preserve identities, Hopf algebras, and bare algebras, which are associative algebras, not considering any other structure that might exist there (identity, coproduct, counit, etc). There is a pair of {adjoint} functors, where $U$ is the obvious underlying bare algebra (forgetful) functor:
 \begin{equation}\label{adbullet}
 \xymatrix{&\{\textrm{bare algebras} \}\ar@/^1pc/[r]^{()_\bullet} \ar@2[r]^{\xi_\bullet}  &\{\textrm{unital algebras} \} \ar@/^1pc/[l]^U_\bot} .
 \end{equation}
(Identity functors will not appear in the diagrams we draw. However they are always implicit, normally being targets of natural transformations, as above. Natural transformations will be denoted by double arrows $\implies$.)  
The left adoint $()_\bullet$ of the adjunction \eqref{adbullet} adjoins an identity to a bare algebra $A$, yielding $A_\bullet=\C\times A$, with product $(\l,a)(\l',a)=(\l \l',\l a'+\l'a +aa')$, where $a,a' \in A$ and $\l,\l' \in \C$. The counit $\ve_\bullet\colon ()_\bullet \circ U \implies \id$, of the adjunction, associates to a unital algebra $A$ the map ${\xi_\bullet}_A \colon A_\bullet \to A$, such that ${\xi_\bullet}_A(\l,a)=\l 1_A +a$, where $a \in A$ and $\l \in \C$.  (Here $1_A$ is the identity of $A$.)

We recall the well known adjunction between unital algebras and Lie algebras:
\begin{equation}\label{aduni}
\xymatrix{&\{\textrm{Lie algebras}\}\ar@/^1pc/[r]^\U  \ar@2[r]^{\xi} &\{\textrm{unital algebras}\} \ar@/^1pc/[l]^{{\rm Lie}}_\bot}.
\end{equation}
Here $\U$ is the universal enveloping algebra functor, and ${\rm Lie}$ associates to any unital algebra (or, more generally, bare algebra) $A$ {its} commutator Lie algebra: $[a,b]=ab-ba$, for each $a,b \in A$. Denote by $\xi\colon \U \circ \Lie \implies \id$ the counit of this adjunction, defining, for each unital algebra $A$, an algebra map $\xi_A\colon \U(\Lie(A)) \to A$. Recall that $\U(\lg)$, where $\lg$ is a Lie algebra, is a quotient of the tensor algebra $\T(\lg)=\C \oplus \bigoplus_{n =1}^{+\infty} \lg^{\otimes n}$ by the biideal generated by $\{XY-YX-[X,Y], X,Y \in \lg\}$. It is well known that $\U(\lg)$ is a Hopf algebra, where every element of $\lg\subset \U(\lg)$ is primitive: $\D(X)=X \tn 1+1 \tn X$, for each $X \in \lg$.  The (Hopf algebra) counit satisfies $\e(X)=0, \forall X \in \lg$.  

The  adjunction \eqref{aduni}, between unital algebras and Lie algebras, factors through another adjunction, between the categories of bare algebras and of Lie algebras, namely:
\begin{equation}\label{baread}
\xymatrix{
&\{\textrm{Lie algebras}\}\ar@/^1pc/[r]^{{\U^0}} \ar@2[r]^{\xi^{\flat}} &\{\textrm{bare algebras} \} \ar@/^1pc/[l]^{{\rm Lie}}_\bot
}.
\end{equation}
Here $\U^0(\lg)=\ker(\e \colon \U(\lg) \to \C)$ is the augmentation ideal of $\U(\lg)$. Theforore ${\U^0}(\lg)$ does not have a unit, but, consequence of the Poincar\'e-Birkhoff-Witt Theorem, we obtain $\U(\lg)$ by formally adjoining a unit to $\U^0(\lg)$. Given a bare algebra $A$, we  have an obvious map (the flat counit) $\xi^\flat_A \colon {\U^0}(\Lie(A)) \to A$.  The top and the bottom of the following diagram commute, and this goes over to (adjunction) units and counits:
  \begin{equation}
  \xymatrix{&\{\textrm{Lie algebras}\}\ar@/^2pc/[rr]!<-3ex,2ex>^\U \ar@/^1pc/[r]_{{\U^0}}
  &\{\textrm{bare algebras} \} \ar@/^1pc/[l]^{{\rm Lie}}_\bot \ar@/^1pc/[r]_{()_\bullet}  &\{\textrm{unital algebras} \} \ar@/^1pc/[l]^U_\bot   
  \ar@/^3pc/[ll]!<-3ex,1ex>^{\Lie}_{\bot} } 
  \end{equation}

   If $A$ is a bare algebra (without unit) there cannot exist a surjective algebra map $\U(\Lie(A)) \to A$. However, clearly, there is a natural transformation $\xi_\sharp$ (the sharp counit) between the functors $\U \circ \Lie$ and $()_\bullet$:   
   \begin{equation}
   \xymatrix{& &\ar@2[dr]!<-10ex,-2ex>^{\xi_\sharp} \\ &\{\textrm{Lie algebras}\}\ar@/^2.5pc/[rr]!<-3ex,2ex>^\U     &\{\textrm{bare algebras} \} \ar@/^1pc/[l]^{{\rm Lie}} \ar@/_1pc/[r]_{()_\bullet}  &\{\textrm{unital algebras} \}  }         
   \end{equation}
  For a bare algebra $A$,  $(\xi_\sharp)_A \colon \U(\Lie(A)) \to A_\bullet$ on generators sends $a \in A$ to $(0,a)$ and $1_{\U(\Lie(A))}$ to $(1,0)$.

   There is a functor ${\rm Prim}$ from the category of Hopf algebras to the category of Lie algebras, sending a Hopf algebra $H$ to its Lie algebra ${\rm Prim}(H)$  of primitive elements. We  have a functor $()^*_{\rm gl}$ from the category of Hopf algebras to the category of groups, associating to each Hopf algebra $H$ the group $H_{{\rm gl}}^*\subset H[[h]]$ of group like invertible elements of $H[[h]]$ (the algebra of formal power series over $H$), such that the first term {of the power series} is $1_H$, the identity of $H$. (The holonomy of a $\lg$ valued connection form in $M$  essentially takes values in $\U(\lg)^*_{\rm gl}$.) There is another functor $()^*$, from the category of unital algebras to the category of groups, associating to a unital algebra $A$ the group $A^*$ of formal power series in $A[[h]]$, {starting in} $1_A$.
   
   In the uncategorified setting, the diagram of all categories, functors and  natural transformations which we will use is below; \eqref{dia}. Note that $U$ always denotes 
the obvious forgetful functor.

   \begin{equation}\label{dia}
   \xymatrix{&\{\textrm{Lie algebras}\} \ar@2[r]^{\xi^{\flat}}  \ar@/^2pc/[rr]!<-3ex,2ex>^\U\ar@<-1ex>[dr]_\U \ar@/^1pc/[r]^>>>>{{\U^0}}  &\{\textrm{bare algebras} \}\ar@2[r]^{\xi_\bullet} \ar@/^1pc/[l]^<<<{{\rm Lie}}_\bot \ar@/^1pc/[r]^>>>>>>>{()_\bullet}  &\{\textrm{unital algebras} \}\ar@/_4pc/[ll]_{\rm Lie}^{\top} \ar@2[d]^{\iota}\ar[dr]^{()^*} \ar@/^1pc/[l]^U_\bot\\ & &\{\textrm{Hopf Algebras}\}\ar[rr]_{()^*_{\rm gl}} \ar@<-1ex>[ru]_U \ar@/^2pc/[lu]^<<<{\rm Prim} && \{{\rm Groups}\}
   }
   \end{equation}
   In \eqref{dia}, sub-diagrams where all arrows are left pointing commute. The same holds for right-pointing diagrams, except for the triangle in the bottom right.   We have  a natural transformation $\iota: ()^* \circ U \implies ()^*_{\rm gl}$, filling the triangle below. As before,  $\xi$ is the counit of the adjunction \eqref{aduni} between Lie algebras and unital algebras. 
\begin{equation}\label{poi}
   \xymatrix{&\{\textrm{Lie algebras}\}\ar@<-1ex>[dr]_\U  \ar@2[rr]^{\xi}  &&\{\textrm{unital algebras} \}   \ar@/_2pc/[ll]_{{\rm Lie}}^{\top} \ar@2[d]^{\iota} \ar[dr]^{()^*} \\ & &\{\textrm{Hopf Algebras}\}\ar[rr]_{()^*_{\rm gl}} \ar@<-1ex>[ru]^U 
    & &\{{\rm Groups}\} 
   }
   \end{equation}

   Here, if $H$ is an Hopf algebra, then $\iota_H$ is the group inclusion $H^*_{\rm gl} \to H^*$.  Composing the natural transformations $\xi$ and $\iota$ in \eqref{poi}, we have a natural transformation of functors $\{\textrm{unital algebras}\} \to \{ \textrm{groups} \}$:
  $${\rm proj}\colon ()^*_{\rm gl} \circ \U\circ {\rm Lie} \implies ()^*.$$

Explicitly, if $A$ is a unital algebra, one has the counit $\xi_A\colon \U(\Lie(A)) \to A$, extending to an algebra map  $\xi_A\colon \U(\Lie(A))[[h]] \to A[[h]]$. This clearly sends $\U(\Lie(A))^*_{\rm gl}$ to $A^*$.

  If we have a connection form $\w\in \Omega^1(M,A)$, where $M$ is a manifold and $A$ a unital  algebra, then we can address the holonomy of $\w$, either directly over $A[[h]]$, taking values in $A^*$, or over $\U(\Lie(A))[[h]]$, taking values in $\U(\Lie(A))^*_{\rm gl}$. In both cases, given a path $\g\colon [0,1] \to M$ we solve the differential equation $\frac{d}{d t} P_\w (\gamma, [t, t_0 ])=h\w(\g'(t))\,  P_\w (\gamma, [t, t_0 ])$, with $P_\w(\g,[t_0,t_0])=1$, however firstly over $\U(\Lie(A))[[h]]$, and in the latter case over $A[[h]]$. In the first case the holonomy takes values in $\U(\Lie(A))^*_{\rm gl}$, whereas in the second case it takes values over $A^*$. These two holonomies are related by the map ${\rm proj}_A\colon \U(\Lie(A))^*_{\rm gl} \to A^*$. 
   \subsubsection{Algebraic and categorical setting: two-dimensional holonomy}

  In order to consider two dimensional holonomies, we first extend the previous diagrams \eqref{dia} of categories,  functors and natural transformations to the categorified (crossed module) setting. We will consider in total five categories of crossed modules, namely: crossed modules of Lie algebras  and of (Lie) groups (already mentioned), and, in addition, crossed modules of Hopf algebras, crossed modules of bare algebras and unital crossed modules of bare algebras. The patern of the construction is always the same. For example, a crossed module  $\H=(\d \colon I \to H, \rho)$, of Hopf algebras  \cite{MajidHA},  is given by a Hopf algebra map $\d\colon I \to H$, and a linear map $\rho\colon H \tn I \to I$, making $I$ into an $H$-module algebra and coalgebra (and $\rho$ should satisfy the obvious Hopf analogues of the first and second Peiffer relations.) 
  
  A crossed module $\A=(\d\colon A \to B, {\tra,\tla})$ of bare algebras is given by a bare algebra map $\d\colon A \to B$, together with two actions $\tra$ and $\tla$, on the left and on the right, of $B$ on $A$, by (left or right) multipliers. The left and right actions  are to be compatible, each furthermore satisfying the first and second Peiffer relations. A crossed module of bare  algebras is said to be unital if the bottom algebra $B$ is unital, with an identity $1_B$ satisfying $1_B \tra a = a= a \tla 1_B$, for each $a \in A$. If (possibly) only the first Peiffer relation is satisfied by $\tra$ and $\tla$ then $\hat{\A}=(\d\colon A \to B, {\tra,\tla})$ will be called a pre-crossed module of bare algebras, and similarly for Hopf algebras. As in the case of groups and Lie algebras, a pre-crossed module  $\hat{\A}=(\d\colon A \to B, {\tra,\tla})$, of bare algebras, can be converted (by a reflection functor ${\cal R}$) into a crossed module of bare algebras ${\cal R}(\hat{A})=(\d\colon \underline{A} 
\to B, {\tra,\tla})$. We obtain $\underline{A}$ from $A$ 
by dividing out the ideal 
generated by the second Peiffer identities. 
  
  We have the following diagram of categories and functors, where the acronyms ``C.M.'' and ``P.C.M''mean crossed modules and pre-crossed modules, extending some of the instances of the previous discussion to the ``crossed module'' context,  with some new categories and functors added; c.f. \eqref{dia}. Full explanation is below.
  \begin{equation}\label{BIG}
  \hskip-2.8cm\xymatrix@R=30pt@C=50pt{&\{\textrm{Lie algebras C.M.}\}
  \ar@/_1pc/[d]_\U 
  &\{\textrm{bare algebras C.M.} \}\ar[l]_{{\rm Lie}}
  &\{\textrm{unital bare algebras C.M.} \} 
  \ar@/_2pc/[ll]!<-3ex,2ex>_{\rm Lie}^{\top}\ar[ddl]^{()^*_\bullet}
  \ar[l]_U
   \ar[dd]^{{\cal C}^+() }
  \\ 
  &\{\textrm{Hopf Algebras C.M.}\} \ar[rru]_{{\cal BA}}  \ar[r]_{\widehat{\cal BA}}\ar[dr]_{()^*_{\rm gl}}  \ar@/_1pc/[u]_{\rm Prim}  &\{\textrm{unital bare algebras P.C.M.} \}\ar[ru]_{\cal R} &\\  
    & &\{\textrm{Groups C.M.}\} \ar[r]^{{\cal C}^\times()}
& \{\textrm{2-groupoids}\}
    }
    \end{equation}
    (Left pointing arrows commute. Right pointing arrows commute.) 
  Let us briefly explain the diagram (leaving full details for the paper itself). {We have a forgetful functor $U$} from the category of unital crossed module of bare algebras to the category of crossed modules of bare algebras. A crossed module of bare algebras $\A=(\d\colon A \to B, {\tra,\tla})$ gives rise to a crossed module of Lie algebras, simply as: 
    $\Lie(\A)=(\d\colon \Lie(A) \to \Lie(B), \t)$, where $b \t a=b \tra a-a\tla b$, for $b \in B$ and $a \in A$. Suppose that $\A$ is unital, thus that $B$ is a unital algebra. Let $\Ad^* \subsetneq {\rm Units}(\Ad\h)$, where ``${\rm Units}$'' denotes the group of invertible elements, be the group of formal power series in $\Ad=\C \times A$, whose first term is $1_\C=(1_\C,0)$, and such that all other terms are $h^n$ ($n \ge 1$) multiples of elements of $A=\{0\} \times A$.  We have a crossed module of groups: $(\A)^*_\bullet=(\d\colon \A^*_\bullet \to B^*, \t)$, where, as before, $B^*\subset {\rm Units}(B[[h]])$ is the group of formal power series starting in $1_B$.

    A crossed module of Lie algebras $\mathfrak{X}=(\d\colon \ae \to \lg, \t)$ gives rise to a crossed module of Hopf algebras $\U(\mathfrak{X})=(\d\colon \U(\ae ) \to \U(\lg), \rho)$, \cite{MajidHA}, in a natural way. Passing to the Lie algebras of primitive elements, if $\H=(\d \colon I \to H, \rho)$ is a Hopf algebra crossed module, then $\Prim(\H)=(\d \colon \Prim(I) \to \Prim(H), \t)$ is a differential crossed module. Considering the groups of group-like elements in $H[[h]]$ and $I[[h]]$, whose first term of the power series is the identity, yields a crossed module of groups $(\H)^*_{\rm gl}=(\d \colon I^*_{\rm gl} \to H^*_{\rm gl}, \t)$.

    There is also a functor $\widehat{\cal BA}$ sending a crossed module  $\H=(\d \colon I \to H, \rho)$, of Hopf algebras, to a pre-crossed module of bare algebras. Explicitly: $\widehat{{\cal \BA}} \big((\d \colon I \to H, \rho)\big)= \left (t\colon I^0 \tn_{\rho} H  \to H, \tra, \tla \right), $ where $I^0$ is the augmentation ideal (the kernel of the  counit) of $I$, and  $I^0 \tn_{\rho} H $ is the smash product as an algebra and a product coalgebra. Also $t(a \tn x)=\d(a)x$, where $a \in I$ and $x \in H$.  The natural inclusion $H \to I \tn_\rho H$ turns $I^0 \tn_{\rho} H$ into a left and a right module over $H$. Composing $\widehat{{\cal BA}}$ with ${\cal R}$ gives a functor ${\cal BA}$, associating a crossed module of bare algebras to a crossed module of Hopf algebras. We have an adjunction between the categories of unital crossed modules of bare algebras and of crossed modules of Lie algebras, where ${\cal BA} \circ \U= {\cal R} \circ \widehat{ {\cal BA}} \circ \U$ is the left adjoint: 
    \begin{equation}\label{AAA}
    \hskip-2.8cm\xymatrix@R=30pt@C=50pt{&\{\textrm{Lie algebras C.M.}\}
  \ar@/_1pc/[d]_\U 
  & &\{\textrm{unital bare algebras C.M.}\}
  \ar@/_2pc/[ll]!<-3ex,2ex>_{\rm Lie}^{\top}
  \\ 
  &\{\textrm{Hopf Algebras C.M.}\} \ar[rru]_{{\cal BA}}  \ar[r]_{\widehat{\cal BA}} &\{\textrm{unital bare algebras P.C.M.} \}\ar[ru]_{\cal R}}
  \end{equation}
  
  Given a unital crossed module of bare algebras $\A=(\d\colon A \to B, {\tra,\tla})$, we have a map of pre-crossed modules of bare algebras:
  $\hat{K}_\A=(\kappa_\A, \xi_B)\colon \widehat{\cal BA}\big(\U(\Lie(\A))\big) \to \A$, where $\xi_B\colon \U(\Lie(B)) \to B$ is the usual counit, and  $\kappa_\A( a \tn  x)= \xi_A^\flat(a) \tla \xi_B(x)$, where $a \in \U^0(\Lie(A))$ and $x \in \U(\Lie(B))$. (Note that $A$ need not have a unit).  Given that the second Peiffer relations are satisfied by $\A$, $\hat{K}_\A$ descends to a map $K_\A \colon {\cal BA}(\Lie(\A)) \to \A$, of crossed modules of bare algebras, called the  bare counit. It defines (see below) a natural transformation of functors $K\colon {\cal BA} \circ \U \circ {\rm Lie} \implies \id$, the latter being the identity functor of the category of unital crossed modules of bare algebras. The natural transformation $K$ is the counit of the adjunction of diagram \eqref{AAA}:
  \begin{equation}\label{bareco}
  \xymatrix{
  &\{\textrm{Lie algebras C.M.}\}\ar@2[rr]!<1ex,2ex>_{{K}}   \ar@/_1pc/[d]_\U & 
    &\{\textrm{unital bare algebras C.M.} \}
  \ar@/_1.3pc/[ll]!<-3ex,2ex>_{\rm Lie}^{\top}
  \\ & \{\textrm{Hopf Algebras C.M.}\} \ar[rru]_{{\cal BA}}  && } .
  \end{equation}
   
  Given a unital crossed module of bare algebras $\A=(\d\colon A \to B, {\tra,\tla})$, the pair of algebra maps $\big(\xi^\sharp_A\colon \U(\Lie(A)) \to A_\bullet,\xi_B\colon \U(\Lie(B)) \to B\big)$ yields  a group crossed module map  ${\rm Proj}_\A\colon    \big(\U(\Lie(\A)\big)^*_{\rm gl} \to (\A)^*_\bullet $, called the  crossed module counit. This yields a natural transformation of functors, as below:
  \begin{equation}\label{cmcounit}
  \xymatrix{&\{\textrm{Lie algebras C.M.}\}   \ar@/_1pc/[d]_\U 
  & \ar@2[dr]!<-19ex,1ex>_{\rm Proj} 
  &\{\textrm{unital bare algebras C.M.} \} 
  \ar@/_1.3pc/[ll]!<-3ex,2ex>_{\rm Lie}\ar[ddl]^{()^*_\bullet}
  \\ & \{\textrm{Hopf Algebras C.M.}\} \ar[dr]_{()^*_{\rm gl}}  &&\\   &\
    & \{\textrm{Groups C.M.}\}
    } .
    \end{equation}
  The crossed module counit  ${\rm Proj}_\A\colon    \big(\U(\Lie(\A)\big)^*_{\rm gl} \to (\A)^*_\bullet $, is the composition of the natural transformation in \eqref{bareco} with the  natural transformation ${\rm Inc}=(J,\iota) \colon ()^*_\bullet \circ {\cal BA} \implies ()^*_{{\rm gl}}$, the crossed module inclusion; \eqref{IOTA}.
  \begin{equation}\label{IOTA}
  \xymatrix{
  & & 
  &\{\textrm{unital bare algebras C.M.} \} \ar@{<=}[dll]!<-5ex,-8ex>^>>>>>>>>>>>>>>{\rm Inc}
 \ar[ddl]^{()^*_\bullet}
  \\ & \{\textrm{Hopf Algebras C.M.}\}\ar[rru]^{{\cal BA}} \ar[dr]_{()^*_{\rm gl}}  &    &\\   &
    & \{\textrm{Groups C.M.}\} &
    } 
    \end{equation}
  Given a crossed module    $\H=(\d \colon I \to H, \rho)$, then $\iota_\H\colon H^*_{\rm gl} \to H^*$ is simply group inclusion. But not for $J_\H$.

  In the bottom right corner of the diagram \eqref{BIG}, we have the usual functor $\Cc^\times$, mentioned in the beginning of the Introduction, from the category of group crossed modules to the category of 2-groupoids; \cite{BHS}. If $\X=(\d\colon E \to G, \t)$ is a crossed module of groups, all compositions in the 2-groupoid $\Cc^\times(\A)$ are induced by the group multiplications of $E$ and $G$. We can also go directly from the category of unital crossed modules of bare algebras to the category of 2-groupoids, through a functor ${{\cal C}^+()}$, which, given a unital crossed module of bare algebras $\A=(\d\colon A \to B, {\tra,\tla})$, uses the sum in $A$ to define the vertical composition of 2-morphisms.  What is far from obvious is that there exists also a natural transformation $\T$ connecting the functor ${{\cal C}^\times()}\circ ()^*_\bullet$ to the functor  ${{\cal C}^+()}$; see the diagram below. 
  \begin{equation}\label{defT}
  \xymatrix{&
  &\{\textrm{unital bare algebras C.M.} \}  
  \ar[dl]_{()^*_\bullet}
   \ar[d]<2.5ex>^{{\cal C}^+() }
   \\
    & \{\textrm{Groups C.M.}\}\ar@2[ur]!<2ex,-3ex>|\T \ar[r]_{{\cal C}^\times()}
& \{\textrm{2-groupoids}\}
    } 
    \end{equation}
    
    We note that a pre-crossed module $\hat{\A}$ of bare algebras also defines a sesquigroupoid  ${\hat{\cal C}^+(\hat{\A})}$, the latter being a 2-groupoid if, and only if, $\hat{\A}$ is a crossed module. A sesquigroupoid  \cite{St} is defined exactly like a 2-groupoid, except that we do not impose the interchange law to hold, and therefore the horizontal composition of 2-morphisms in  ${\hat{\cal C}^+(\hat{A})}$ is not necessarily well defined. 
   
  \subsubsection{Main differential geometric results}
  Let $M$ be a differential manifold. {Let $\P(M)=(\P^2(M),\P^1(M),M)$ be disjoint union of the sets of (piecewise smooth) 2-paths $\Gamma \colon [0,1]^2 \to M$, 1-paths $\g\colon [0,1] \to M$, as well as the underlying set of $M$. There are obvious vertical and horizontal   composition (partial) operation in $\P^2(M)$, as well as vertical and horizontal reverses, composition of paths in $\P^1(M)$, and reverses, and also two whiskering (left and right) of 2-paths by 1-paths. We additionally consider several, however obvious,  boundary and identity maps.} 
  The main results of this paper as far as two-dimensional holonomy is concerned are:
  \begin{itemize}
   \item[]{\bf The Exact Holonomy  ${\cal F}^\times_{(\bw,\bm)}$:} Let $\H=(\d \colon I \to H, \rho)$ be a crossed module of Hopf algebras. Let $(\w,m_1,m_2)$ be a fully primitive Hopf 2-connection over $\H$. This means that $\w\in \Omega^1(M,H)$ is a 1-form in $M$, with values in $H$, and $m_1,m_2 \in\Omega^2(M,I)$ are  2-forms in $M$, with values in $I$, compatible in the sense that  $\d(m_1)=d\w$ and $\d(m_2)=-\frac{1}{2} \w \wedge \w$. Considering the graded forms $\bw=h\w$ and $\bm=hm_1+h^2 m_2$, the compatibility condition is equivalent to imposing that $\d(\bm)=\F_{\bw}=h d \w -h^2\frac{1}{2} \w \wedge \w$. 
   The prefix ``fully primitive'' means that both $\w$, $m_1$ and $m_2$ take values in the primitive spaces of $H$ and $I$. We have a map ${\cal F}^\times_{(\bw,\bm)}\colon \P(M) \to \Cc^\times(\H^*_{\rm gl})$, the latter being a 2-groupoid, preserving all compositions,  boundaries and reverses, in $\P(M)$ and in  $\Cc^\times(\H^*_{\rm gl})$. We call ${\cal F}^\times_{(\bw,\bm)}$ the ``Exact Holonomy''. 
   
   \item[]{\bf The Bare Holonomy   ${\cal F}^+_{(\bw,\bm)}$:} Let $\hat{\A}=(\d\colon A \to B, {\tra,\tla})$, be a pre-crossed module of bare algebras, where ${\tra,\tla}$ are left and right actions of $B$ on $A$, by multipliers.  Let $(\w,m_1,m_2)$ be a bare 2-connection over $\A$. This means that $\w \in \Omega^1(M,B)$ and $m_1,m_2 \in \Omega^2(M,A)$ are such that, if we put $\bw=h\w$ and $\bm=hm_1+h^2 m_2$, it holds that $\d(\bm)=\F_{\bw}=h d \w - h^2\frac{1}{2}\w \wedge \w$. We have a map ${\cal F}^+_{(\bw,\bm)}\colon \P(M) \to \hat{\Cc}^+(\hat{\A})$, preserving all boundaries, composition of 1-morphisms, vertical reverses and composition of 2-morphisms, but not necessarily the horizontal composition of 2-morphisms (which is not even well defined in $\hat{\Cc}^+(\hat{\A})$, unless   $\hat{\A}$ is a crossed module, of bare algebras.)   If $\hat{\A}$ is a crossed module, thus $\hat{\Cc}^+(\hat{\A})={\Cc}^+(\hat{\A})$ is a 2-groupoid, then the map ${\cal F}^+_{(\bw,\bm)}\colon \P(M) \to {\Cc}^+(\hat{\A})$ 
also preserves horizontal 
composition (and reverses)  of 2-morphisms, in $\P(M)$ and in ${\Cc}^+(\hat{\A})$.
   
   \item[]{\bf The Fuzzy Holonomy ${{\Fz}^+_{(\bw,\bm)}}$ and the Blur Holonomy ${\widehat{\Fz}^+_{(\bw,\bm)}}$} Let $\H=(\d \colon I \to H, \rho)$ be a crossed module of Hopf algebras. Let $(\w,m_1,m_2)$ be a fully primitive Hopf 2-connection over $\H$. We can see $(\w,m_1,m_2)$ as being a bare 2-connection with values in $\widehat{{\cal BA}}(\H)=\left (t\colon I^0 \tn_{\rho} H  \to H, \tra, \tla \right)$, a pre-crossed module of bare algebras, or, projecting into  ${\cal BA}(\H)=\left (t\colon \underline{I^0 \tn_{\rho} H }  \to H, \tra, \tla \right)$, as a bare 2-connection with values in the crossed module ${{\cal BA}}(\H)$. The two dimensional holonomy of the former is called the blur holonomy, and is denoted by  ${\widehat{\Fz}^+_{(\bw,\bm)}}\colon \P(M) \to \hat{\Cc}^+(\widehat{\cal BA}(\H))$. The two dimensional holonomy of the latter is called the fuzzy holonomy and is denoted by:  ${\Fz^+_{(\bw,\bm)}}\colon \P(M) \to \Cc^+({\cal BA}(\H))$. By construction, and the item above, the 
fuzzy 
holonomy preserves all compositions and reverses. The blur holonomy does not, necessarily, {preserve} horizontal compositions (and reverses) of 2-morphisms. 
   
   \item[]{\bf Compatibility between formulations of two dimensional holonomy-I} This is a non-trivial relation between the bare and exact formulations of two dimensional holonomy (a relation between the exact and the seemingly much finer blur holonomy remains elusive). Let $\H=(\d \colon I \to H, \rho)$ be a crossed module of Hopf algebras. Let $(\w,m_1,m_2)$ be a fully primitive Hopf 2-connection over $\H$. The fuzzy and exact holonomies of $(\w,m_1,m_2)$ are compatible, in the sense that the following diagram commutes (recall the natural transformations ${\rm Inc}$ and $\T$ of diagrams \eqref{IOTA} and \eqref{defT}):
   \begin{equation}\label{comp1}
  \hskip-2cm \xymatrix{
   &&\Cc^\times\Big(\big (\H\big)^*_{\rm gl}\Big) 
   \ar[d]^{\Cc^\times\left ({\rm Inc}_{\H}\right)}&&&\P(M)\ar[lll]_{ {\cal F}^\times_{(\bw,\bm)}  }   \ar[d]^{ \Fz^+_{(\bw,\bm)}  }   \\
   && \Cc^\times\Big({\cal BA}\left(\H\right)^*_\bullet \big) \ar[rrr]_{ {\cal T}_{ {\cal BA}\left (\H\right)}}  & & &\Cc^+\Big({\cal BA}\left(\H\right) \Big)
   } 
   \end{equation}
   
   \item[]{\bf Compatibility between formulations of two dimensional holonomy-II} Let $\A=(\d\colon A \to B, {\tra,\tla})$, be a crossed module of bare algebras.  Let $(\w,m_1,m_2)$ be a bare 2-connection over $\A$.  We can also see $(\w,m_1,m_2)$ as being a fully  primitive Hopf 2-connection over $\U(\Lie(\A))$, and consequently also as a bare 2-connection over 
   ${\cal BA}\big(\U(\Lie(\A))\big)$. Recall \eqref{BIG},\eqref{bareco}, \eqref{cmcounit} and \eqref{IOTA}.    The  diagram commutes:
  \begin{equation}\label{Comp2}
  \hskip-2cm \xymatrix{
   &&\Cc^\times\Big(\U\big(\Lie(\A)\big)^*_{\rm gl}\Big) \ar@/_1pc/[ddl]_{{\Cc^\times ({\rm Proj}_\A)}}
   \ar[d]^{\Cc^\times\left ({\rm Inc}_{\U(\Lie(\A)}\right)}&&&\P(M)\ar[lll]_{ {\cal F}^\times_{(\bw,\bm)}  }   \ar[d]^{ \Fz^+_{(\bw,\bm)}  } \ar[drr]^{ {\cal F}^+_{(\bw,\bm)}  }    \\
   && \Cc^\times\Big({\cal BA}\left(\U\big(\Lie(\A)\big)\right)^*_\bullet \big) \ar[rrr]_{ {\cal T}_{ {\cal BA}\left (\U\big(\Lie(\A)\big)\right)}}\ar[dl]^{{\cal C}^\times\Big(\big (K_\A\big)^*_{\bullet}\Big) }]  & & &\Cc^+\Big({\cal BA}\left(\U\big(\Lie(\A)\big)\right) \Big)\ar[rr]_<<<<<<<<<{\Cc^+({ K}_\A)}&& \Cc^+(\A)  \\
   & \Cc^\times(\A^*_\bullet)\ar@/_2.5pc/[rrrrrru]^{\Tc_\A}  &&&& &
   } .
   \end{equation}
(This essentially follows from the previous diagram \eqref{comp1}, either by definition, functorially or naturally.)

If we unpack what {the} commutativity of the outer polygon of \eqref{Comp2} means, in the  particular case of crossed modules derived from chain-complexes of vector spaces, we obtain an equality between two versions of two dimensional holonomy (of 2-connections and of representations up to homotopy), previously independently  noticed in \cite{CFM2,AF2}, being fully proved in the latter reference. 
    \end{itemize}
    
\section{Algebraic preliminaries}
All algebras will by default be over $\C$, sometimes over $\C\h$, the latter always being explicitly mentioned. 
\subsection{Crossed modules of groups and of Lie algebras}
\begin{Definition}[(Pre-) Crossed module of (Lie) groups]
 A modern treatment of group crossed modules appears in \cite{Baues,BL,BHS,Br}. A crossed module (of groups) ${\X= ( \d\colon E \to  G,\t)}$ is given by a group morphism $\d\colon E \to G$, together with a left action $\t$ of $G$ on $E$ by automorphisms, {such that:} 
\begin{enumerate}
 \item $\d(g \t e)=g \d(e)g^{-1}; \fo g \in G \an  e \in E,$ (First Peiffer relation).
  \item $\d(e) \t e'=ee'e^{-1};\fo  e,e' \in E.$ (Second Peiffer relation).
\end{enumerate} If both $G$ and $E$ are Lie groups, $\d\colon E \to G$ is a smooth morphism,  and the left action of $G$ on $E$ is smooth, then $\X$ will be called a Lie crossed module. 
If $\hat{\X}$ satisfies all relations for it to be a group crossed module, except possibly for the second Peiffer relation, then $\hat{\X}=( \d\colon E \to  G,\t)$ will be called a pre-crossed module.
\end{Definition}

 A morphism $(\psi,\phi)\colon \X \to \X'$ between the crossed modules ${\X= ( \d\colon H \to  G,\t)}$  and {$\X'=(\d'\colon H' \to G',\t')$} is given by a pair of group morphisms $\left ( \psi\colon H \to H',\phi\colon G \to G' \right)$, such that the diagram:
 $$ \begin{CD}
   H @>\d>> G \\
  @V \psi VV   @VV \phi V \\
   H' @>\d'>> G' \\
    \end{CD}
 $$
 is commutative. In addition, we must have $\psi(g \t h)=\phi(g) \t' \psi(h)$, for each $ h \in H$ and each $ g \in G$.

Given a Lie crossed module ${\X= ( \d\colon H \to  G,\t)}$, then the induced Lie algebra map $\d\colon \ae \to \lg$, together with the derived action of $\lg$ on $\ae$ (also denoted by $\t$) is a differential crossed module, in the sense below:

\begin{Definition}[Differential (pre-) crossed module]  A differential crossed module  ${\mathfrak{X}=(\d \colon \ae \to  \lg,\t )}$, see \cite{BL,BC,W}, also called crossed module of Lie algebras, is given by a Lie algebra map $\d\colon \ae \to \lg$, together with a left action (by derivations, see below) of $\lg$ on the underlying vector space of  $\ae$, {such that}:
\begin{enumerate}
 \item For any $X \in \lg$ the map $\zeta \in \ae \mapsto X \t \zeta \in \ae$ is a derivation of $\ae$, in other words:
$$X \t [\zeta,\nu]=[X \t \zeta,\nu]+[\zeta, X \t \nu]\,,\quad \forall\, X \in \lg\,, \forall\, \zeta ,\nu\in \ae \,.$$
\item The map $\lg \to \Der(\ae)$, from $\lg$ into the derivation algebra of $\ae$, induced by the action of $\lg$ on $\ae$, is a Lie algebra morphism, in other words:
$$[X,Y] \t \zeta=X \t (Y \t \zeta)-Y \t(X \t \zeta)\,,\quad \forall\, X,Y \in \lg \,,\forall\, \zeta \in \ae \,,$$
\item $\d( X \t \zeta)= [X,\d(\zeta)]\,,\quad \forall X \in \lg\,, \,\forall\, \zeta\in \ae\,,\quad 
\textrm{(First Peiffer relation).}$
\item $ \d(\zeta) \t \nu=[\zeta,\nu]\,,\quad  \forall\, \xi,\nu\in \ae\,,\quad \textrm{(Second Peiffer relation).}$
\end{enumerate}
If ${\hat{\mathfrak{X}}=(\d \colon \ae \to  \lg,\t )}$ satisfies all conditions for it to be a differential crossed module, except possibly for the second Peiffer relation, then  ${\hat{\mathfrak{X}}=(\d \colon \ae \to  \lg,\t )}$ will be called a differential pre-crossed module. 
\end{Definition}
\noindent Given a differential crossed module ${\mathfrak{X}=(\d \colon \ae \to  \lg,\t )}$, where $\ae$ and $\lg$ are finite dimensional,  there exists a unique, up to isomorphism, crossed module of simply connected Lie groups ${\X=(\d \colon H \to  G,\t )}$ whose differential form is $\mathfrak{X}$.

\subsection{Crossed modules of bare algebras}
\subsubsection{Definition of crossed modules and pre-crossed modules of bare algebras}\label{dcmba}
\begin{Definition}[Bare algebra]
 A bare algebra $A$ (as opposed to a bi-algebra, a unital algebra, or a Hopf algebra) is an associative algebra (not necessarily abelian) over a field (in this paper always $\C$) with or without identity, and without taking into account any other structure which may possibly exist in $A$.
\end{Definition}
We propose the following {definition} of a crossed module of bare algebras. This is a working definition, completely bypassing its (unknown, and possibly not-existing) relation with simplicial algebras with Moore complex of length one. We note that crossed modules of associative algebras are defined in \cite{A,CIKL,DIKL}, and our {definition} is more restrictive. Notice that, in the case of commutative algebras (where left and right actions are the same thing) our {definition} of crossed modules coincides with the one in \cite{PA}.

We recall that if $A$ is a bare algebra, then a left multiplier of $A$ is a linear map $f\colon A \to A$ such that $f(ab)=f(a) b$, for each $a,b \in A$. A right multiplier is a linear map $g\colon A \to A$ such that $g(ab)=ag(b)$, for each $a,b \in A$. The composition of linear maps $A \to A$ gives the algebras of left and right multipliers of $A$. 

\begin{Definition}[Crossed modules of bare algebras]
 A crossed module $\A=(\d\colon A \to B, {\tra,\tla})$ of bare algebras,  is given by an algebra map $\d\colon A \to B$, together with two bilinear maps:
 \begin{align*}
& (a,b) \in A\times B  \mapsto a \tla b \in A,  &(b,a) \in B \times A \mapsto b \tra a \in A .
\end{align*}
The remaining conditions are:
\begin{itemize}
 \item  $\tla$ and $\tra$ are right and left actions of algebras on vector spaces, namely we have:
$$(b b') \tra a=b\tra( b' {\tra} a) \an a \tla( bb')=(a \tla b) \tla b',  \fo b,b' \in B \an a \in A.$$
\item $\tla$ and $\tra$ are (right and left) actions by (left and right multipliers) and therefore are compatible with the multiplication in $A$:
$$b \tra (aa')=(b \tra a)a' \an (aa') \tla b=a (a' \tla b), \fo a,a'\in A \an b \in B. $$
 \item The left and right actions commute:
$$b \tra (a \tla b')=(b \tra a) \tla b', \fo b,b' \in B \an a \in A. $$

 \item The left and right actions are fully interchangeable, in the sense that:
$$ (a \tla b)\, (b' \tra a') =a\,  \big ( (bb') \tra {a'}\big ) =\big (a \tla (bb') \big)\, a' , \fo b,b' \in B \an a,a' \in A. $$
\item The first Peiffer relations hold: $$\d(b  \tra a)=b\d(a) \an \d(a \tla b)=\d(a)b, \fo a \in A \an b \in B.$$

\item The second Peiffer relations hold: $$\d(a') \tra a =a'a \an a \tla \d(a')=aa', \fo a,a' \in A.$$
\end{itemize}
As it has been customary, if  $\hat\A=(\d\colon A \to B, {\tra,\tla})$ is such that all conditions for it to be a crossed module of bare algebras hold, except possibly for the second Peiffer relation, then  $\hat\A=(\d\colon A \to B, {\tra,\tla})$ will be called a pre-crossed module, of bare algebras. 
\end{Definition}
\begin{Definition}[Unital crossed module of bare algebras]
A crossed (or pre-crossed) module $\A=(\d\colon A \to B, {\tra,\tla})$ of bare algebras is said to be unital, if  $B$ has an identity $1_B$,  satisfying $1_B \tra a=a=a \tla 1_B , \forall a \in A$.
\end{Definition}
\subsubsection{From pre-crossed modules of bare algebras to crossed modules of bare algebras}\label{refl}

 A useful identity holding in any crossed module of bare algebras  $\A=(\d\colon A \to B, {\tra,\tla})$ is the following:
\begin{equation}\label{excl}
\d(a)\tra a'=aa'=a \tla \d(a'), \fo a,a' \in A.
\end{equation}
This fails to be true, in general, in a pre-crossed module. Given  a pre-crossed module  $\hat\A=(\d\colon \hat{A} \to B, {\tra,\tla})$, of bare algebras, consider the following bilinear map in $A$, vanishing in the crossed module case:
\begin{equation}\label{ppba}
 (a,a') \in A \times A \mapsto \{a,a'\}\doteq \d(a) \tra a' -a \tla \d(a') \in A.
\end{equation}
Let $L$ be the biideal of $A$ generated by the image of the linear map $a \tn a' \in A \otimes A   \mapsto \{a,a' \}.$ The pre-crossed module axioms imply that $L$ stays invariant under the actions $\tra$ and $\tla$ of $B$ on $A$. 
Put $\underline{A}=A/L$. Then the actions $\tra$ and $\tla$ descend to actions $\tra$ and $\tla$ of $B$ on $\underline{A}$, by multipliers. By using the first Peiffer relation $\d\colon A \to B$ descends to a map $\d\colon \underline{A} \to B$. Clearly $\big(\d\colon \underline{A} \to B, \tra,\tla)\doteq {\cal R}\big ( (\d\colon {A} \to B, \tra,\tla)\big) $ is a crossed module of bare algebras. Moreover we  have a functor ${\cal R}$, from the category of pre-crossed modules of bare algebras, to the category of crossed modules of bare algebras, making the latter a reflexive subcategory.  
\subsubsection{From crossed modules of bare algebras to differential crossed modules}\label{bdga}
{The following generalizes the construction of a Lie algebra $\Lie(A)$ as being the commutator algebra of a bare algebra $A$, namely $[a,b]=ab-ba, \forall a,b \in A$. It follows from simple calculations. The fact that the left and right actions are fully interchangeable has a primary role in the proof that the action $\t$, below, is by derivations. That $\t$ is itself an action uses solely the fact that $\tra$ and $\tla$ commute.}  
\begin{Lemma}
Let $\A=(\d\colon A \to B, {\tra,\tla})$ be a crossed module of bare algebras. Then, looking at the associated Lie algebras $\Lie(A)$ and $\Lie(B)$, and the induced Lie algebra map {$\d \colon \Lie(A) \to \Lie(B)$}, the bilinear map:
$$(b,a) \in B \otimes A {\mapsto}b \t a\doteq b \tra a -a\tla b \in A,   $$
is a Lie algebra action of $B$ on $A$ by derivations, and, moreover, $\Lie(\A)=(\d \colon \Lie(A) \to \Lie(B), \t)$ is a crossed module of Lie algebras.
\end{Lemma} 

\subsubsection{From unital crossed modules of bare algebras to group crossed modules}\label{bag}
{There exists a functor from the category of unital crossed modules of bare algebras  $\A=(\d\colon A \to B, {\tra,\tla})$, to the category of crossed modules of groups. In the finite dimensional case, this is related to the composition of the functors, sending a crossed module of bare algebras to a differential crossed module, and sending a differential crossed module to a crossed module of simply connected Lie groups. In order to work with explicit models, and in the infinite dimensional case, we will do all constructions mostly over  $\Ch$ (except for the initial construction, only done for motivating purposes), making this relation  appear less clear} 

\begin{Definition}[From $A$ to $A_\bullet$: adjoining an identity]
Let $A$ be bare algebra. We denote by $A_\bullet$ the algebra derived from $A$ by adjoining an identity, which is different from an eventually existing identity of $A$. Explicitly $A_\bullet = \C \oplus A$, with product $(\l,a)(\l',a')=(\l \l', \l a'+\l' a + aa')$, where $\l,\l'\in \C$ and $a,a' \in A$. We note that adjoining an identity is functorial, being  left adjoint to the forgetful functor that assigns a bare algebra to an algebra with identity. (In the category of unital algebras morphisms are assumed to preserve identities.) 
\end{Definition}

 Let  $\A=(\d\colon A \to B, {\tra,\tla})$ be a unital crossed module of bare algebras. There is a morphism of unital algebras $\d\colon A_\bullet \to B$, where $\d(\l,a)=\l1_B +\d(a)$, for each $\l \in \C$ and $a \in A$. Let ${\rm Units}(B)$ denote the group of invertible elements (units) of $B$.
The following follows easily from the {definition} of a crossed module of bare algebras.  The fact that $\tla$ and $\tra$ are fully interchangeable plays a fundamental role in proving that the action of ${\rm Units}(B)$ on $A_\bullet$ is by algebra isomorphisms.
\begin{Proposition}
 The map $\big(b,(\l, a)\big) \in {\rm Units}(B) \times A_\bullet \mapsto b \t (\l,a) \doteq (\l,b\tra a \tla b^{-1})$ is an action of the group ${\rm Units}(B)$ on $A_\bullet$ by unital algebra isomorphisms. Moreover, if we pass to the group ${\rm Units}(A_\bullet)$ of units of $A_\bullet$, then ${\rm Units}(\A)\doteq (\d\colon {\rm Units}(A_\bullet) \to {\rm Units}(B), \t)$ is a crossed module of groups. 
\end{Proposition}
\begin{Proof}
 The only unclear bit is the second Peiffer relation. Let $(\l,a)$ and $(\l',a')$ be invertible elements of $A_\bullet$. We want to prove that $\d(\l,a) \t (\l',a')=(\l,a)\, (\l',a') \,(\l,a)^{-1}$. 
 We have, putting $(\l,a)^{-1}=(\l^{-1},\overline{a})$:
 \begin{align*}\allowdisplaybreaks
 \d(\l,a) \t (\l',a') \,\, &=\left (\l' , \d(\l,a) \tra a' \tla \d(\l,a)^{-1}\right)=\big (\l' , (\l1_B+\d(a)) \tra a' \tla (\l^{-1}1_B+\d(\overline{a})\big)\\
 &=\left (\l',  a' +\l^{-1} \d(a) \tra a' + \l a' \tla  \d(\overline{a}) +\d(a) \tra a'  \tla\d(\overline{a})  \right).
 \end{align*}
Whereas:
 \begin{align*}\allowdisplaybreaks
  (\l,a)\, (\l',a') \,(\l,a)^{-1}&=(\l\l', \l'a + \l a'+ aa') \, (\l^{-1},\overline{a}) \\
   &=\big (\l',\l^{-1} \l'a + \l^{-1}\l a'+ \l^{-1} aa' + \l'a\overline{a} + \l a'\overline{a}+ aa'\overline{a}+\l\l'\overline{a}\big ).
 \end{align*}
By using the second Peiffer relation for crossed module of bare algebras, we have:
\begin{align*}\allowdisplaybreaks
 \d(\l,a) \t (\l',a')  - (\l,a)\, (\l',a') \,(\l,a)^{-1}=-\big(0, \l^{-1} \l'a +\l' a \overline{a}+ \l\l'\overline{a}\big)=(0,0),
\end{align*}
since we have $(1,0)=(\l,a) (\l^{-1}, \overline{a})=\big(1,a\overline{a}+ \l \overline{a}+\l^{-1} a\big)$.
\end{Proof}

For this paper we will only need a (slightly restricted) version of the previous construction, however over $\C\h$, by adjoining a formal parameter $h$ to our construction. 
Recalling that $A$ is a bare algebra, we let $A\h$ be the $\Ch$ algebra of formal power series over $A$, with the obvious product, the $\Ch$-linear extension of the product over $A$. We let $\Ad\h$ be the algebra obtained from $A\h$ by (functorially) adjoining an identity element. Therefore $\Ad\h=\Ch \oplus A\h$, as a $\Ch$-module, with product:
$(p,a)(p',a')=\left(pp',p a'+p'a+aa'\right),$ where $p,p'\in \Ch$ and $a,a'\in A\h$. The identity is: $1=\big(1+0h+0h^2+0h^3+\dots,0\big)$.

Recall that any power series, over an algebra, whose first element is the identity is invertible. Let $\Ad^* \subsetneq {\rm Units}(\Ad\h)$ be the group of formal power series in $\Ad$, whose first term is $1_\C=(1_\C,0)$, and such that all other terms are $h^n$ multiples of elements of $A=\{0\} \times A$. Write each element in $\Ad^*$  as $a=1+\sum_{n=1}^{+\infty} a_n h^n$, where $a_n \in A, \forall n$. 
Let $\A=(\d\colon A \to B, \tla,\tra)$ be a unital crossed module of bare algebras. Let $B^*\subsetneq {\rm Units} (B\h)$ be the group of formal power series over $B$ of the form $b=1_B +\sum_{n=1}^{+\infty} h^n {b_n}$, where $b_n \in B, \forall n$. Consider the actions $\tla$ and $\tra$ of $B\h$ on $A\h$ obtained by extending $\tla$ and $\tra$ $\Ch$-linearly. 
\begin{Proposition}
  Let $\A=(\d\colon A \to B, \tra,\tla)$ be a crossed module of bare algebras, with $B$ unital. Consider the induced $\Ch$-algebra  map $\d\colon \Ad\h \to B\h$, restricting to a group map $\d\colon \Ad^* \to B^*.$ With the left action
 $b \t (1,a) \doteq (1,b\tra a \tla b^{-1}),\, \forall b \in B^*, \forall a \in A\h$ we have  a crossed module of groups $\A_\bullet^*=\big(\d \colon  \Ad^* \to B^*, \t \big)$. 
\end{Proposition}
\subsubsection{Crossed modules of bare algebras from chain complexes of vector spaces}\label{cccvs}
The most suggestive non-trivial class of  examples of unital crossed modules of bare algebras, comes from considering chain complexes ${\cal V}=(\dots \ra{\b} V_n \ra{\b} V_{n-1}\ra{\b}\dots)$, of vector spaces. 
We let $\gl_0(\V)$ be the unital algebra of chain maps $f\colon \V \to \V$, with the usual composition of chain maps, and commutator  $\{f,g\}=f\circ g-g \circ f$, where $f,g\in \gl_0(\V)$. We let $\hom^i(\V)$ denote the vector space of degree $i$ maps $\V \to \V$. Composition yields a bilinear map ${\hom^i}(\V) \times {\hom}^j(\V) \to {\hom}^{i+j}(\V)$.  A degree-1 map $s\colon \V \to \V$, will be called a homotopy, and is given by a sequence of vector space maps $\{s_i\colon V_i \to V_{i+1}\}_{i \in \mathbb{Z}}$, with no compatibility relation with respect to the boundary maps $\beta$. There is a map $\d'\colon {\hom^1(\V)} \to {\gl}_0(\V)\subset \hom^0(\V)$ with $\d'(s)=\beta s+s \beta, \, \forall s \in \hom^1(\V)$.
If $f\in \gl_0(\V)$ and $s \in \hom^1(\V)$:
\begin{align}\allowdisplaybreaks
\d'(f \circ s)&=f\circ \d'(s), && \d'( s \circ f)=\d'(s) \circ f.
\end{align}

A degree-2 map $h\colon \V \to \V$ is a sequence of vector space maps  $\{h_i\colon V_i \to V_{i+2}\}_{i \in \mathbb{Z}}$. There is a map $\d'\colon \hom^2(\V) \to \hom^1(\V)$, with $\d'(h)=\beta h - h \beta,\forall h \in \hom^2(\V)$. If $f\in \gl_0(\V)$ and $h \in \hom^2(\V)$, then:
\begin{align}\allowdisplaybreaks
\d'(f \circ h)&=f\circ \d'(h), && \d'( h \circ f)=\d'(h) \circ f.
\end{align}
Put $\gl_1(\V)=\hom^1(\V)/\d'(\hom^2(\V))$. Then clearly $\d'\colon \hom^1(\V) \to \gl_0(\V)$ descends to a linear map $\d\colon \gl_1(\V) \to \gl_0(\V)$. Morever the composition maps: 
\begin{align*}\allowdisplaybreaks
& (f,s) \in \hom_0(\V) \times \hom_1(\V) \mapsto f \circ s \in \hom_1(\V),  & (s,f) \in \hom_1(\V) \times \hom_0(\V) \mapsto s \circ f \in  \hom_1(\V) \end{align*}
descend to maps:
\begin{align*}\allowdisplaybreaks
& (f,s) \in  \gl_0(\V) \times \gl_1(\V) \mapsto s \circ f \in  \gl_1(\V) ,  & (s,f) \in \gl_1(\V) \times \gl_0(\V) \mapsto s \circ f \in  \gl_1(\V) .\end{align*}
\begin{Lemma}
For each $s,t \in \hom^1(\V)$, the following holds:
\begin{equation}
\d'(s\circ t)=\d'(s)\circ t-s \circ \d'(t).
\end{equation}
In particular in $\gl_1 (\V)$ we have:
$$\d(s) \circ t=s\circ  \d(t). $$
\end{Lemma}

We consider the algebra structure in $\gl_0(\V)$ given by the composition of chain-maps. In $\gl_1(\V)$, we consider the following associative algebra structure:
$$s*t\doteq s \circ \d(t)=\d(s) \circ t.$$
\begin{Theorem}
Let $\V=(\dots \ra{\b} V_n \ra{\b} V_{n-1}\ra{\b}\dots)$ be a chain complex of vector spaces. Then $\d\colon \gl_1(\V) \to \gl_0(\V)$, with $\d(s)=s\beta +\beta s, \forall s \in \gl_1(\V)$ is an algebra map. Moreover, with the left and right actions:
\begin{align*}\allowdisplaybreaks
(f,s) \in \gl_0(\V) \times \gl_1(\V) \to \gl_1(\V) &\mapsto f \tra s\doteq f \circ s, \\
(s,f)\in  \gl_1(\V) \times \gl_0(\V) \to \gl_1(\V) &\mapsto s \tla f\doteq s \circ f,
\end{align*}
this defines a unital crossed module of bare algebras, ${\cal HOM}(\V)=\big(\d\colon \gl_1(\V) \to \gl_0(\V),\tra, \tla\big)$.
\end{Theorem}
\begin{Remark}[Crossed modules of groups and Lie algebras from chain-complexes]
Combining the previous {theorem} with  \ref{bdga} and \ref{bag}, yields,  given a chain complex  $\V$ of vector spaces,  a differential crossed module: $$\lgl(\V)=\big(\d\colon \gl_1(\V) \to \gl_0(\V),\t)=\Lie\big({\cal HOM}(\V)\big)$$ and group crossed modules, previously defined (essentially)  in \cite{BL,KP,FMM,CFM1,CFM2}: 
\begin{align*}\allowdisplaybreaks
\GL(\V)&=\big(\d\colon {\rm Units}\big(\gl_1(\V)_\bullet\big) \to {\rm Units}\big(\gl_0(\V)\big),\t\big)={\rm Units}\big({\cal HOM}(\V)\big),\\
\GL(\V,h)&=\big(\d \colon  \gl_1(\V)_\bullet^* \to \gl_0(\V)^*, \t \big)={\cal HOM}(\V)_\bullet^*.
\end{align*}
\end{Remark}
\subsection{Two-groupoids and crossed modules}
\subsubsection{The group crossed module case}\label{gcmc}
For the {definition} of a 2-groupoid and of a sesquigroupoid see \cite{HKK,McLane,St}. Let $\X=(\d \colon E \to G,\t)$ be a group crossed module. We can define a 2-groupoid $\Cc^\times \left ( \X \right)$ out of it, with a single object $*$, which is therefore a categorical group \cite{BaMa}. This is an old construction; see, for example, \cite{BrS} and, more recently \cite{BHS,BL}. We will take as being the ``fundamental'' compositions in $\Cc^\times \left ( \X \right)$  the vertical composition of 2-morphisms, as well as the whiskerings of 2-morphisms by 1-morphisms, as in \cite{St,KP}, since the horizontal composition of 2-morphisms can be derived from these. 
The set of 1-morphisms of  $\Cc^\times\left (\X\right) $ is given by all arrows of the form $*\ra{g} *$, with $g \in G$, the composition being:
$$*\ra{g}*\ra{g'} *=*\ra{gg'} *, \fo g,g' \in G. $$
 The set of 2-morphisms of  $\Cc^\times \left ( \X\right)$  is given by all 2-arrows of the form:
$$\mortimes{g}{e}{\d(e)^{-1} g},\we g \in G, \an e \in E ,$$
each of these having vertical and horizontal inverses of the form (respectively):
$$\mortimes{{\d(e)^{-1} g}}{e^{-1}}{g} \qquad \textrm{ and } \mortimes{g^{-1}}{g^{-1} \t e^{-1}}{g \d{(e)}} .$$
{The horizontal composition of 2-morphisms appears in \eqref{horc}.} The vertical composition of 2-morphisms is:
$$ \mortimesvcomp{g}{e}{\d(e)^{-1}g\,}{f}{\d(f)^{-1}\d(e)^{-1}g} \quad = \nss \mortimes{g}{ef}{\d(ef)^{-1} g},\we g \in G, \an e,f \in E.$$
The left and right whiskerings are (for $g,h \in G$ and $e \in E)$, respectively:
$$\rwtimes{g}{e}{\d(e)^{-1} g}{h} \quad = \hspace*{-.6cm}\mortimes{gh}{e}{\d(e)^{-1} gh},$$
$$\lwtimes{h}{g}{e}{\d(e)^{-1} g}\quad = \hspace*{-.6cm} \mortimes{hg}{h\t e}{h\d(e)^{-1}g} .$$
These whiskering maps are distributive with respect to the vertical composition of 2-morphisms, and are associative, whenever compositions are defined. In addition, the interchange law is   satisfied, namely:
\begin{equation}\label{il}
{\xymatrix{ & {*}\ar[rr]_g \ar@/^3pc/[rr]^{\d(e)^{-1} g}  \ar @/^/ @{{}{ }{}} [rr]^{\boxed{\times}\Uparrow {e}} & &{*}  \ar[rr]^{\d(e')^{-1} g'} \ar@/_3pc/[rr]_{g'}  \ar @/_/ @{{}{ }{}} [rr]_{\boxed{\times}\Uparrow {e'}} & &{*} }}\quad=\quad\hskip-1cm{\xymatrix{ &{*} \ar[rr]^{\d(e)^{-1} g} \ar@/_3pc/[rr]_{g}  \ar @/_/ @{{}{ }{}} [rr]_{\boxed{\times}\Uparrow {e}}  && {*}\ar[rr]_{g'} \ar@/^3pc/[rr]^{\d(e')^{-1} g'}  \ar @/^/ @{{}{ }{}} [rr]^{\boxed{\times}\Uparrow {e'}} & &{*} }},\, \forall  g,g' \in G , \forall  e,e' \in E.
\end{equation}
This interchange law is a well known consequence of the second Peiffer law $\d(e) \t e'=ee'e^{-1},\forall e,e' \in E$, for crossed modules of groups.  
Therefore we can define the horizontal composition of two 2-morphisms as being:
\begin{align}\allowdisplaybreaks&\hspace*{-.6cm}\hspace*{-.6cm}  \mortimeshcomp{g}{e}{\d(e)^{-1}g}{g'}{e'}{\d(e')^{-1} g'}\label{horc}  \\
&\quad\quad  = {\xymatrix{ & {*}\ar[rr]_g \ar@/^3pc/[rr]^{\d(e)^{-1} g}  \ar @/^/ @{{}{ }{}} [rr]^{\boxed{\times}\Uparrow {e}} & &{*}  \ar[rr]^{\d(e')^{-1} g'} \ar@/_3pc/[rr]_{g'}  \ar @/_/ @{{}{ }{}} [rr]_{\boxed{\times}\Uparrow {e'}} & &{*} }} \label{goodexp}\quad \quad 
= \hspace*{-1cm} \mortimes{gg'}{(g\t e')\,e}{\d(e)^{-1} g\d(e')^{-1} g'} \\&= {\xymatrix{ &{*} \ar[rr]^{\d(e)^{-1} g} \ar@/_3pc/[rr]_{g}  \ar @/_/ @{{}{ }{}} [rr]_{\boxed{\times}\Uparrow {e}}  && {*}\ar[rr]_{g'} \ar@/^3pc/[rr]^{\d(e')^{-1} g'}  \ar @/^/ @{{}{ }{}} [rr]^{\boxed{\times}\Uparrow {e'}} & &{*} }}= \hspace*{-1cm} \mortimes{gg'}{e\, (\d(e)^{-1}g) \t e' }{\d(e)^{-1} g\d(e')^{-1} g'}.
 \end{align}

\begin{Remark}\label{referf}
Let $\hat{\X}=(\d \colon E \to G,\t)$ be a pre-crossed module of groups. We can define a sesqui-category  $\hat{\Cc}^\times({\hat{\X}})$ of out ${\hat{\X}}$; see \cite{St} for the {definition} of a sesqui-category, which is similar to the {definition} of a 2-category, except that the interchange law does not necessarily hold. The objects, morphisms and 2-morphisms of $\hat{\Cc}^\times  ({\hat{\X}})$ are constructed exactly in the same way as before, and so are their vertical compositions and whiskerings. However, the interchange law does not hold in general, and therefore the horizontal composition of 2-morphisms is not well defined, given that the two sides of \eqref{il} may yield different results.
\end{Remark} 

\subsubsection{The bare algebras crossed module case}\label{bacmc}
 Let $\A=(\d\colon A \to B, {\tra,\tla})$ be a crossed module of bare algebras. Consider the induced algebra  map $\d\colon A\h \to B\h$.
We define a 2-groupoid $\Cc^+\left (\A\right)$ out of $\A$. The set of 1-morphisms of  $\Cc^+\left (\A\right)$ is given by all arrows $*\ra{b} *$, with $b \in  B\h$. The compositions is:
$$* \ra{b} * \ra{b'} *=* \ra{bb'} *, \we b,b' \in B\h. $$
The set of 2-morphisms of  $\Cc^+\left (\A\right)$  is given by all 2-arrows of the form:
\begin{equation}\label{2mb}
\morplus{b}{a}{b +\d(a)}, \textrm{ where } a \in A\h \textrm{ and } b \in B\h ,
\end{equation}
each of these having a vertical inverse of the form:
$$\morplus{b+\d(a)}{(-a)}{b}.$$
If both $b$ and $b +\d(a)$ are invertible, we also have an horizontal inverse, of the form:
$$\morplus{b^{-1}}{-(b+\d(a))^{-1} \tra a \tla b^{-1}}{\big(b+\d(a)\big)^{-1}} $$
The vertical composition of 2-morphisms is:
$$\morplusvcomp{b}{a}{b+\d(a)}{a'}{b+\d(a)+\d(a')}=\hskip-1cm\morplus{b}{(a+a')}{b+\d(a)+\d(a')}, \we b,b' \in B\h \an a,a' \in A\h.  $$
The left and right whiskerings have the form:
$$\lwplus{b}{b'}{a}{b'+\d(a)}=\hskip-1cm\morplus{bb'}{b\tra a}{bb'+b\d(a)}, \we b,b' \in B\h \an a \in A\h , $$
$$\rwplus{b}{a}{b+\d(a)}{b'}=\hskip-1cm\morplus{bb'}{ a\tla b'}{bb'+ \d(a)b'} , \we b,b' \in B\h \an a \in A\h . $$ 
The interchange law follows from equation \eqref{excl}, which holds in any crossed module of bare algebras. We can therefore define the horizontal composition of two 2-morphisms, as being:
\begin{align}\allowdisplaybreaks
\hspace*{-.6cm}\hspace*{-.6cm}  \morplushcomp{b}{a}{\d(a)+b}{b'}{a'}{\d(a')+ b'}  \quad = \hspace*{-1cm}
\morplus{bb'}{b \tra a'+a \tla b'+aa'}{b+ \d (b \tra a'+a \tla b'+aa') }
\label{phc}
\end{align}

\begin{Remark}Similarly to Remark \ref{referf}, if $\hat{\A}=(\d\colon A \to B, {\tra,\tla})$   is a pre-crossed module of bare algebras,  the construction just outlined does not define a 2-groupoid. We can still define 0, 1 and 2-cells and whiskerings, but the interchange law does not hold, in general, and therefore we can only build a sesqui-category $\hat{\Cc}^+({\hat{\A}})$. 
\end{Remark} 

{This result will be essential for relating the exact and fuzzy formulations of two-dimensional holonomy. }
\begin{Theorem}\label{rel}Let $\A=(\d\colon A \to B, {\tra,\tla})$ be a unital crossed module of bare algebras. Consider the group crossed module $\A_\bullet^*=\big(\d \colon  \Ad^* \to B^*, \t \big)$, defined in \ref{bag}. 
 There exists a 2-functor $\Tc_\A\colon \Cc^\times ( {\A_\bullet^*} )\to \Cc^+\left ({\A}\right)$, defining a natural transformations of functors $\Tc\colon \Cc^\times \Longrightarrow \Cc^+$. It has the form:
$$\Tc_\A\left ( * \ra{b}*  \right )=*\ra{b}*, \textrm{ on 1-cells}$$
and, on 2-cells:
\begin{equation}\label{relf}
\Tc_\A\left( \hspace*{-1cm} \mortimes{g}{e}{\d(e)^{-1} g} \right ) = \hspace*{-1cm} \morplus{g}{(e^{-1}-1)\tla g}{\d(e)^{-1} g}.
\end{equation}
\end{Theorem} 
\begin{Remark}Note that if {$e^{-1}=1+\sum_{n=1}^{+\infty} e_n h^n$, where $e_n \in A, \forall n \in \N$, and where $1$ denotes the identity adjoined to $A\h$,  then $e^{-1}-1=\sum_{n=1}^{+\infty} e_n h^n$.}
\end{Remark}
\begin{Proof}
 The sets of 1-morphisms of $\Cc^\times  ( {\A_\bullet^*})$ and $\Cc^+\left ({\A}\right)$ coincide.  That vertical compositions are preserved by $\Tc_\A$ is consequence of the fact (for each $e,e' \in A_\bullet\h$ and $g \in B\h$, and putting $f=e^{-1}-1$, $f'={e'}^{-1}-1$):
\begin{align*}\allowdisplaybreaks
 ((ee')^{-1}-1)\tla g&= ({e'}^{-1}e^{-1}-1)\tla g =(f' f+f'+f) \tla g = f\tla g + f'f \tla g+f'\tla g \\&=f\tla g + f'\tla (\d(f) g)+f'\tla g=f\tla g + f'\tla\big (\d(f) +1_B) g\big)\\&=(e^{-1}-1)\tla g + (e'-1)\tla\big (\d(e)^{-1}g\big).
\end{align*}
(Note that the second Peiffer condition is explicitly used, and therefore this is a property exclusive to crossed modules of bare algebras, as opposed to pre-crossed modules.)
The preservation of left and right whiskerings (respectively) follows from the fact (for $g,h \in B\h$ and $e \in A_\bullet^*$ and putting $f=e^{-1}-1$):
\begin{align*}\allowdisplaybreaks
 \big((h \t e^{-1})-1\big) \tla (hg)&\doteq  (h \tra f \tla h^{-1}) \tla (hg)= h \tra (f \tla g)= h \tra \big ((e^{-1}-1)\tla g\big)\\
(e^{-1}-1) \tla (gh)&=f \tla(gh)=(f \tla g)\tla h=\big((e^{-1} -1)\tla g\big) \tla h. 
\end{align*}
\end{Proof}

\subsection{Crossed modules of Hopf algebras}
\subsubsection{Hopf algebraic conventions. }
Recall some Hopf algebra conventions; see \cite{Kassel,MajidQG}. 
   Let $H=(H,\D,\e)$ be a bialgebra, with identity $1=1_H$ and counit $\e\colon H \to \C$. The coproduct in $H$ will be written in  {Sweedler}'s 
notation, as in \cite{Kassel}: $$\D(x)= \dis \sum_{(x)} x' \tn x'', \textrm{ for each } x \in H.$$

Suppose that $I$ is a bialgebra, with identity $1_I$. We say that $I$ is a $H$-module algebra  if there exists a left-action $\rho\colon H \times I \to I$, explicitly $\rho\colon (x,v) \in H \times V \mapsto x \t_\r v \in V,$ such that:
\begin{itemize}
  \item  $1_H \tr v=v$ for each $v \in I$.
 \item $\dis x \tr (uv)=\sum_{(x)} (x' \tr u) \, (x''\tr v)$, for each $u,v \in I$, and each $x \in H$.
\item $x \tr 1_I=\e(x) 1_I$, for each $x \in H$.
 \end{itemize}
 With the algebra action $\rho\colon H \tn I \to I$, $I$ is said to be an $H$-module coalgebra if we have:
 \begin{itemize}
 \item $\dis \Delta (x \tr v)=\sum_{(x)} \sum_{ (v)} (x' \tr v') \tn (x'' \tr v'')$ for each $x \in H$ and each $v \in V$.
\item $\epsilon ( x \tr v)=\e(x) \e(v)$ for each $x \in H$ and each $v \in V$.
\end{itemize}
The following is well known and a proof appears in \cite{MajidQG}.
\begin{Theorem}[Majid]\label{majmain0}
 If $I$ is an $H$-module algebra and coalgebra, under the left action $\rho$, then we can define an algebra $I \otimes_\rho H$, with  underlying vector space $I \otimes H$, with identity $1_I \otimes 1_H$, 
and with product:
$$(u \tn x)(v \tn y)=\sum_{(x)}\big( u \,\, x' \tr v \big)\tn\big( x'' y \big), \we u,v \in I \an x,y \in H. $$
If, additionally, the following compatibility condition holds:
\begin{equation}\label{comp}
\sum_{(x)} x' \tn ( x''\tr v)=\sum_{(x)} x'' \tn ( x'\tr v), \textrm{ for each } x \in H \textrm{ and } v \in I, 
\end{equation}
 then we have a bialgebra structure in $I \otimes_\rho H$, with coproduct:
$$\D(u \tn x)= \sum_{(u)} \sum_{(x)} (u'\tn x') \tn (u'' \tn x'') , \we x \in H \an u \in I,$$
(the compatibility condition \eqref{comp} is only used for proving that the coproduct is an algebra morphism),
and: 
$$\epsilon (u \tn x)= \e(u) \e(x), \we x \in H \an u \in I.$$
If, in addition, $H$ and $I$ are Hopf algebras, then $I \otimes_\rho H$ is a Hopf algebra, with antipode:
$$S(u \tn x)=\big ( 1_I \tn S(x) \big)\, (S(v) \tn 1_H),  \we x \in H \an u \in I.  $$
Moreover the inclusions of $H$ and $I$ in  $I \otimes_\rho H$, namely:
\begin{align*}\allowdisplaybreaks
& i_H \colon x \in H  \mapsto 1_I \tn x \in I \otimes_\rho H, 
 &i_I \colon v \in I \mapsto v \tn 1_I \in I \otimes_\rho H 
\end{align*}
are bialgebra morphisms (and Hopf algebra morphisms if $H$ and $I$ are Hopf algebras). 
\end{Theorem}
\subsubsection{The adjoint representation. Primitive and group like elements}

All of this is well known. 
If $H$ is a Hopf algebra, then $H$ has a $H$-module algebra structure, given by:
$$\rho(x \tn y)=\sum_{(x)} x' y S(x'')\doteq  x \ad y,  $$
called the adjoint representation. (In general we do not have an $H$-module coalgebra structure, unless, for example,  $H$  is cocommutative. In fact $\D(x \t y)=\sum_{(x)}\sum_{(y)} x' y' S(x'''') \tn x'' {y''} S(x''') $, for each $x,y \in H$.)

Recall that an element $x \in H$ is said to be:
\begin{itemize}
 \item Grouplike, if $\D(x)=x \tn x$.
\item Primitive, if $\D(x)=x \tn 1 + 1 \tn x$. 
\end{itemize}
It is easy to see that $\epsilon(x)=1$, if $x$ is group-like, and $\epsilon(x)=0$, if $x$ is primitive. Any group-like element $x$ is invertible, being $S(x)=x^{-1}$. If $x$ is primitive then $S(x)=-x$. Therefore:
\begin{Lemma}[The adjoint action of group-like and of primitive elements]
Let $y \in H$. Then:
\begin{itemize}
 \item If $g\in H$ is group-like then $g \ad y=gyg^{-1}$.
 \item If $x \in H$ is primitive then $x \ad y=xy-yx$. 
\end{itemize}

\end{Lemma}

\subsubsection{Crossed modules and pre-crossed modules of Hopf algebras}\label{cmh}
The following {definition} appeared in \cite{MajidHA}. It is more restrictive than the one appearing in \cite{FW,FLV} in that condition \eqref{comp} is imposed from the beginning.  
\begin{Definition}[Crossed modules and pre-crossed module of Hopf algebras]
 A Hopf algebra crossed module $\H=(\d \colon I \to H, \rho)$ is given by:
\begin{itemize}
 \item A map $\d\colon I \to H$ of Hopf algebras.
 \item A left action   $\rho\colon x \tn v \in H \otimes I \mapsto x \tr v \in I,$ which turns $I$ into an $H$-module algebra and coalgebra, satisfying moreover the compatibility condition \eqref{comp}.
\end{itemize}
We furthermore impose the first and second Peiffer laws, as displayed below:
\begin{enumerate}
 \item $\d(x \tr v)=x \ta \d(v)$ \quad  (first Peiffer Law), where $x \in H$ and $v \in I$. 
 \item  $\d(u) \tr v=u \ta v$ \quad (second Peiffer Law), where $u,v \in I$.
\end{enumerate}
If the first Peiffer law is satisfied (but not necessarily the second) for $\hat{\H}=(\d \colon I \to H, \rho)$, then $\hat{\H}$ is called a pre-crossed module, of Hopf algebras. 
\end{Definition}
The following theorem appeared in \cite{MajidQG}, being a simple extension of Theorem \ref{majmain0}
\begin{Theorem}[Majid]\label{majmain}
Suppose that $\H=(\d \colon I \to H, \rho)$ is a pre-crossed module of Hopf algebras. We have Hopf algebra maps (called source and target maps) $s,t\colon I \otimes_\rho H \to H$, which read on generators:
\begin{align*}\allowdisplaybreaks
&s\colon v \tn x \in  I \otimes_\rho H \mapsto \e(v) x \in H, 
&t\colon v \tn x \in  I \otimes_\rho H \mapsto \d(v) x \in H.
\end{align*}
Moreover, together with the inclusion $i_H \colon H \to I \otimes_\rho H$ (called identity map) these define a  graph in the category of Hopf algebras, namely we have $s \circ i_H=\id_H \textrm{ and } t \circ i_H=\id_H. $
\end{Theorem}
\begin{Remark} Even though we do not assume our Hopf algebras to be cocommutative, the conditions we impose are still quite restrictive, however perfectly adapted to the type of examples we have in mind: crossed modules of Hopf algebras arising from differential crossed modules. In \cite{MajidQG} a much more general definition of crossed modules of (braided) Hopf algebras appears, enabling for {\it genuinely} non-cocommutative examples. 
\end{Remark}

\subsubsection{From crossed module of Hopf algebras to group crossed modules}\label{hag}
Let $H$ be a Hopf algebra. We let $ H[[h]]$ be the ring of formal power series over $H$. This inherits an obvious algebra structure, with ground ring $\C[[h]]$.
Namely
$\left ( \sum_{n=0}^{+\infty} h^n a_n \right)\left (\sum_{m=0}^{+\infty} h^m b_m \right)= \sum_{n=0}^{+ \infty} h^n \sum_{i+j=n} a_i b_j, $
clearly  $\C[[h]]$-bilinear. 
Note, however, that $H[[h]]$ is not necessarily a bialgebra, with ground ring $\C[[h]]$. The (very well known) point here is that, given a formal power series $\sum_{n=0}^{+\infty} h^n a_n \in H\h$, its coproduct would be $\sum_{n=0}^{+\infty} h^n \Delta(a_n)$, which lives in $(H \otimes H)[[h]]$, strictly containing  $H[[h]] \otimes_{\C[[h]]} H[[h]].$
Nevertheless, the $\Ch$-algebra $\Hh$ can be seen as being a Hopf algebra, however such that the coproduct $\D$ is a map $\D\colon \Hh \to (H \tn H)[[h]] \supsetneq  H[[h]] \otimes_{\C[[h]]} H[[h]].$ (This can be dealt with fully precisely by considering a completion of the tensor product of $\Ch$-modules, as explained, for example, in \cite[XVI.1]{Kassel}. We will not use that language since a completely algebraic proof can always be given.)

\begin{Definition}[Group-like element in $H\h$] Let $H$ be a Hopf algebra. 
 An element of $x \in \Hh$ is said to be group like if $\D(x)=x\tn x$. Explicitly, $ x=\sum_{n=0}^{+\infty} h^n x_n $ is said to be group like if
$\sum_{n=0}^{+\infty} h^n \D(x_n) =\sum_{n=0}^{+ \infty} h^n\sum_{i+j=n} x_i \tn x_j, $
that is, using {Sweedler}'s notation:
$$ \D(x)=x \tn x\quad  \Longleftrightarrow \quad \sum_{n=0}^{+\infty} h^n \sum_{(x_n) } x_n' \tn x_n'' =\sum_{n=0}^{+ \infty} h^n\sum_{i+j=n} x_i \tn x_j. $$
\end{Definition}
 \noindent For example the identity $1_H$ of $H$ (and therefore of $\Hh$) is group like. Group-like elements of $\Hh$ are closed under multiplication.  Let $a \in H$ be a primitive element. Then $ \exp(ha)=\sum_{n=0}^{+ \infty} h^n \frac{a^n}{n!}$ is a group like element of $\Hh$.
If $ x=\sum_{n=0}^{+ \infty} x_n h^n \in H\h$ is group like, then, similarly to the case of Hopf algebras over a field, $\epsilon(x)=\sn h^n \epsilon (x_n)=1_H $ and  $x S(x)=1_H=S(x)x$. The antipode $S$ is extended $\C\h$-linearly.

Let $\H=(\d \colon I \to H, \rho)$ be a crossed module of Hopf algebras. Extend the action of $H$ on $I$, $\Ch$-linearly:
$$\left (\sum_{n=0}^{+ \infty} x_n h^n \right) \tr \left(\sum_{m =0}^{+\infty} v_m h^m \right)= \sum_{n=0}^{+ \infty} h^n\sum_{i+j=n} x_i \tr v_j, \we x_n \in H \an v_n \in I, \fo n \in \N.$$
Let $H^*_{\rm gl}$ and $I^*_{\rm gl}$ be the groups of invertible group like elements of $\Hh$ and of  $\Ih$, whose underlying power series start at the identity (of $H$ or $I$). Consider the induced map {$\d\colon \Ih \to \Hh$,} sending group like elements to group-like elements, and invertible elements to invertible elements. The restriction {$\d\colon I^*_{\rm gl} \to H^*_{\rm gl}$} of $\d\colon \Ih \to \Hh$ is a map of groups. In addition: 
\begin{Lemma}
The action $\tr$ of $\Hh$ on $\Ih$ restricts to a group action of $H^*_{\rm gl}$ on $I^*_{\rm gl}$, by automorphisms. 
\end{Lemma}
\begin{Proof}
 Let $x,y \in H^*_{\rm gl}$ and $u,v \in I[[h]]$. Put $ x=\sum_{n=0}^{+ \infty} h^n x_n$, $ v=\sum_{n=0}^{+ \infty} h^n v_n$, $ u=\sum_{n=0}^{+ \infty} h^n u_n$. That $(xy) \tr u=x\tr (  y \tr u)$ follows since $\rho$ is an algebra action. 
That $x \tr (uv)=(x \tr u) (x \tr v)$ follows from:
\begin{align*}\allowdisplaybreaks
 &x \tr (uv)=\sum_{n=0}^{+\infty} h^n \sum_{i+j+k=n} x_i \tr (u_j v_k) =\sum_{n} h^n \sum_{i+j+k=n}\,\,\sum_{(x_i)} \left( x_i' \tr  u_j\right) \, \left (x_i'' \tr  v_k \right) \\
         &=\sum_{n=0}^{+\infty} h^n \sum_{i+j+k=n} \,\,\sum_{o+p=i} \left( x_o \tr  u_j\right) \, \left (x_p \tr  v_k \right)=\sum_{n} h^n \sum_{o+p+j+k=n} \,\,\left( x_o \tr  u_j\right) \, \left (x_p \tr  v_k \right)=(x \tr u) \,\,(x \tr v).
\end{align*}
 In fact what we proved is that the action of $H^*_{\rm gl}$ on $\Ih$ is a group action by  algebra isomorphisms. 
\end{Proof}

\begin{Lemma}
 Let $x,y \in A[[h]]$, where $A$ is a Hopf algebra, with $x$ invertible and group-like.  Then:
\begin{equation}\label{act}
 x \ad y=x\,\, y\,\,x^{-1}. 
\end{equation}
\end{Lemma}
\begin{Proof}  Put $x=\sum_{n=0}^{+ \infty} h^n x_n$ and $y=\sum_{n=0}^{+ \infty} h^n y_n$. We have: 
\begin{align*}\allowdisplaybreaks
 x \ad y  &= \left ( \sum_{n=0}^{+ \infty} h^n x_n \right) \ad \left (\sum_{n=0}^{+ \infty} h^n y_n \right)= \sum_{n=0}^{+ \infty} h^n \sum_{i+j=n} x_i \ad y_j =\sum_{n=0}^{+ \infty} h^n \sum_{i+j=n} \sum_{(x_i)} x_i' \, y_j\, S(x_i'')\\&=\sum_{n=0}^{+ \infty} h^n \sum_{i+j=n}\,\, \sum_{o+p=i} x_o \, y_j\, S(x_p)=\sum_{n=0}^{+ \infty} h^n \sum_{o+p+j=n}\,\, x_o \, y_j\, S(x_p)=xyx^{-1}.
\end{align*}
\end{Proof}

Combining the two previous lemmas we see that.
\begin{Theorem}\label{rrr}Let $\H=(\d \colon I \to H, \rho)$ be a crossed module of Hopf algebras.
 We have a crossed module of groups $\H^*_{\rm gl}=(\d\colon I^*_{\rm gl} \to H^*_{\rm gl},\tr)$. 
\end{Theorem}
\begin{Remark}\label{impcor}
Consider a Hopf algebra crossed module $\H=(\d \colon I \to H, \rho)$.
 Looking at  the two previous lemmas,  given $x \in H^*_{\rm gl}$, then, for all $y \in I\h$, we have
$\d(x \t y)=x \t_{{\rm ad}} \d(y)= x\,\d(y)\, x^{-1}. $
Define $\hat{I}^*_{\rm gl}\supset {I}^*_{\rm gl} $ as being the group of invertible elements of $I\h$,  mapping, through $\d\colon I\h \to H\h$,  to group-like elements of $H\h$, of the form $1_H +\sum_{n=1}^{+\infty} h^n x_n$. Clearly
 $\hat{\H}^*_{\rm gl}=(\d\colon \hat{I}^*_{\rm gl} \to H^*_{\rm gl},\tr)$ is a  pre-crossed module of groups.
\end{Remark}

\subsubsection{From crossed modules of Hopf algebras to differential (pre)-crossed modules}\label{pcm}
Let $H$ be an Hopf algebra. Let $g,x,y \in H$. If $g$ is group-like then $g \ad y=gyg^{-1}$. If $x$ is primitive then $x \ad y=[x,y]=xy-yx$. The vector space  of primitive elements of $H$ is a Lie algebra, denoted by $\Prim(H)$, and the adjoint action by a  group-like element $g \in H$ is a Lie algebra map $x \in \Prim(H) \mapsto g \ad x \in \Prim(H)$.   We denote by $\widehat{\Prim}(I)$ the vector space of elements of $I$ mapping through $\d$ to $\Prim(H)$. 

\begin{Proposition}
  Let $\H=(\d \colon I \to H, \rho)$ be a Hopf algebra crossed module. Then we have a differential crossed module $\Prim(\H)=(\d\colon \Prim(I) \to \Prim(H),\t_\rho)$ and a differential pre-crossed module $\widehat{\Prim(\H)}=(\d\colon \widehat{\Prim(I)} \to \Prim(H),\t_\rho)$
\end{Proposition}
We will consider the following bilinear map in $\widehat{\Prim(I)}$, vanishing in $\Prim(I) \tn \Prim(I)$:
\begin{equation}\label{ppps}
u \tn v \mapsto \langle u, v \rangle  \doteq [u,v] -\d(u) \t v, \we u,v \in  \widehat{\Prim(I)}.
\end{equation}
{This map measures how  far the differential pre-crossed module $\widehat{\Prim(\H)}$ is from being a crossed module. It will be used when discussing the relation between holonomy and curvature; see Subsection \ref{2cccc}.}  
\subsection{Universal enveloping algebras}
\subsubsection{Review on universal enveloping algebras}

Consider the underlying vector space functor, from the category of (unital) associative algebras to the category of vector spaces, and let $\T$ be its left adjoint, the tensor algebra functor, thus: 
$$\T(V)=\bigoplus_{n=0}^{+\infty} V^{n \otimes}, \textrm{ where  } V^{n \otimes}= \underbrace{V \otimes \dots \otimes V}_{n \textrm{ times}}, \textrm{ and } V^0\doteq \C, $$
with multiplication being the tensor product. Given a vector space $V$,  the tensor algebra $\T(V)$ can be given a Hopf algebra structure by putting $\Delta(x)=x \tn 1 + 1 \tn x$ and $\e(x)=0$, on generators $x \in V$, with antipode $S(x)=-x, \forall x \in V$. If $B$ is a basis of $V$, then $\T(V)$ is isomorphic to the algebra $\C\langle B\rangle$ of non commutative polynomials, with coefficients in $\C$, and with a formal variable for each element of $B$. Therefore $\T(\C)\cong \C[x]$.

We follow \cite{Var}, closely. If $A$ {is} a bare algebra, recall that an algebra derivation $f\colon A \to A$ is a linear map such that $f(aa')=f(a)\,a'+a\,f(a'), \forall a,a' \in A$. Given a Lie algebra $\lg$, a Lie algebra derivation is a map $f\colon \lg \to \lg$, such that $f([x,y])=[f(x),y]+[x,f(y)], \forall x,y \in \lg$. We denote by $\Der(A)$ and $\Der(\lg)$ the Lie algebras of algebra derivations of $A$ and of Lie algebra derivations of $\lg$.
Any vector space map $f\colon V \to V$ gives rise to a unique derivation $\der(f)$ of $\T(V)$, extending $f$. This yields a Lie algebra map $\der\colon f \in \gl(V) \mapsto \der(f) \in \Der(\T(V))$, where $\gl(V)$ is the Lie algebra of linear maps $V \to V$. 

Let $\lg$ be a Lie algebra. We denote its universal enveloping algebra by $\U(\lg)$. This is the quotient of $\T(\lg)$ by the bilateral ideal $L$ generated by $\{XY-YX-[X,Y], \forall X,Y \in \lg\}\subset {\T(\lg)}$. This is a Hopf-algebra ideal: $\D(L) \subset L \tn \T(\lg)+\T(\lg) \tn L$, $\e(L)=\{0\}$ and $S(L)\subset L$. Therefore the universal enveloping algebra $\U(\lg)$ is a Hopf algebra. Given $X \in \lg$ the coproduct is $\D(X)=X \tn 1 + 1 \tn X$, and the counit is $\e(X)=0$. 
If $f\colon \lg \to \lg$ is a Lie algebra derivation, then $\der(f)(L)\subset L$, thus we have a Lie algebra map $\der\colon \Der(\lg) \to \Der(\U(\lg))$.

Recall that the structure of $\U(\lg)$ is unraveled by the Poincar\'{e}-Birkhoff-Witt (PBW) theorem. One of the numerous equivalent ways to state {it} is the following: if $\{X_i, i \in I\}$ is a basis of $\lg$, where $I$ is a totally ordered set,  then $\{1\} \cup \{X_{i_1} X_{i_2}\dots X_{i_n}, n \in \N, i_1 \leq i_2 \leq \dots \leq i_n\}$ is a basis of $\U(\lg)$. The PBW theorem allows us to consider $\lg$ as being, canonically a subvector space of $\lg$. (A fact that will be used without explanation.)

\subsubsection{The restricted universal enveloping algebra $\U^0(\lg)$}
Given a Lie algebra $\lg$, we will also consider the restricted universal enveloping algebra $\U^0(\lg)$, an algebra without unit, which is defined as being the kernel $\ker(\e) \subset \U(\lg)$ of the counit $\e\colon \U(\lg) \to \C$, therefore being the augmentation ideal of $\U(\lg)$. By the PBW theorem, $\{X_{i_1} X_{i_2}\dots X_{i_n}, n \in \N, i_1 \leq i_2 \leq \dots \leq i_n\}$ is a basis of $\U^0(\lg)$. Clearly (from the PBW theorem) ${\U^0(\lg)}$ is  the smallest non-unital subalgebra of $\U(\lg)$ containing $\lg$. One recovers $\U(\lg)$ from ${\U^0(\lg)}$ by formally adjoining an identity; in other words $\U^0(\lg)_\bullet\cong \U(\lg)$.

\subsubsection{The sharp $\xi^\sharp_A\colon \U(\Lie(A)) \to A_\bullet$ and flat $\xi^\flat_A\colon \U^0(\Lie(A))\to A$ counits for a  bare algebra $A$}\label{mcu}

Recall the adjunction between unital algebras and Lie algebras, whose right adjoint associates to a unital algebra $A$ its commutator algebra ${\rm Lie}(A).$ The left adjoint is precisely the universal enveloping algebra functor. The counit of this adjuntion is the natural transformation $\xi$, which, to a unital  algebra $A$, associates the unital algebra map $\xi_A\colon \U(\Lie(A)) \to A$, which, on generators $x \in A$,  is $x \mapsto x$. Given a Lie algebra $\lg$ and a unital algebra $A$, any Lie algebra map $f\colon \lg \to {\rm Lie}(A)$ yields a unique unital algebra map $\hat{f}\colon \U(\lg) \to A$, which, on generators $x \in \lg$, is $x \mapsto f(x)$. Let us now suppose that $A$ does not necessarily have an identity.  
 \begin{Definition}[Sharp counit]
  Let $A$ be a bare algebra. Let $A_\bullet=\C\oplus A$ be the  unital algebra obtained from $A$ by adjoining a unit. Let $\i\colon A \to A_\bullet$ be the inclusion map $x \mapsto (0,x)$, inducing a Lie algebra map $\i\colon {\rm Lie}(A) \to {\rm Lie}(A_\bullet)$.  The sharp counit is the unital algebra map:
$$\xi^\sharp_A \colon \U(\Lie(A)) \to A_\bullet,  $$
universally provided by the Lie algebra map  $\i\colon {\rm Lie}(A) \to {\rm Lie}(A_\bullet)$.
Therefore, on generators $x \in A$, we have $\xi^\sharp_A(x)=(0,x)$, and finally $\xi^\sharp_A(1)=(1,0)$. Therefore, clearly $\xi^\sharp_A\big(\U^0(\Lie(A))\big)=\{0\}\times A\subsetneq A_\bullet.$ 
 \end{Definition}

 Let $A$ be a bare algebra. Consider the obvious algebra isomorphism $p_A\colon \{0\} \times A \subset A_\bullet \to A$.
 \begin{Definition}[Flat counit]
  Let $A$ be a bare algebra. The flat counit is the bare algebra map: 
  $$\xi^\flat_A \colon {\U^0(\Lie(A))} \to A,  $$
  given by $\xi^\flat_A=p_A \circ \epsilon^\sharp_A$. Therefore, on generators $a \in A$, we have: $\xi^\flat_A(a)=a$. Note that $\xi_A^\flat$ is surjective. 
 \end{Definition}

The categorical explanation of the flat counit is the following: we have an adjunction between the  category of Lie algebras and the category of (not-necessarily unital) algebras; i.e. bare algebras. The {right adjoint} assigns to any associative algebra its commutator Lie algebra, whereas the left adjoint is the restricted universal enveloping algebra. The flat counit is exactly the counit of this adjunction.

\subsubsection{From differential crossed modules to crossed module of Hopf algebras}\label{bcvcv}

As proved in \cite{MajidHA}, a differential crossed module  ${\mathfrak{X}=(\d \colon \ae \to  \lg,\t )}$
gives rise to a crossed module of Hopf algebras
$\U(\mathfrak{X})=(\d \colon\U( \ae) \to  \U(\lg),\rho )$. Let us spell out the construction. The algebra map $\d \colon\U( \ae) \to  \U(\lg)$ is the map functorially given by the Lie algebra map $\d\colon \ae \to \lg$. The action $\t$ of $\lg$ on $\ae$, by derivations, yields a Lie algebra map $\lg\to \Der(\ae)$, and therefore, {composing with ${\rm der}\colon \Der(\ae) \to \Der(\U(\ae))$,} a Lie algebra map $\lg \to \Der(\U(\ae))\subset \gl(\U(\ae))$,  the  Lie algebra of linear maps $\U(\ae) \to \U(\ae)$. Thence we have  an algebra map $\U(\lg) \to {\rm hom}(\U(\ae))$, {where ${\rm hom}(\U(\ae))$ is the algebra} of linear maps $\U(\ae) \to \U(\ae)$, inducing an  action:
$$\rho\colon (x,v) \in \U(\lg) \otimes 	\U(\ae) \mapsto x \t_\r v \in 	\U(\ae).$$
\begin{Theorem}[Majid] If ${\mathfrak{X}=(\d \colon \ae \to  \lg,\t )}$ is a differential crossed module then $\U(\mathfrak{X})=(\d \colon\U( \ae) \to  \U(\lg),\r )$ is a crossed module of Hopf algebras.
\end{Theorem}
\subsubsection{The crossed module counit}\label{cmcu}

Suppose that $\A=(\d\colon A \to B, {\tra,\tla})$ is a unital crossed module of bare algebras. Pass to the associated crossed module of Lie algebras $\Lie(\A)=(\d\colon \Lie(A) \to \Lie(B), \t)$; see \ref{bdga}. This will yield a crossed module of Hopf algebras   $\U(\Lie(\A))=(\d\colon \U(\Lie(A)) \to \U(\Lie(B)), \r)$. Recall the notation of \ref{bag}, \ref{hag} and \ref{mcu}. Consider the flat counit $\xi_A^\flat\colon \U^0(\Lie(A)) \to A$ and the usual counit $\xi_B \colon \U(\Lie(B)) \to B$. Let us look at the pair $\left(\xi^\flat_A,\xi_B\right)$. We will denote by the same symbols the extensions   $\xi_A^\flat\colon \U^0(\Lie(A))[[h]] \to A[[h]]$ and  $\xi_B \colon \U(\Lie(B))[[h]] \to B[[h]]$. Clearly the diagram below  (in the category of bare algebras) commutes:
$$\xymatrix{& \U^0(\Lie(A))[[h]]\ar[r]^\d\ar[d]_{\xi_A^\flat} & \U(\Lie(B))[[h]] \ar[d]^{\xi_B}  \\ &A[[h]] \ar[r]_\d  &B[[h]] } $$
We cannot conclude that {$(\xi^\flat_A,\xi_B)$} is a map of crossed modules of bare algebras, among other reasons since $\U(\Lie(\A))=(\d\colon \U^0(\Lie(A)) \to \U(\Lie(B)), \t_\r)$ does not, a priori, have the structure of such; this will be expanded in  \ref{qwer}. However, passing to the associated crossed modules of groups, the pair $\left(\xi^\flat_A,\xi_B\right)$ induces a group crossed module map, the crossed module counit, ${\rm Proj}_\A\colon \U(\Lie(\A))^*_{\rm gl} \to \A^*_\bullet.$ We now address this.

\begin{Lemma}\label{www}
For each $x \in B \subset \U(\Lie(B))$, and each  $a \in {\U^0(\Lie(A))}$, we have:
$$ {\xi^\flat_A( x \t a)=\xi_B(x) \tra \xi^\flat_A(a) - \xi^\flat_A(a) \tla \xi_B(x).}$$ 
\end{Lemma}
\begin{Proof}
It suffices  proving it for $a=v_1*\dots *v_n$, where $v_i \in A$. Use $*$ to denote the product in $\U(\Lie(A))$.
\begin{align*}\allowdisplaybreaks
 \xi^\flat_A(x \t a)&=\xi^\flat_A\big(x \t (v_1* v_2 *\dots *v_n)\big)\\
                   &=\xi^\flat_A\big((x \t v_1)*v_2* \dots* v_n)\big)+\xi^\flat_A\big( v_1* (x \t v_2)* \dots * v_n)\big)+\dots + \xi^\flat_A\big( v_1* v_2 *\dots *(x \t v_n) )\big)\\
 &=(x \t v_1)v_2 \dots v_n+ v_1 (x \t v_2) \dots v_n+\dots + v_1 v_2 \dots (x \t v_n) \\
&=(x \tra v_1-v_1 \tla x)v_2 \dots v_n+ v_1 (x \tra v_2-v_2 \tla x) \dots v_n+\dots + v_1 v_2 \dots
(x \tra v_n-v_n \tla x)\\
&=x \tra (v_1 v_2 \dots v_n)-(v_1 v_2 \dots v_n) \tla x=\xi_B(x) \t \xi^\flat_A(a). 
\end{align*}
(Note that the left and right actions $\tra$ and $\tla$ are fully interchangeable.)
\end{Proof}

We emphasize  that, in general, if $w\in \U(\Lie(B))$ then $ \xi^\flat_A( w \t a)\neq \xi_B(w) \tra \xi^\flat_A(a) - \xi^\flat_A(a) \tla \xi_B(w).$
Consider the algebra $\hom(A)$ of linear maps $A \to A$. Then we have a unital algebra map $R\colon \U(\Lie(B)) \to \hom(A)$, which, on algebra generators $x \in B$, is such that $R(x)(a)=x \tra a -a \tla x=x \t a$. (Considering the Lie algebra action $\t$ of $\Lie(B)$ on $A$, this is just the universal enveloping algebra version of $\t$.)
\begin{Lemma}
 For each $w$ in $\U(\Lie(B))$, and each  $a \in A$,  we have:
\begin{equation}\label{je}
 R(w)(a)=\sum_{(w)} \xi_B(w')\tra\, a\,\, \tla\xi_B(S(w'')).
 \end{equation}
\end{Lemma}
\begin{Proof}
Let us see that the right-hand-side of \eqref{je} defines an algebra map $T\colon \U(\Lie(B)) \to \hom(A)$, where $$T(w)(a)= \sum_{(w)} \xi_B(w')\tra\, a\,\, \tla\xi_B(S(w'')), \textrm{ where } a \in A \textrm{ and } w \in \U(\Lie(B)) .$$ If $x,y \in  \U(\Lie(B))$ and $a \in A$ we have:
 \begin{align*}
   T(xy)(a)&\doteq \sum_{(xy)} \xi_B((xy)')\tra\, a\,\, \tla\xi_B(S((xy)'')= \sum_{(x)} \sum_{ (y)} \xi_B(x'y')\tra\, a\,\, \tla\xi_B\big (S(y'') S(x'')\big)\\
   &= \sum_{(x)} \sum_{(y)} \xi_B(x') \tra\left( \xi_B(y')\tra\, a\,\, \tla\xi_B\big (S(y'')\big)\right) \tla \xi_B\big( S(x'')\big)=T(x)\big( T(y)(a)\big). 
 \end{align*}
Now just use the fact that, clearly, $R$ and $T$ coincide on the algebra generators $x \in B$, all primitive. 
\end{Proof}

Now consider the associated algebras  over $\C[[h]]$, $\U(\Lie(B))[[h]], \U^0(\Lie(A))[[h]], B[[h]]$ and $A[[h]]$.
\begin{Lemma}
For each $w \in B \subset \U(\Lie(B))^*_{\rm gl}$, and each  $a \in {\U^0(\Lie(A))}[[h]]$, it holds that:
$$ \xi^\flat_A( w \t a)=\xi_B(w) \tra \xi^\flat_A(a) \tla  \xi_B(w)^{-1}.$$ 
\end{Lemma}
\begin{Proof}
If $w$ is group like,  then $\D(w)=w \tn w$, being $S(w)=w^{-1}$. In particular we have a bare algebra map $a \in {\U^0(\Lie(A))}[[h]] \mapsto w \t a \in  {\U^0(\Lie(A))}[[h]]$. Thus also  $a \in {\U^0(\Lie(A))}[[h]] \mapsto \xi^\flat_A\big(w \t a\big) \in A[[h]] $ is a bare algebra map. Clearly $a  \in {\U^0(\Lie(A))}[[h]] \mapsto\xi_B(w) \tra \xi^\flat_A(a) \tla  \xi_B(w)^{-1} \in A[[h]]$ is a bare algebra map. The latter two maps coincide on generators $a \in A \subsetneq {\U^0(\Lie(A))}[[h]]$, by the previous lemma.  
\end{Proof}

\begin{Theorem}[Crossed module counit]\label{zxzx}
{Let $\A=(\d\colon A \to B, {\tra,\tla})$ be a unital crossed module of bare algebras over $\C$.
  The algebra maps $\xi^\sharp_A\colon  \U(\Lie(A)) \to A_\bullet$ and  $\ve_B \colon  \U(\Lie(B)) \to B$ induce a group crossed module map ${\rm Proj}_\A\colon \U(\Lie(\A))^*_{\rm gl} \to \A^*_\bullet,$ called the crossed module counit.}
\end{Theorem}
\begin{Proof}{
We must prove that ${\rm Proj}_\A$ intertwines the group actions. Let $1_{\U(\Lie(A))}+\sum_{n=1}^{+\infty} h^n a_n \in \U(\Lie(A))^*_{\rm gl}$. Let $x \in \U(\Lie(B))^*_{\rm gl}$. By using the previous lemma, and the fact $\epsilon(x)=1$, we have:}
\begin{align*}
 \xi_A^\sharp&\Big(x \t \big (1 _{\U(\Lie(A))} +\sum_{n=1}^{+\infty} h^n a_n \big)  \Big)= 
 \xi_A^\sharp\Big(\epsilon(x) 1 _{\U(\Lie(A))} + x \t \sum_{n=1}^{+\infty} h^n a_n \big)  \Big)=
 \left ( 1,  \ve^\flat_A\Big( x \t \sum_{n=1}^{+\infty} h^n a_n \Big)  \right)\\&=  \left ( 1, \ve_B(x) \tra \ve^\flat_A\Big(\sum_{n=1}^{+\infty} h^n a_n \Big)\tla \ve_B(x)^{-1}  \right)\doteq \xi_B(x) \t \left ( 1,  \ve^\flat_A\Big(\sum_{n=1}^{+\infty} h^n a_n\Big) \right)\doteq\xi_B(x) \t \xi_A^\sharp\left( 1_{\U(\Lie(A))}+\sum_{n=1}^{+\infty} h^n a_n \right).
\end{align*}

\end{Proof}

\subsubsection{The bare counit}\label{qwer}
Recall Majid's theorems \ref{majmain0} and \ref{majmain}.  Let  $\H=(\d \colon I \to H, \rho)$ be a Hopf algebra crossed module. Consider $I^0\doteq \ker (\epsilon \colon I \to \C)\subset I$. This is a non-unital bare algebra, closed under the action $\rho$ of $H$ on $I$, since $I$ is an $H$-module coalgebra. This permits us to consider the bare algebra $\ker (\epsilon \colon I \to \C) \otimes_\rho H=I^0\otimes_\rho H$. Recall that we have an inclusion $i_H \colon H \to I \otimes_\rho H$, with $i_H(x)=1_I \tn x$. If $v \tn x \in I^0 \tn H$ then, if $y \in H$, we have:
$$(v \tn x) (1_H \tn y)= \sum_{(x)} (v x' \t 1_H) \tn( x'' y )=\sum_{(x)} (\e(x')v ) \tn( x'' y )= v \tn {(xy)},$$
and
$$(1_H \tn y) ( v \tn x)=\sum_{(y)} y' \t v \tn y'' x. $$
Therefore, clearly:
\begin{Proposition}\label{qwer2} Let  $\H=(\d \colon I \to H, \rho)$ be a Hopf algebra crossed module. We have a pre-crossed module: $$\widehat{\BA}(\H)=\big(t \colon I^0\otimes_\rho H \to H, \tra, \tla\big),$$ of bare algebras, where:
\begin{itemize}
 \item $t( v \tn x)=\d(v)x$, where $v \in I^0$ and $x \in H$. Recall the target map $t\colon I\otimes_\rho H \to H$ of {Lemma} \ref{majmain}.
 \item $x \tra a=i_H(x)\, a$ and $a \tla x=a\, i_H(x)$,  where $a \in  I^0\otimes_\rho H$, $x \in H$, and $i_H\colon H \to  I \otimes_\rho H$ is the inclusion. 
\end{itemize}
By considering the reflection functor ${\cal R}$; see \ref{refl}, from the category of pre-crossed modules of bare algebras to the category of crossed modules of bare algebras, one therefore has a crossed module of bare algebras:
 $${\BA}(\H)={\cal R}\big(\widehat{\BA}(\H)\big )= \big(t \colon \underline{I^0\otimes_\rho H} \to H, \tra, \tla\big).$$
\end{Proposition}

Suppose that $\A=(\d\colon A \to B, {\tra,\tla})$ is a unital crossed module of bare algebras. Pass to the associated crossed module of Hopf algebras  $\U(\Lie(\A))=(\d\colon \U(\Lie(A)) \to \U(\Lie(B)), \r)$; \ref{mcu}. We therefore have a pre-crossed module of bare algebras $\big(t\colon {\U^0(\Lie(\A))}\tn_\rho \U(\Lie(B)) \to \U(\Lie(B)), \tra,\tla\big)$. 
\begin{Theorem}[bare counit]\label{fgfgfg}
 Let $\kappa_\A\colon  {\U^0(\Lie(A))}\tn_\rho \U(\Lie(B)) \to A$ be the linear  map such that (see \ref{mcu}): $$\kappa_\A(v \tn x)= \xi^\flat_A(v) \tla \xi_B(x) , \textrm{ for each } x \in \U(\Lie(B)) \textrm{ and each } v \in \U(\Lie(A)).$$
 Then $\kappa_\A$ is an algebra map, and, in addition, the pair: 
$$\hat{K}_\A=(\kappa_\A, \xi_B)\colon \widehat{\BA}(\U(\Lie(\A)))=\big(t\colon {\U^0(\Lie(A))}\tn_\rho \U(\Lie(B)) \to \U(\Lie(B)), \tra,\tla\big) \to (\d\colon A \to B, {\tra ,\tla})$$ 
is a map of pre-crossed modules of bare algebras $\hat{K}_\A \colon \widehat{\BA}(\U(\Lie(\A))) \to \A$. Moreover, since $\A$ is a unital crossed module of bare algebras, then $\hat{K}_\A$ descend to a crossed module map, the bare counit,  denoted as: $$K_\A \colon \BA\big(\U(\Lie(\A))\big) \to \A.$$
\end{Theorem}
\begin{Proof} 
 Let $u \tn x, v \tn y \in {\U^0(\Lie(\A))}\tn_\rho \U(\Lie(B))$. Suppose that $x \in B$, thus $\D(x)=x \tn 1 + 1 \tn x$. Then:
\begin{align*}\allowdisplaybreaks
 \kappa_\A\big ( (u \tn x)(v \tn y) \big) &=\sum_{(x)} \kappa_\A\big ( (u\,   x' \t v) \tn (x'' y) \big)=\sum_{(x)} \xi^\flat_A( u\,   x' \t v) \tla \xi_B(x'' y)\\ &=\sum_{(x)} \big ( \xi^\flat_A( u) \, \xi^\flat_A(  x' \t v)\big) \tla\big( \xi_B(x'') \xi_B( y)\big)\\
&=\xi^\flat_A(u) \xi^\flat_A (x \t v) \tla \xi_B(y) +  \xi^\flat_A(u) \xi^\flat_B (v) \tla \xi_B(x)\xi_B(y)\\
&=\xi^\flat_A(u) \big (\xi_B(x) \tra\xi^\flat_A ( v) -\xi^\flat_A ( v)  \tla \xi_B(x)\big)  \tla \xi_B(y) +  \xi^\flat_A(u) \xi^\flat_A (v) \tla \xi_B(x)\xi_B(y)\\
&=\big(\xi^\flat_A(u) \tla \xi_B(x)\big)\, \big(\xi^\flat_A(v) \tla \xi_B(y)\big)=\kappa_\A\big ( u \tn x		\big)\, \kappa_\A \big(v \tn y\big) ,  
\end{align*}
where we have used lemma \ref{www} and the fact that the actions $\tla$ and $\tra$ are fully interchangeable. Suppose, by induction, that we have proved that: $$\kappa_\A\big ( (u \tn x)(v \tn y) \big)=\big(\xi^\flat_A(u) \tla \xi_B(x)\big)\, \big(\xi^\flat_A(v) \tla \xi_B(y)\big)=\kappa_\A\big ( u \tn x		\big)\, \kappa_\A \big(v \tn y\big),$$ in general, for all $x$ of filtered degree $\leq n$ in $\U(\Lie(B))$. If $b \in B$  we have, since $\Delta(b)=b\tn 1 +1 \tn b$:
\begin{align*}\allowdisplaybreaks
 \kappa_\A&\big ( (u \tn xb) (v \tn y) \big)=\kappa _\A\big ( (u \tn x) \, (1 \tn b)\, (v \tn y) \big)= \kappa _\A\big ( (u \tn x) \, ( b \t v \tn y)+ (u \tn x) \, ( v \tn by)    \big )\\
&=\kappa _\A(u \tn x)  \, \kappa_\A( b \t v \tn y)+ \kappa(u \tn x) \, \kappa( v \tn by)  \\
&=\big(\xi^\flat_A(u) \tla \xi_B (x) \big)\, \big(\xi^\flat_A(b \t v) \tla \xi_B (y) \big) + \big(\xi^\flat_A(u) \tla \xi_B (x) \big)\, \big(\xi^\flat_A( v) \tla \xi_B (by) \big)\\
&=\big(\xi^\flat_A(u) \tla \xi_B (x) \big)\, \Big( \xi_B(b) \tra \xi^\flat_B(v) - \xi^\flat_A(v) \tla \xi_B(b) \Big) \tla\xi_B (y)  + \big(\xi^\flat_A(u) \tla \xi_B (x) \big)\, \big(\xi^\flat_A( v) \tla \xi_B (by) \big)\\
&=\big(\xi^\flat_A(u) \tla \xi_B (x) \big)\, \xi_B(b) \tra \big(\xi^\flat_A( v) \tla \xi_B (y) \big) =\big(\xi^\flat_A(u) \tla \xi_B (x) \big)\tla\xi_B(b) \, \big(\xi^\flat_A( v) \tla \xi_B (y) \big) \\
&=\big(\xi^\flat_A(u) \tla \xi_B (xb) \big)\,  \big(\xi^\flat_A( v) \tla \xi_B (y) \big)=\kappa_\A\big ( (u \tn xb) \big) \, \kappa _\A\big((v \tn y) \big).
\end{align*}
Therefore $\kappa _\A$ is an algebra morphism. Now note that, if $y \in B$ and $v \tn x \in {\U^0(\Lie(A))}\tn_\rho \U(\Lie(B))$:
\begin{align*}\allowdisplaybreaks
\kappa_\A\big( (v \tn x) \tla y\big)&\doteq\kappa_\A( (v \tn x) (1 \tn y) \big) = \kappa _\A\big( v \tn (xy)\big)   =\xi^\flat_A(v) \tla \xi_B(xy)\\&=\big(\xi^\flat_A(v) \tla \xi_B(x)\big) \tla \xi(y)=\kappa _\A(v  \tn x) \tla \xi_B(y).
\end{align*}
And, since $\Delta(y)=y \tn 1 + 1 \tn y$:
\begin{align*}\allowdisplaybreaks
\kappa_\A\big( y \tra (v \tn x) \big)&=\kappa_\A(  (1 \tn y) (v \tn x) \big)= \kappa _\A\big (y \t v \tn x + v \tn yx \big)= \xi^\flat_A ( y \t v) \tla \xi_B(x)+ \xi^\flat_A(v) \tn \xi_B(yx)\\
 &= \xi(y) \tra \xi^\flat_A ( v) \tla \xi_B(x)-\xi^\flat_A(v) \tla \xi_B(y) \xi_B(x)+ \xi^\flat_A(v) \tn \xi_B(yx)\\&=\xi_B(y) \tra \xi^\flat_A ( v) \tla \xi_B(x)=\xi_B (y) \tra \kappa _\A(v \tn x). 
\end{align*}
Since $\U(\Lie(B))$ is generated, as an algebra, by $B$,  follows that $\hat{K}_\A=(\kappa_\A,\xi_\A)$ preserves the actions $\tla$ and $\tra$. The rest of the assertions follow immediately. 
\end{Proof}

We have an adjunction between the category of crossed modules of bare algebras and the category of crossed modules of Lie algebras; diagram \eqref{bareco}, where ${\rm Lie}$ is the right adjoint. The bare counit $K$ is its counit.

\subsubsection{The crossed module inclusion}\label{cminc}
   Let $\H=(\d \colon I \to H, \rho)$ be a crossed module of Hopf algebras. Consider the crossed modules of groups $\H^*_{\rm gl}$ and $\BA(\H)^*_\bullet$; see \ref{bag}, \ref{hag} and \ref{qwer}. We have a crossed module map, called the crossed module inclusion: 
   $${\rm Inc}_\H=(J_\H,\iota_H) \colon  (\H)^*_{{\rm gl}} \to  \big({\cal BA}(\H)\big)^*_\bullet.$$ 
   The map $\iota_\H\colon H^*_{\rm gl} \to H^*_\bullet$ is simply  group inclusion. Noting that $\widehat{{\cal BA}}(\H)=\left (t\colon {I^0 \tn_{\rho} H}  \to H, \tra, \tla \right )$, where ${I^0 \tn_{\rho} H}$ is a subalgebra of the Hopf algebra ${I \tn_{\rho} H}$,
   being that $I$ is a Hopf subalgebra of  ${I \tn_{\rho} H}$, follows that we have  group inclusions:
  $$I^*_{\rm gl} \to {\left (I \tn_{\rho} H\right )}^*_{\rm gl} \to \left (I \tn_{\rho} H\right )^*= \left (I^0 \tn_{\rho} H\right )^*_\bullet. $$ 
  The group map $J_\H\colon I^*_{\rm gl} \to \underline{\left(I^0 \tn_{\rho} H\right )}^*_\bullet$, is obtained by composing with the map $\left (I^0 \tn_{\rho} H\right )^*_\bullet \to \underline{\left (I^0 \tn_{\rho} H\right )}_\bullet^*$, induced by the algebra projection  $I^0 \tn_{\rho} H \to  \underline{I^0 \tn_{\rho} H}$. Clearly ${\rm Inc}_\H \colon  (\H)^*_{{\rm gl}} \to  \big ({\cal BA}(\H)\big)^*_\bullet$ is a map of crossed modules of groups.

The crossed module counit  ${\rm Proj}_\A\colon \U(\Lie(\A))^*_{\rm gl} \to \A^*_\bullet$ (Theorem \ref{zxzx}) can be obtained from the
bare counit $K_\A \colon \BA(\U(\Lie(\A))) \to \A$ and the crossed module inclusion as the composition; c.f. \eqref{bareco}, \eqref{cmcounit}  and \eqref{IOTA}:
$$\U(\Lie(\A))^*_{\rm gl} \ra{{\rm Inc}_{\U(\Lie(\A))} }  \big({\cal BA}(\U(\Lie(\A)))\big)^*_\bullet \ra{ (K_A)^*_\bullet} \A^*_\bullet.$$

\section{Holonomy over Hopf and bare algebras}

\subsection{Some facts on iterated integrals}

\subsubsection{Differentiating and integrating over general algebras}
We will consider algebras $H$ without any topology in  it, justifying the need for the discussion now. Recall that all normed spaces of finite dimension are equivalent to some $\mathbb{R}^n$, and any linear map $\mathbb{R}^n \to \mathbb{R}^m$ is $C^{\infty}$. 

\begin{Definition}[Smooth map, partial differentiation and integration]
Let $H$ be an algebra over $\C$. Given an open or closed set  (usually a simplex or a parallelepiped)  $U \subset 	\mathbb{R}^n$, a map $f\colon U \to H$ is called smooth if:
\begin{enumerate}
 \item The range of $f$ is contained in a finite dimensional subspace $V$ of $H$.
 \item The induced map $f\colon U \to V$ is smooth (in other words $C^{\infty}$). 
\end{enumerate}
Let $C^{\infty}(U,H)$ be the vector space of smooth maps $f\colon U \to H$. Given that the range of each $f \in C^\infty(U,H)$ is contained in a finite dimensional vector space $V \subset H$, we can consider the  partial differentiation operators $\frac{\partial}{\partial x_i} \colon C^{\infty}(U,H)\to C^{\infty}(U,H),$ as well as, whenever $U$ is compact, the integration $\displaystyle{{\idotsint\limits_U} f(x)d^nx\in V \subset H }$.

\noindent Note that products and sums of smooth functions $U \to H$ are smooth. 

\end{Definition}
Given that linear maps between finite dimensional vector spaces are always smooth it follows that:
\begin{Lemma}\label{fd}
 If $f\colon U \to H$ is smooth and $T\colon H \to H'$ is a linear map then $T\circ f\colon U \to H'$ is smooth. Moreover $T$ commutes with partial differentiation and integration.  
\end{Lemma}

\begin{Lemma}[Uniqueness of integration]\label{ui}
 Let $H$ be an algebra over $\C$. Let $f\colon [a,b] \to H$ be a smooth map. Given $g \in H$ there exists a unique smooth function $F\colon [a,b] \to H$ such that $F'(x)=f(x), \forall x \in [a,b]$,  and with $f(a)=g$. In fact $F(t)=g+\int_a^t f(u) du, \forall t \in  [a,b]$. 
\end{Lemma}
\begin{Proof} Let $W\subset H $ be a finite dimensional vector space containing both $f([a,b])$ and $\{g\}$. Then use uniqueness of integration for smooth functions with values in finite dimensional vector spaces. 
\end{Proof}

\subsubsection{Definition of iterated integrals }
We present a small part of the  theory of iterated integrals appearing in \cite{Chen}. Let $H$ be an algebra. Given smooth functions $u \in [0,1]  \mapsto f_i(u) \in H$, where $i=1,2,\dots, n$, we define, given $t_0,t \in [0,1]$, with $t_0\leq t$:
$$\oint_{t_0}^t f_1(u)du =\int_{t_0}^t f_1(u_1) du_1.  $$
In general,   we put:
\begin{align*}
\oint_{t_0}^t f_1(u)*\dots* f_n (u)\, du&=
\int_{t_0}^t \int_{t_0}^{u_1} \int_{t_0}^{u_2}\int_{t_0}^{u_3} \dots \int_{t_0}^{u_{n-1}}  f_1(u_1)\,  f_2(u_2)\,f_3(u_3)\dots f_{n}(u_n)\, du_{n}\, \dots \, du_3 \, d u_2\, d u_1\\
&=\idotsint\limits_{\{t\ge u_1 \ge u_2 \ge \dots \ge u_n \ge t_0 \}}  f_1(u_1)\,  f_2(u_2)\,f_3(u_3)\dots f_{n}(u_n)\, du_{n}\, \dots \, du_3 \, d u_2\, d u_1 \,. 	
 \end{align*}
Clearly:
\begin{equation}
\frac{ d }{d t}\oint_{t_0}^t f_1(u)*f_2(u)\dots* f_n (u)  \,d u=f_1(t) \oint_{t_0}^t  f_2(u)\dots* f_n (u) \, d u.
\end{equation}
We will also consider the following alternative presentation of iterated integrals, which, for reasons which will be clearer later, we prefer to denote differently.
Given $t_1, t \in [0,1]$, with $t\leq t_1$, we put:
$$\sqint_{t}^{t_1} f_1(u)du =\int_{t}^{t_1} f_1(u_1) du_1, $$
and, in general:
\begin{align*}
\sqint_{t}^{t_1} f_1(u)*\dots* f_n (u) du&=  
\int_{t}^{t_1} \int^{t_1}_{u_1}\int^{t_1}_{u_2} \int^{t_1}_{u_3} \dots \int^{t_1}_{u_{n-1}}  f_1(u_n)   \, f_2(u_{(n-1)})\, f_3(u_{(n-2)}) \dots f_{n}(u_1)   \, du_{n}\, \dots \, du_3 \, d u_2\, d u_1\\
&= \idotsint\limits_{\{t \leq u_1 \leq u_2 \leq \dots \leq u_n\leq t_1\}}   f_1(u_n)   \, f_2(u_{(n-1)})\, f_3(u_{(n-2)}) \dots f_{n}(u_1)   \, du_{n}\, \dots \, du_3 \, d u_2\, d u_1\,.
 \end{align*}
We have:
$$\frac{d}{dt} \sqint_{t}^{t_1} f_1(u)*\dots* f_n (u) du =-\left (\sqint_{t}^{t_1} f_1(u)*\dots* f_{(n-1)} (u) du  \right) f_n(t)  .$$
Clearly (by inverting the order of integration, and then renaming the variables.)
\begin{Lemma}\label{sqcir}
 For all $0<a<b<1$ we have:
$$\oint_{a}^b f_1(u)*\dots* f_n (u) du=\sqint_{a}^{b} f_1(u)*\dots* f_n (u) du.$$
\end{Lemma}

\subsection{The holonomy of algebra valued connections}

\subsubsection{1-paths, 2-paths and 3-paths. Thin 2-paths.}\label{123paths}

Nothing here is new; see for example \cite{CP,SW1,SW2,FMP1,FMP3}. 
Let $M$ be a smooth manifold. We denote by $\P^1(M)$ the space of continuous and piecewise smooth {maps} $\g\colon [0,1] \to M$, called 1-paths, and normally represented in the form below, emphasizing initial and end points:
$$x=\g(0) \ra{\g} \g(1)=y. $$
If $\g_1$ and $\g_2$ are 1-paths with $\g_1(1)=\g_2(0)$, we define their concatenation, written as: $$\left (x \ra{\g_1} y \ra{\g_2} z\right)=x\ra{\g_1 \g_2 } z,$$ where $x=\g_1(0)$, $y=\g_1(1)=\g_2(0)$ and $z=\g_2(1)$, in the usual way:
$$(\g_1\g_2)(t)=\left\{\begin{CD} \g_1(2t), \quad t \in [0,1/2]  \\
                                 \g_2(2t-1), \quad t \in [1/2,1]   
                                \end{CD}\right. $$
The reversion of the {1-path} $x \ra{\g} y$ is the {1-path} $y \ra{\g^{-1}} x$, where $\g^{-1}(t)=\g(1-t),\forall t \in [0,1]$. 
The composition of 1-paths is not associative. This will not pose any problem for this paper.  This can be fixed, for example, by considering a thin homotopy relation between 1-paths; see \cite{CP} and below.

We analogously define a 2-path as being a map $\G\colon [0,1]^2 \to M$, which is to  be continuous and piecewise smooth for some paving of the square by polygons, transverse to the boundary of the square. Moreover, we suppose that $\G(\{0,1\} \times [0,1])$ has at most two elements. Putting $\{x\}=\G(\{0\} \times [0,1])$, $\{y\}=\G(\{1\} \times [0,1])$ (the left and right boundaries of $\G$), and also $\g_1(t)=\G(t,1)$, $\g_0(t)=\G(t,0)$, where $t \in [0,1]$ (the bottom and top boundaries of $\G$), defining {1-path}s $\g_0,\g_1\colon [0,1] \to M$, we denote $\G$ in the form $\G=\morts{x}{\g_0}{\G}{\g_1}{y} $, emphasizing its boundaries. The set of 2-paths in $M$ is denoted by $\P^2(M)$.

Clearly 2-paths compose horizontally and vertically in the obvious way, whenever pertinent boundaries coincide. We also have whiskerings of 2-paths by 1-paths: consider a 2-path $\G$, and 1-paths $\g_1$ and $\g_2$, such that the left boundary of $\G$ is $\g_1(1)$ and the right boundary of $\G$ is $\g_2(0)$; see below. We put:
$$ \hskip-1cm{\xymatrix{ & {x} \ar@/^1pc/[rr]^{{\phi_0} }\ar@/_1pc/[rr]_{{\phi_1}} & \Uparrow {\G} & {y} \ar[r]^{{\g_2}} &{z} }} =\hskip-1cm\mort{x}{\phi_1\g_2}{\G\g_2}{\phi_1\g_2}{z}, \we (\G \g_0)(t,s)= \begin{cases}\G(2t,s), \textrm{ if }0\leq t \leq 1/2,  \\ \g_2(2t-1) \textrm{ if } 1/2 \leq t \leq 1 \end{cases},	 $$
$$\hskip-1cm{\xymatrix{ & {w} \ar[r]^{{\g_1}}  & x \ar@/^1pc/[rr]^{{\phi_1}} \ar@/_1pc/[rr]_{{\phi_0}} & \G & {y}  } }  =\hskip-1cm\mort{w}{\g_1\phi_0}{\g_1\G} {\g_1\phi_1}{y}, \we (\g_1 \G )(t,s)= \begin{cases}\g_1(2t), \textrm{ if }0\leq t \leq 1/2,  \\ \G(2t-1,s) \textrm{ if } 1/2 \leq t \leq 1 \end{cases} . $$
These compositions are not associative and do not satisfy the interchange law. This can be fixed, in order to define a 2-category, by considering thin homotopy relations between 2-paths; see \cite{FMP1,FMP2,SW1}.
\begin{Definition}[Thin 2-path. Thin homotopy]
 A 2-path $\hskip-1cm {\xymatrix{&x\ar@/^0.5pc/[r]^{\g_2}\ar@/_0.5pc/[r]_{\g_2}\ar@{{}{ }{}} [r]|{\Uparrow\G} &y}}$
 is called thin if, at each $(s,t) \in  [0,1]^2$, where $\G\colon [0,1]^2 \to M$ is smooth, the rank of the derivative of $\G$ is $\leq 1$. Such a 2-path $\G$ is called a thin homotopy from $\g_0$ to $\g_1$, and the 1-paths $\g_0$ and $\g_1$ are said to be thin-homotopic. 
\end{Definition}
\begin{Remark}[The thin homotopy groupoid of a smooth manifold]\label{tth2g}
Thin homotopy is an equivalence relation, the composition of 1-paths descends to the quotient, and we have a groupoid, the thin homotopy groupoid of $M$, whose set of objects is $M$ and whose morphisms are {1-path}s up to thin-homotopy; see \cite{FMP1,CP}. 
\end{Remark}

\begin{Definition}[3-path. Thin 3-path.]\label{3path}
Consider two 2-paths \smash{$\morts{x}{\g_0}{\G}{\g_1}{y}$ and $\morts{x}{\phi_0}{\Phi}{\phi_1}{y}.$}
 A 3-path from the first to the second  is given by a map $J\colon [0,1]^3 \to M$, such that:
\begin{itemize}
 \item The map $J$ is piecewise smooth for some paving of the cube $[0,1]^3$ by polygons, which is transverse to the boundary of the cube.
 \item $J(\{0\}\times [0,1]^2)=\{x\}$ and $J(\{1\}\times [0,1]^2)=\{y\}.$
 \item $J(t,s,0)=\G(t,s)$ and $J(t,s,1)=\Phi(t,s)$, for all $t,s\in [0,1]$. 
\end{itemize}
Each 3-path $J$ has top and bottom 2-paths $J^t$ and $J^b$, where $J^t(t,x)=J(t,1,x)$ and $J^b(t,x)=J(t,0,x)$, where $t,x \in [0,1]$. We say that $J$ is thin if both $J^t$ and $J^b$ are thin 2-paths, and also, given any $(t,s,x)\in [0,1]^3$ where $J\colon [0,1]^3\to M$ is smooth, the derivative of $J$ at $(t,s,x)$ has rank $\leq 2$.
\end{Definition}
\begin{Remark}[The thin homotopy 2-groupoid of a smooth manifold]
Thin homotopy of 2-paths is an equivalence relation. All operations already defined for 1-paths and 2-paths (whiskerings and horizontal/vertical compositions) descend to the quotients by thin homotopy. Moreover, we have a 2-groupoid, the thin homotopy 2-groupoid of $M$, with objects being the points of $M$, with morphisms being the {1-path}s connecting points, up to thin homotopy, and 2-morphisms being 2-paths, up to thin homotopy. We note that this result is not immediate and the details appear {spelled} out (using two different arguments) in \cite{FMP1} and \cite{FMP2}. 
\end{Remark}
\subsubsection{The  holonomy of a 1-connection}\label{h1c}
\begin{Definition}[Algebra-valued smooth form on a manifold] Let $M$ be a smooth manifold.
 Let $H$ be  a bare algebra. A smooth $n$-form $\w$ is an $n$-form in $M$, with values in $H$, whose range is included in a finite dimensional subspace $V \subset H$, being that $\w$ is smooth when considered to have values in $V$.   
 \end{Definition}

\noindent We consider holonomy operators with values in $H\h$, where $H$ is an algebra. This is due to Chen \cite{Chen}.
\begin{Definition}
 Let $U$ be an open or closed set of $\mathbb{R}^n$. A function $ x\in U \mapsto \sum_{n=1}^{+\infty} f_n(x) h^n \in H\h$ is said to be smooth if each function $f_n\colon U \to H$ is smooth.  Differentiation and integration of smooth functions will be done term-wise, in the obvious way.
\end{Definition}
\begin{Lemma}\label{ivp}
 Let $H$ be an algebra. Let $x \in H\h$. Let $f\colon [0,1] \to H\h$ be a smooth function. The initial value problems below have unique smooth solutions, $g_1,g_2\colon [0,1] \to H\h$:
 \begin{align*}
 & g_1'(t)=h\, g_1(t)\, f(t); & \textrm{with } g_1(0)=x , \\  & g_2'(t)=h \, f(t) \, g_2(t); &  \textrm{with } g_2(0)=x  .
 \end{align*}
\end{Lemma}
\begin{Proof}
 Follows directly by induction and by uniqueness  of integration; Lemma \ref{ui}. Recall that  products and sums of smooth functions $[0,1]\to H$ are  smooth. 
\end{Proof}

\begin{Definition}[(local) connection and curvature. Primitive connection]
 Let $H$ be a {unital algebra}, with unit $1_H$.  A (local) connection, or local connection form,  in $M$ is a smooth 1-form $\w$ in $M$ with values in $H$. If $H$ is a Hopf algebra, a  connection $\w$ is called primitive if it takes values inside the primitive space  $\Prim(H)$ of $H$. Given a connection $\w$ its curvature (note our conventions) is the $H$-valued smooth  2-form:
 $$\F_\w=d \w -\frac{1}{2} \w \wedge \w .$$
 Here, if $X$ and $Y$ are vector fields in $M$, we put
 $$ \frac{1}{2}(\w \wedge^{[,]} \w) (X,Y)= \frac{1}{2}(\w \wedge \w) (X,Y)=[\w(X),\w(Y)].$$ 
 From a connection  form $\w$ we can build the graded connection form $\bw\doteq h\w$, and   the  graded curvature 2-form: $$\F_{\bw}\doteq h d \w -h^2\frac{1}{2} \w \wedge \w= d \bw-\frac{1}{2}\bw\wedge \bw. $$
\end{Definition}
\noindent Let $\g\colon [0,1] \to M$ be a piecewise smooth curve in $M,$ a smooth manifold. Given an $H$-valued connection form in $M$, where $H$ is a {unital algebra}, consider the following $H[[h]]$ valued function  in $[0,1]$, where $\g'$ denotes the derivative of $\g$. Here $t,t_0 \in [0,1]$:
\begin{align*}\allowdisplaybreaks
{P_{\bw}}(\g,[t, t_0])&=\sum_{n=0}^{+ \infty} h^n \oint_{t_0}^t \underbrace{\w(\g'(u))*\dots*\w(\g'(u))}_{n \textrm{ terms}} du=\sum_{n=0}^{+ \infty} \oint_{t_0}^t \underbrace{\bw(\g'(u))*\dots*\bw(\g'(u))}_{n \textrm{ terms}} du.
\end{align*}
Here we have conventioned that $\oint_{t_0}^t \underbrace{\w(\g'(u))*\dots*\w(\g'(u))}_{0 \textrm{ terms}} du=1_H$, where $t,t_0 \in [0,1]$. From {Lemma} \ref{ui}:
\begin{Proposition}\label{uii}
The function $t \mapsto {P_{\bw}}(\g,[t,t_0])$ is the unique smooth solution of the differential equation:
$$f'(t)= \bw\big(\g'(t)\big)\, f(t), \quad f(t_0)=1_H. $$ 
\end{Proposition}
In general we have, for $t,t_0 \in [0,1]$ (note {Lemma} \ref{sqcir}):
\begin{align}\label{h1d}
  &\dt {P_{\bw}}(\g,[t,t_0])= \bw\big(\g'(t)\big)\,{P_{\bw}}(\g,[t,t_0]),
  &\frac{\d} {\d t_0} {P_{\bw}}(\g,[t,t_0])= -{P_{\bw}}(\g,[t,t_0])\,\, \bw\big(\g'(t_0)\big).
  \end{align}

Standard iterated integral theory (or simple calculations) yields.
\begin{Proposition}
 Let $\g_1$ and $\g_2$ be piecewise smooth curves in $M$, such that $\g_1(1)=\g_2(0)$. Let $\g=\g_1 \g_2$ be the concatenation of $\g_1$ and $\g_2$. Then, for each $t \in [1/2,1]$, we have:
$${P_{\bw}}(\g, [t,0])={P_{\bw}}(\g_2,[2t-1,0])\,{P_{\bw}}(\g_1,[1,0]).  $$
\end{Proposition}
\begin{Proposition}
 Let $\g$ be a piecewise smooth  curve in $M$. Then, if $0 \leq t_0 \leq s \leq t_1 \leq 1$, we have:
\begin{equation}\label{compr} {P_{\bw}}(\g,[t_1,t_0])={P_{\bw}}(\g,[t_1,s])\,{P_{\bw}}(\g,[s,t_0]) . \end{equation}
Furthermore, this is valid for arbitrary $t_0,t_1, s\in [0,1]$. 
\end{Proposition}

The following is well known.
\begin{Lemma}\label{glp} Suppose that $H$ is an Hopf algebra and that $\w$ is a primitive connection. Then $t \mapsto {P_{\bw}}(\g,[t,t_0])$ takes values in $H^*_{\rm gl}$ (the space of grouplike elements of $H[[h]]$ whose power series starts at $1_H$; see \ref{hag}.)
\end{Lemma}
\begin{Proof}
The map $t \mapsto \Delta({P_{\bw}}(\g,[t,t_0]))\in (H\tn H)\h$ is smooth. Moreover $\Delta({P_{\bw}}(\g,[t_0,t_0]))=\D(1_H)=1_H \tn 1_H. $ 
Also:
\begin{align*}\allowdisplaybreaks
 \frac{d}{dt}\Delta\big({P_{\bw}}(\g,[t,t_0])\big)&=\Delta\left(  \frac{d}{dt} {P_{\bw}}(\g,[t,t_0]) \right)  =\Delta\big( \bw(\g'(t)\,\,   {P_{\bw}}(\g,[t,t_0])\big)\\&=\Delta\big(\bw(\g'(t)) \big)\,\,\Delta\big(  {P_{\bw}}(\g,[t,t_0])\big) =\big(\bw(\g'(t)) \tn 1_H + 1_H \tn \bw(\g'(t))   \big)\,\, \Delta\big(  {P_{\bw}}(\g,[t,t_0])\big). 
\end{align*}
The map $t\mapsto {P_{\bw}}(\g,[t,t_0])\tn   {P_{\bw}}(\g,[t,t_0])$ is smooth, and ${P_{\bw}}(\g,[t_0,t_0])\tn   {P_{\bw}}(\g,[t_0,t_0])=1_H \tn 1_H$. Moreover:
\begin{align*}\allowdisplaybreaks
 \frac{d}{dt} \big( {P_{\bw}}(\g,[t,t_0])\tn   {P_{\bw}}(\g,[t,t_0]) \big)&=   \Big( \frac{d}{dt} {P_{\bw}}(\g,[t,t_0])\Big) \tn   {P_{\bw}}(\g,[t,t_0]) +  {P_{\bw}}(\g,[t,t_0])\tn  \Big( \frac{d}{dt} {P_{\bw}}(\g,[t,t_0])\Big)\\
 &=   \Big( \w(\g'(t) )\ {P_{\bw}}(\g,[t,t_0])\Big) \tn   {P_{\bw}}(\g,[t,t_0]) +  {P_{\bw}}(\g,[t,t_0])\tn  \Big( \w(\g'(t) )\ {P_{\bw}}(\g,[t,t_0])\Big)\\
 &= \big(\w(\g'(t) )\tn 1_H + 1_H \tn \w(\g'(t))   \big)\,\, \big( {P_{\bw}}(\g,[t,t_0])\tn   {P_{\bw}}(\g,[t,t_0]) \big).
\end{align*}
Now apply {Lemma} \ref{ivp} to conclude that $\Delta\big({P_{\bw}}(\g,[t,t_0])\big)= {P_{\bw}}(\g,[t,t_0])\tn   {P_{\bw}}(\g,[t,t_0]),\forall t \in [t_0,1].	$
\end{Proof}

The following frequently used compatibility between holonomy and algebra maps has a trivial proof. 
\begin{Lemma}\label{opopo}
 Let $f\colon H \to H'$ be an algebra map. Let $\w\in \Omega^1(M,H)$. We have, for each $t_0,t_1 \in [0,1]$.
  $$f\left({P_{\bw}}(\g,[t_1,t_0])\right) = {P_{f(\bw)}}(\g,[t_1,t_0]).$$
\end{Lemma}

\subsubsection{The curvature and holonomy}\label{ch}
Given a manifold $M$, let now $\G\colon [0,1]^2 \to M$ be a 2-path. Therefore $s \mapsto \G(0,s)$ and $s \mapsto \G(1,s)$ are constant {1-path}s.  Define, given $(t,s) \in [0,1] \times [0,1]$:
$${P_{\bw}}(\G,[t,0],s)= \sum_{n=0}^{+ \infty} h^n \oint_0^t \underbrace{\w\left( \frac{\d}{\d u} \G(u,s)  \right) *\dots*  \w\left( \frac{\d}{\d u} \G(u,s)  \right) }_{n \textrm{ terms}} du.$$
If $s \in [0,1]$, considering the 1-path $\g_s$ with $\g_s(t)=\G(s,t), \forall s,t \in [0,1]$, we have ${P_{\bw}}(\G,[t,0],s)={P_{\bw}}(\g_s,[t,0])$.

For each  $s,t \in [0,1]$ we have:\small
\begin{align*}\allowdisplaybreaks
&\frac{\d}{\d s} {P_{\bw}}(\G,[t,0],s)=\sum_{n=0}^{+ \infty} h^n \\&\sum_{i+j+1=n}\oint_0^t \underbrace{\w\left( \frac{\d}{\d u} \G(u,s)  \right) *\dots* \w\left( \frac{\d}{\d u} \G(u,s)  \right) }_{i \textrm{ terms}  } *\ds  \w\left( \frac{\d}{\d u} \G(u,s)  \right) *\underbrace{\w\left( \frac{\d}{\d u} \G(u,s)  \right) *\dots* \w\left( \frac{\d}{\d u} \G(u,s)  \right) }_{j \textrm{ terms}  }  d u.
\end{align*}
\normalsize
Note that: $$\ds  \w\left( \frac{\d}{\d u} \G(u,s)\right ) =d \w \left ( \frac{\d}{\d s} \G(u,s), \frac{\d}{\d u} \G(u,s) \right )+ \du \w \left (\frac{\d}{\d s} \G(u,s) \right).$$
Therefore (integrating by parts, and noting that $\frac{\d}{\d s} \G(0,s)=0$):
\small
\begin{align*}\allowdisplaybreaks
&\frac{\d}{\d s} {P_{\bw}}(\G,[t,0],s)=\sum_{n=0}^{+ \infty} h^n\sum_{i+j+1=n} \Big ( \\&\oint_0^t \underbrace{\w\left( \frac{\d}{\d u} \G(u,s)  \right) *\dots* \w\left( \frac{\d}{\d u} \G(u,s)  \right) }_{i \textrm{ terms}  } * \,d \w \left ( \frac{\d}{\d s} \G(u,s), \frac{\d}{\d u} \G(u,s) \right )  *\underbrace{\w\left( \frac{\d}{\d u} \G(u,s)  \right) *\dots* \w\left( \frac{\d}{\d u} \G(u,s)  \right) }_{j \textrm{ terms}  }du  \\
&-\oint_0^t \underbrace{\w\left( \frac{\d}{\d u} \G(u,s)  \right) *\dots* \w\left( \frac{\d}{\d u} \G(u,s)  \right) }_{i \textrm{ terms}  } *    \underbrace{\left (\w\left( \frac{\d}{\d s} \G(u,s)\right) \w \left (\frac{\d}{\d u} \G(u,s)  \right) \right)  * \w \left (\frac{\d}{\d u} \G(u,s) \right)*\dots*  \w\left( \frac{\d}{\d u} \G(u,s)  \right) }_{j \textrm{ terms}  } du\\
&+\oint_0^t \underbrace{\w\left( \frac{\d}{\d u} \G(u,s)  \right) *\dots* \w\left( \frac{\d}{\d u} \G(u,s)  \right)*\left (\w\left( \frac{\d}{\d u} \G(u,s)\right) \w \left (\frac{\d}{\d s} \G(u,s)  \right) \right)   }_{i \textrm{ terms}  } *    \underbrace{  \w \left (\frac{\d}{\d u} \G(u,s) \right) * \dots *\w\left( \frac{\d}{\d u} \G(u,s)  \right) }_{j \textrm{ terms}  } du \Big)\\
&+\sum_{n=0}^{+ \infty} h^n \w \left(\frac{\d}{\d s} \G(t,s)\right) \oint_0^t \underbrace{\w\left( \frac{\d}{\d u} \G(u,s)  \right) *\dots*  \w\left( \frac{\d}{\d u} \G(u,s)  \right) }_{(n-1) \textrm{ terms}} du.
\end{align*}
\normalsize
Thus:
\small
\begin{align*}\allowdisplaybreaks
&\frac{\d}{\d s} {P_{\bw}}(\G,[t,0],s)=\sum_{n=0}^{+ \infty} h^n\sum_{i+j+1=n}  \\&\oint_0^t \underbrace{\w\left( \frac{\d}{\d u} \G(u,s)  \right) *\dots* \w\left( \frac{\d}{\d u} \G(u,s)  \right) }_{i \textrm{ terms}  } * \,d \w \left ( \frac{\d}{\d s} \G(u,s), \frac{\d}{\d u} \G(u,s) \right )  *\underbrace{\w\left( \frac{\d}{\d u} \G(u,s)  \right) *\dots* \w\left( \frac{\d}{\d u} \G(u,s)  \right) }_{j \textrm{ terms}  }du  \\
&-
\sum_{n=0}^{+ \infty} h^n\sum_{p+q+2=n} \\&\oint_0^t \underbrace{\w\left( \frac{\d}{\d u} \G(u,s)  \right) *\dots* \w\left( \frac{\d}{\d u} \G(u,s)  \right) }_{p \textrm{ terms}  } * \frac{1}{2} \w \wedge \w \left ( \frac{\d}{\d s} \G(u,s), \frac{\d}{\d u} \G(u,s) \right )  *\underbrace{\w\left( \frac{\d}{\d u} \G(u,s)  \right) *\dots* \w\left( \frac{\d}{\d u} \G(u,s)  \right) }_{q \textrm{ terms}  }du \\
&+\sum_{n=0}^{+ \infty} h^n \w \left(\frac{\d}{\d s} \G(t,s)\right) \oint_0^t \underbrace{\w\left( \frac{\d}{\d u} \G(u,s)  \right) *\dots*  \w\left( \frac{\d}{\d u} \G(u,s)  \right) }_{(n-1) \textrm{ terms}} du\\
&=\sum_{i,j =0}^{+\infty}\Big (  \\&\oint_0^t \underbrace{\bw\left( \frac{\d}{\d u} \G(u,s)  \right) *\dots* \bw\left( \frac{\d}{\d u} \G(u,s)  \right) }_{i \textrm{ terms}  } * \,  \F_{\bw} \left ( \frac{\d}{\d s} \G(u,s), \frac{\d}{\d u} \G(u,s) \right )  *\underbrace{\bw\left( \frac{\d}{\d u} \G(u,s)  \right) *\dots* \bw\left( \frac{\d}{\d u} \G(u,s)  \right) }_{j \textrm{ terms}  }du\Big)\\
&+\sum_{n=0}^{+ \infty}  \bw \left(\frac{\d}{\d s} \G(t,s)\right) \oint_0^t \underbrace{\bw\left( \frac{\d}{\d u} \G(u,s)  \right) *\dots*  \bw\left( \frac{\d}{\d u} \G(u,s)  \right) }_{(n-1) \textrm{ terms}} du =A+B. 
\end{align*}
\normalsize
Now note that:
\small
\begin{align*}\allowdisplaybreaks 
A&=\sum_{i \in \N} \oint_0^t \underbrace{\bw\left( \frac{\d}{\d u} \G(u,s)  \right) *\dots* \bw\left( \frac{\d}{\d u} \G(u,s)  \right) }_{i \textrm{ terms}  } *\,\, \left (\F_{\bw} \left ( \frac{\d}{\d s} \G(u,s), \frac{\d}{\d u} \G(u,s) \right ) {P_{\bw}}(\G,u,s) \right ) du\\
&=\sum_{i \in \N} \sqint_0^t \underbrace{\bw\left( \frac{\d}{\d u} \G(u,s)  \right) *\dots* \bw\left( \frac{\d}{\d u} \G(u,s)  \right) }_{i \textrm{ terms}  } *\,\, \left (\F_{\bw} \left ( \frac{\d}{\d s} \G(u,s), \frac{\d}{\d u} \G(u,s) \right ) {P_{\bw}}(\G,u,s) \right ) du\\
&=\int_0^t {P_{\bw}}(\G,[t,u],s) \, \F_{\bw} \left ( \frac{\d}{\d s} \G(u,s), \frac{\d}{\d u} \G(u,s) \right )\, {P_{\bw}}(\G,[u,0],s) \,\, du. 
\end{align*}
\normalsize
The first and third equations are tautological. The second follows from Lemma \ref{sqcir}.
This well know result follows:
\begin{Lemma}\label{um}
In the conditions stated in the beginning of \ref{ch}, we have, for each $t,s \in  [0,1]$:
\begin{multline}\label{dsi} 
\frac{\d}{\d s} {P_{\bw}}(\G,[t,0],s)\\=\int_0^t {P_{\bw}}(\G,[t,u],s) \, \F_{\bw} \left ( \frac{\d}{\d s} \G(u,s), \frac{\d}{\d u} \G(u,s) \right ) {P_{\bw}}(\G,[u,0],s)  du+  \bw \left(\frac{\d}{\d s} \G(t,s)\right)  {P_{\bw}}(\G,[t,0],s). 
\end{multline}
In particular, since $\frac{\d}{\d s} \G(1,s)=0,\forall s \in [0,1]$, because $\G$ is a 2-path:
\begin{align}\allowdisplaybreaks \label{dsj}
\frac{d}{d s} {P_{\bw}}(\G,[1,0],s)&=\int_0^1 {P_{\bw}}(\G,[1,u],s) \, \F_{\bw} \left ( \frac{\d}{\d s} \G(u,s), \frac{\d}{\d u} \G(u,s) \right )\, {P_{\bw}}(\G,[u,0],s) \,du,
\end{align}
or, by \eqref{compr},
\begin{align}\allowdisplaybreaks\label{1dd}
\frac{d}{d s} {P_{\bw}}(\G,[1,0],s)&= \left(\int_0^1 {P_{\bw}}(\G,[1,u],s)   \F_{\bw} \left ( \frac{\d}{\d s} \G(u,s), \frac{\d}{\d u} \G(u,s) \right ) {P_{\bw}}(\G,[1,u],s)^{-1}  du\right){P_{\bw}}(\G,[1,0],s).
\end{align}
\noindent In particular it follows that thin homotopic 2-paths have the same holonomy, a fact observed in \cite{CP}. 
\end{Lemma}
Since $\frac{d}{d s} \left ( {P_{\bw}}(\G,[1,0],s) \,  {P_{\bw}}(\G,[1,0],s)^{-1} \right) =0$,   using Leibniz law, $\ds \big (a(s)\, b(s)\big)=\big ( \ds a(s)\big)\, b(s) + a(s) \big( \ds b(s)\big)$:
\begin{align}\allowdisplaybreaks\label{2dd}
\hskip-2mm\frac{d}{d s} {P_{\bw}}(\G,[0,1],s)&= -\left(\int_0^1\hskip0mm {P_{\bw}}(\G,[0,u],s)   \F_{\bw} \left ( \frac{\d}{\d s} \G(u,s), \frac{\d}{\d u} \G(u,s) \right )\hskip-1mm {P_{\bw}}(\G,[0,u],s)^{-1} du\right) \hskip-1mm {P_{\bw}}(\G,[0,1],s).
\end{align}
(Note  ${P_{\bw}}(\G,[1,0],s)^{-1}={P_{\bw}}(\G,[1,0],s)$.) Also, proceeding in the same way as before:
\begin{multline}\label{dsi2} 
\frac{\d}{\d s} {P_{\bw}}(\G,[0,t],s)=- \left ( \int_0^t {P_{\bw}}(\G,[0,u],s) \, \F_{\bw} \left ( \frac{\d}{\d s} \G(u,s), \frac{\d}{\d u} \G(u,s) \right )\, {P_{\bw}}(\G,[u,0],s)\right) {P_{\bw}}(\G,[0,t],s)\\ - {P_{\bw}}(\G,[0,t],s) \bw \left(\frac{\d}{\d s} \G(t,s)\right). 
\end{multline}
\begin{multline}\label{dsi3} 
\frac{\d}{\d s} {P_{\bw}}(\G,[t,1],s)=-\int_t^1 {P_{\bw}}(\G,[t,u],s) \, \F_{\bw} \left ( \frac{\d}{\d s} \G(u,s), \frac{\d}{\d u} \G(u,s) \right )\, {P_{\bw}}(\G,[u,1],s) \,\, du\\+  \bw \left(\frac{\d}{\d s} \G(t,s)\right)  {P_{\bw}}(\G,[t,1],s). 
\end{multline}
\section{The holonomy of  Hopf  and bare 2-connections}
\subsection{Local Hopf 2-connections and their 2-curvature}\label{lhopf}
Let $\H=(\d \colon I \to H, \rho)$ be a crossed module of Hopf algebras.   Let $M$ be a manifold. 
\begin{Definition}[Local Hopf 2-connection]
 A local Hopf 2-connection is given by a triple $(\w,m_1,m_2)$ of differential forms in $M$:
$$\w \in \Omega^1(M,H),  \quad\quad  \quad m_1,m_2 \in \Omega^2(M,I). $$
We impose the following compatibility relation (called the vanishing of the fake curvature; \cite{BS}):
$$\d (m_1)= d \w; \quad \d(m_2)=- \frac{1}{2} \w \wedge^{[,]} \w. $$
Let $h$ be a formal parameter. From the triple $(\w,m_1,m_2)$ we can build the pair $(\bw,\bm)$ of graded forms:
$$\bw=h\w; \quad \bm=h m_1 +h^2 m_2. $$
Thus $\bw$ takes values in $H[[h]]$, $\bm$ takes values in $I[[h]]$, and the vanishing of the fake curvature means that:
$$\d(\bm)=\F_{\bw} $$
where by {definition} $\F_{\bw} \in \Omega^2(M,H)$ is the graded curvature of $\bw$, thus
$\F_{\bw}=d \bw-\frac{1}{2} \bw \wedge^{[,]} \bw=h d\w-h^2 \frac{1}{2} \w \wedge^{[,]} \w. $
\end{Definition}
\noindent This {definition} makes sense for a crossed module of bare algebras (see \ref{bare}) or of Lie algebras (the latter being the usual setting \cite{SW1,SW2}). For the following definition we restrict ourselves to the Hopf algebra case. 
\begin{Definition}[Primitive Hopf 2-connection]\label{fully}
 A local Hopf 2-connection $(\w,m_1,m_2)$  is called:
\begin{itemize}
 \item 1-primitive if $\w$ takes values in the vector space of primitive elements of $H$.
 \item Fully primitive if 1-primitive and both $m_1$ and $m_2$ take values in the space of primitive elements of $I$. 
\end{itemize}
\end{Definition}
\begin{Definition}[2-curvature 3-form]\label{2c3f}
Given a 1-form $\w$ in a manifold $M$ with values in $H$ and a 2-form $m$ with values in $I$, we put $\w \wedge^\t m $, as being the $I$-valued 3-form in $M$, such that, for each vector fields $X,Y,Z$:
$$(\w \wedge^\t m)(X,Y,Z)=\w(X) \t m(Y,Z) + \w(Y) \t m(Z,X) + \w(Z) \t m(X,Y) .$$
 Given a local Hopf 2-connection, described by a triple ${(\w,m_1,m_2)}$, its 2-curvature $\M$ is the graded 3-form:
$$\M=d \bm - \bw \wedge^\t \bm.   $$
Explicitly $\M=h dm_1+h^2dm_2-h^2 \w \wedge^\t m_1-h^3 \w \wedge^\t m_2.$

\end{Definition}

\subsection{The exact {two-dimensional} holonomy}Let $\H=(\d \colon I \to H, \rho)$ be a Hopf algebra crossed module.
Let $M$ be a smooth manifold. 
 Let $(\w,m_1,m_2)$ be a local Hopf 2-connection. We define two (in general) distinct {two-dimensional} holonomies (exact and fuzzy), for a 2-path.
 This subsection is devoted to the first.
 Recall the {definition} of the crossed module of groups  $\H^*_{\rm gl}=(\d\colon I^*_{\rm gl} \to H^*_{\rm gl},\t)$ and of the larger pre-crossed module of groups $\hat{\H}^*_{\rm gl}=(\d\colon \hat{I}^*_{\rm gl} \to H^*_{\rm gl},\t)$; see \ref{hag}. Recall  the constructions of the 2-groupoid \smash{$\Cc^\times \left ( {\H^*_{\rm gl}}\right)$ } and of the larger sesquigroupoid {$\hat{\Cc}^\times({\hat{\H}^*_{\rm gl}})$};  \ref{gcmc}.

\subsubsection{Definition of the exact {two-dimensional} holonomy}\label{exa}
\begin{Definition}[The exact holonomy] Resume the notation of \ref{h1c}  and \ref{ch}. Suppose that {$(\w,m_1,m_2)$} is 1-primitive, thus
 ${P_{\bw}}(\G,[t,0],s)$ is group-like for all $t,s \in [0,1]$. Consider the differential equation in $I\h$:
\begin{equation}\label{eh}
\begin{split}
 \frac{d}{d s} Q_{(\bw,\bm)}(\G,[s,s_0])&=-\left ( \int_0^1 {P_{\bw}}(\G,[0,u],s)\,\,\t \bm \left ( \frac{\d}{\d s} \G(u,s), \frac{\d}{\d u} \G(u,s) \right )du\right)\,  Q_{(\bw,\bm)  }(\G,[s,s_0])\\
Q_{(\bw,\bm)}(\G,[s_0,s_0])&=1_I.
\end{split}
\end{equation}
{By lemma \ref{ivp} its solution  is unique, explicitly being a  {path}-ordered exponential, see below. Put:}
$$Q_{(\bw,\bm)}(\G)\doteq  Q_{(\bw,\bm)}(\G,[1,0]) \in I\h.$$\end{Definition}
\noindent Clearly,  given $s_0,s, s_1 \in [0,1]$, $Q_{(\bw,\bm)}(\G,[s,s_0])$ is invertible, and  we have: 
\begin{equation}\label{comp2}Q_{(\bw,\bm)}(\G,[s_1,s_0])=Q_{(\bw,\bm)}(\G,[s_1,s]) \,\, Q_{(\bw,\bm)}(\G,[s,s_0]).\end{equation}
\begin{Lemma}
 If {$(\w,m_1,m_2)$} is fully primitive then $Q_{(\bw,\bm)}(\G,[s,s_0])\in I^*_{\rm gl}$ is group like for all $s,s_0\in [0,1]$.
\end{Lemma}
\begin{Proof}
This is entirely analogous to the proof of Lemma \ref{glp}. Note also the lemma below.
\end{Proof}
\begin{Lemma}
 Let $\H=(\d \colon I \to H, \rho)$ be a crossed module of Hopf algebras. The given action of $H$ on $I$ is such that if  $x \in H^*_{\rm gl}$ and $a \in \Prim(I)[[h]]$ then 
 $$x \t a \in \Prim(I)[[h]].$$
\end{Lemma}
\begin{Proof}
 We have, since $I$ is an $H$-module coalgebra:
 \begin{align*} \D(x \t a)= \sum_{(x)}\sum_{(a)} x' \t a' \tn x'' \t a'' &=x \t a \tn x \t 1 + x \t 1 \tn x \t a  \\
&=x \t a \tn \e(x) + \e(x) \tn x \t a =x \t a \tn 1 + 1 \tn x \t a.
 \end{align*}
\end{Proof}

\begin{Theorem}[{Exact two-dimensional} holonomy] Suppose that $(\w,m_1,m_2)$ is 1-primitive.
 The assignments:   
 \begin{align}\allowdisplaybreaks\label{F1}
&{\textrm{ for a 1-path } \g \textrm{ in  } M  :} & &&&& &{\Fc^\times_{(\bw,\bm)}}(x \ra{\g} y)=*\ra{P_\bw(\g)^{-1}} *,&&& & \end{align}
 which is a morphism of the sesquigroupod  $\hat{\Cc}^\times({\hat{\H}^*_{\rm gl}})$, and, for a 2-path in $M$,
\begin{equation}\label{F2}
 {\Fc^\times_{(\bw,\bm)}}\left(
{\hskip-1cm
\xymatrix
{ 
& {x} \ar@/^1pc/[rr]^{{\g_2} }\ar@/_1pc/[rr]_{\g_1} & \Uparrow {\G} & {y}
}
} \right)=\hskip-1cm\bbmortimes{\left(P_\bw(\g_1)\right)^{-1}}{ \left( Q_{(\bw,\bm)}(\G)\right)^{-1}}{\left(P_\bw(\g_2)\right)^{-1}}  ,
\end{equation}
{a 2-morphism of the sesquigroupoid  $\hat{\Cc}^\times({\hat{\H}^*_{\rm gl}})$,
preserve composition of {1-path}s,  vertical compositions of 2-paths, whiskerings of 2-paths by 1-paths and vertical reverses or 2-paths. If $(\w,m_1,m_2)$ is fully primitive then the image of ${\Fc^\times_{(\bw,\bm)}}$  lives in $\Cc^\times \left ( {\H^*_{\rm gl}}\right)\subset \hat{\Cc}^\times({\hat{\H}^*_{\rm gl}})$,  ${\Fc^\times_{(\bw,\bm)}}$  preserving horizontal compositions and  reverses.}
\end{Theorem}
\begin{Proof}That the functor preserves the composition of 1-paths and 2-paths follows from \eqref{compr} and \eqref{comp2}. That the right-hand-side of \ref{F2} is a 2-morphism of the sesquigroupoid $\hat{\Cc}^\times({\hat{\H}^*_{\rm gl}})$ follows from Lemma \ref{glp} and equations \eqref{act}, \eqref{2dd} and  \eqref{eh}. The first Peiffer law thus imply that, for each $s_0,s \in  [0,1]$:
\begin{equation}\label{tytyt}
\d\left(Q_{(\bw,\bm)}(\G,[s,s_0])\right)  =P_\bw(\g_{s})^{-1} \, P_\bw(\g_{s_0}),
\end{equation}
 since both of these satisfy differential equation \eqref{2dd}, and have the same initial condition, by {Lemma} \ref{ivp}.

 The most difficult bit is the preservation of horizontal composition and reverses, in the fully primitive case. We only prove the preservation of horizontal compositions, since horizontal reverses are dealt with exactly in the same way. We will   make use of the previous {lemma} and Theorem \ref{rrr}. Let $\G_1$ and $\G_2$ be  2-paths, such that their horizontal composition $\G=\G_1\G_2$ is well defined.  We have:\small
\begin{align*}\allowdisplaybreaks
 \frac{d}{d s} Q_{(\bw,\bm)}(\G,[s,s_0])&=  -\left( \int_0^1 {P_{\bw}}(\G,[0,u],s) \t \bm \left ( \frac{\d}{\d 	s} \G(u,s), \frac{\d}{\d u} \G(u,s) \right )  du\right) \, Q_{(\bw,\bm)}(\G,[s,s_0]) \\
&= - \left (\int_0^1 {P_{\bw}}(\G_1,[0,u],s) \t \bm \left ( \frac{\d}{\d s} \G_1(u,s), \frac{\d}{\d u} \G_1(u,s) \right ) du\right) Q_{(\bw,\bm)}(\G,[s,s_0]) \,\\
& \hskip-1.5cm - \left ( {P_{\bw}}(\G_1,[1,0],s)\t \int_0^1 {P_{\bw}}(\G_2,[0,u]) \t \bm \left ( \frac{\d}{\d s} \G_2(u,s), \frac{\d}{\d u} \G_2(u,s) \right ) du\right)\,\, Q_{(\bw,\bm)}(\G,[s,s_0]).
\end{align*}\normalsize
On the other hand, using the Leibniz law, the second Peiffer relation for group crossed modules, and \eqref{tytyt}:\small
\begin{align*}\allowdisplaybreaks
&\frac{d}{d s} \Big ( Q_{(\bw,\bm)}(\G_1,[s,s_0]) \, {P_{\bw}}(\G_1,[1,0],s_0) \t  Q_{(\bw,\bm)}(\G_2,[s,s_0])\Big)\\
&=-\left ( \int_0^1 {P_{\bw}}(\G_1,[0,u],s)\,\,\t \bm \left ( \frac{\d}{\d s} \G_1(u,s), \frac{\d}{\d u} \G_1(u,s) \right )  du\right)\Big ( Q_{(\bw,\bm)}(\G_1,[s,s_0]) \,\, {P_{\bw}}(\G_1,[0,1],0) \t  Q_{(\bw,\bm)}(\G_2,[s,s_0])\Big)\\
&\quad \quad -  Q_{(\bw,\bm)}(\G_1,[s,s_0]) \,\, {P_{\bw}}(\G_1,[0,1],0) \t   \left ( \int_0^1 {P_{\bw}}(\G_2,[0,u],s)\,\,\t \bm \left ( \frac{\d}{\d s} \G_2(u,s), \frac{\d}{\d u} \G_2(u,[s,s_0]) \right )  du\right)\,\, \\ & \qquad \quad \qquad \quad \qquad \quad \qquad \quad \qquad \quad \qquad \quad\qquad \quad \qquad \quad \qquad \quad \qquad \quad \qquad {P_{\bw}}(\G_1,[0,1],0) \t  Q_{(\bw,\bm)}(\G_2,[s,s_0])
\\
&=-\left ( \int_0^1 {P_{\bw}}(\G_1,[0,u],s)\,\,\t \bm \left ( \frac{\d}{\d s} \G_1(u,s), \frac{\d}{\d u} \G_1(u,s) \right )  du\right)\Big ( Q_{(\bw,\bm)}(\G_1,[s,s_0]) \,\, {P_{\bw}}(\G_1,[0,1],0) \t  Q_{(\bw,\bm)}(\G_2,[s,s_0])\Big)\\
&\quad \quad -  \left (\d \big(Q_{(\bw,\bm)}(\G_1,[s,s_0])\big)  {P_{\bw}}(\G_1,[0,1],0) \right) \t   \left ( \int_0^1 {P_{\bw}}(\G_2,[0,u],s)\,\,\t \bm \left ( \frac{\d}{\d s} \G_2(u,s), \frac{\d}{\d u} \G_2(u,[s,s_0]) \right )\,   du\right)\,\, \\ & \qquad \quad \qquad \quad \qquad \quad \qquad \quad \qquad \quad\qquad \quad \qquad \quad \qquad \quad \qquad   Q_{(\bw,\bm)}(\G_1,[s,s_0]) \left ( {P_{\bw}}(\G_1,[0,1],0) \t  Q_{(\bw,\bm)}(\G_2,[s,s_0])\right)
\\&=-\left ( \int_0^1 {P_{\bw}}(\G_1,[0,u],s)\,\,\t \bm \left ( \frac{\d}{\d s} \G_1(u,s), \frac{\d}{\d u} \G_1(u,s) \right )  du\right)\Big ( Q_{(\bw,\bm)}(\G_1,[s,s_0]) \,\, {P_{\bw}}(\G_1,[0,1],0) \t  Q_{(\bw,\bm)}(\G_2,[s,s_0])\Big)\\
&\quad \quad - \left ( {P_{\bw}}(\G_1,[0,1],s) \right) \t   \left ( \int_0^1 {P_{\bw}}(\G_2,[0,u],s)\,\,\t \bm \left ( \frac{\d}{\d s} \G_2(u,s), \frac{\d}{\d u} \G_2(u,[s,s_0]) \right )\,   du\right)\,\, \\ & \qquad \quad \qquad \quad \qquad \quad \qquad \quad \qquad \quad\qquad \quad \qquad \quad \qquad \quad \qquad   Q_{(\bw,\bm)}(\G_1,[s,s_0]) \left ( {P_{\bw}}(\G_1,[0,1],0) \t  Q_{(\bw,\bm)}(\G_2,[s,s_0])\right).
\end{align*}\normalsize
The two differential equations coincide. So do initial values. Therefore we proved that, for each $s_0,s \in [0,1]$:
\begin{equation}\label{jhg}
 Q_{(\bw,\bm)}(\G,[s,s_0]) =  Q_{(\bw,\bm)}(\G_1,[s,s_0]) \, {P_{\bw}}(\G_1,[0,1],s_0) \t  Q_{(\bw,\bm)}(\G_2,[s,s_0]).
\end{equation}
(Note that we used the second Peiffer relation for group crossed modules, thus  \eqref{jhg} would not be true if non-fully primitive Hopf 2-connections were used.) 
Putting $s=1$ and $s_0=0$ and glancing at equations  \eqref{F2} and \eqref{goodexp}, we see that  $\F$ preserves horizontal compositions of 2-paths.
{Preservation of whiskerings (valid for 1-primitive connections) is dealt with similarly, except that the second Peiffer condition is not needed.}
\end{Proof}
\subsubsection{The exact {two-dimensional} holonomy and 2-curvature}\label{2cccc}
{Consider a 3-path $J$ connecting the 2-paths $\G_0$ and $\G_1$; {Definition} \ref{3path}. Let $\H=(\d \colon I \to H, \rho)$ be a Hopf algebra crossed module. Let $M$ be a smooth manifold, provided with a  primitive local Hopf 2-connection ${(\w,m_1,m_2)}$. The difference between the exact {two-dimensional} holonomies $Q_{(\bw,\bm)}(\G_0)$ and $Q_{(\bw,\bm)}(\G_1)$ is essentially ruled by the 2-curvature 3-form $\M$ of ${(\bw,\bm)}$, {Definition} \ref{2c3f}. In the Lie crossed module case this appeared in \cite{BS,SW2,FMP1}. We  essentially use the argument of \cite{BS}, which was generalized to the pre-crossed module case in \cite{FMP3}.}

Given $t,s,x \in [0,1]$, put ${P_{\bw}}(J,[t,1],s,x)={P_{\bw}}(J_{(s,x)},[t,1])$, where $J_{(s,x)}(t)=J(t,s,x)$, a 1-path. 
Let $[0,1]^2=\big\{(s,x) \in\mathbb{R}^2: s,x \in [0,1]\big \}$. Consider the smooth differential form $a$ in $[0,1]^2$, with values in $I$, with:
\begin{align*}\allowdisplaybreaks a\left (\ds \right)_{(x,s)}& = -\int_0^1 {P_{\bw}}(J,[0,u],s,x)\t \bm \left ( \frac{\d}{\d s} J(u,s,x), \frac{\d}{\d u} J(u,s,x) \right )  du\\
&=-\int_0^1 {P_{\bw}}(J,[0,u],s,x)\t J^*( \bm) \left ( \frac{\d}{\d s} , \frac{\d}{\d u} \right )  du.
\end{align*}
\begin{align*}\allowdisplaybreaks a\left (\dx \right)_{(x,s)}& = -\int_0^1 {P_{\bw}}(J,[0,u],s,x)\t \bm \left ( \frac{\d}{\d x} J(u,s,x), \frac{\d}{\d u} J(u,s,x) \right )  du\\
&=-\int_0^1 {P_{\bw}}(J,[0,u],s,x)\t J^*( \bm) \left ( \frac{\d}{\d x} , \frac{\d}{\d u} \right )  du.
\end{align*}
By using equations \eqref{eh} and \eqref{dsi}, and Leibniz law, its exterior derivative is  seen to be such that:\small
\begin{align*}\allowdisplaybreaks
 & da \left (\dx,\ds \right)=\dx a\left (\ds \right)_{(x,s)} -\ds a\left (\dx \right)_{(x,s)}\\
&=\int_0^1 {P_{\bw}}(J,[0,u],s,x)\t \left (J^*(\bw)\left(\dx\right) \t J^*(\bm) \left ( \frac{\d}{\d s}, \frac{\d}{\d u}  \right )-J^*(\bw)\left(\ds\right) \t J^*\bm \left ( \frac{\d}{\d x}, \frac{\d}{\d u}  \right )\right )du\\
& -\int_0^1 {P_{\bw}}(J,[0,u],s,x)\t \left (\dx  J^*(\bm) \left ( \frac{\d}{\d s}, \frac{\d}{\d u}  \right )-\ds J^*(\bm) \left ( \frac{\d}{\d x}, \frac{\d}{\d u}  \right )\right )  du \\
&+\int_0^1 \int_0^t {P_{\bw}}(J,[0,u],s,x) J^*(\bw) \left ( \frac{\d}{\d x} , \frac{\d}{\d u} \right ) {P_{\bw}}(J,[u,0],s,x)   {P_{\bw}}(J,[0,t],s,x) \t\,J^*\bm \left ( \ds,\dt\right) \, du \, dt\\
&-\int_0^1 \int_0^t {P_{\bw}}(J,[0,u],s,x) J^*(\bw) \left ( \frac{\d}{\d s} , \frac{\d}{\d u} \right ){P_{\bw}}(J,[u,0],s,x)   {P_{\bw}}(J,[0,t],s,x) \t\,J^*(\bm) \left ( \dx,\dt\right)\, du\, dt\\
&=-\int_0^1 {P_{\bw}}(J,[0,u],s,x)  \t J^*(\M)\left (\dx,\ds,\du\right )du\\
& -\int_0^1 {P_{\bw}}(J,[0,u],s,x)\t \left (\du  J^*(\bm) \left ( \frac{\d}{\d x}, \frac{\d}{\d s}  \right )-J^*(\bw)\left(\du \right) \t J^*(\bm) \left ( \frac{\d}{\d x}, \frac{\d}{\d s}  \right )\right )du\\
&+\int_0^1 \int_0^t {P_{\bw}}(J,[0,u],s,x) J^*(\bw) \left ( \frac{\d}{\d x} , \frac{\d}{\d u} \right ) {P_{\bw}}(J,[u,0],s,x)   {P_{\bw}}(J,[0,t],s,x) \t\,J^*(\bm) \left ( \ds,\dt\right) \, du \, dt\\
&-\int_0^1 \int_0^t {P_{\bw}}(J,[0,u],s,x) J^*(\bw) \left ( \frac{\d}{\d s} , \frac{\d}{\d u} \right ){P_{\bw}}(J,[u,0],s,x)   {P_{\bw}}(J,[0,t],s,x) \t\,J^*(\bm) \left ( \dx,\dt\right)\, du\, dt
\end{align*}\normalsize
Note that, integrating by parts, noting  $\frac{\d}{\d s} \G_x(0,s)=0 = \frac{\d}{\d s} \G_x(0,s)=0$, and equation \eqref{h1d}, we have:\small
\begin{align*}\allowdisplaybreaks
 \int_0^1 {P_{\bw}}(J,[0,u],s,x)\t \left (\du  J^*(\bm) \left ( \frac{\d}{\d x}, \frac{\d}{\d s} \right ) \right )du=-\int_0^1 {P_{\bw}}(J,[0,u],s,x)\t \left ( J^*(\bw)\left(\du \right) \t J^*(\bm) \left ( \frac{\d}{\d x}, \frac{\d}{\d s}  \right )\right )du.
\end{align*}\normalsize
This cancels out the second line of the previous formula. Finally:\small
\begin{align*}
 [a,a]\left(\dx,\ds\right)&=\int_0^1\int_0^1 \left [ {P_{\bw}}(J,[0,u],s,x)\t J^*( \bm) \left ( \frac{\d}{\d s} , \frac{\d}{\d u} \right )  du,  {P_{\bw}}(J,[0,t],s,x)\t J^*( \bm) \left ( \frac{\d}{\d x} , \frac{\d}{\d t} \right )\right] du \, dt  \\
 &= \int_0^1\int_0^t \left [ {P_{\bw}}(J,[0,u],s,x)\t J^*( \bm) \left ( \frac{\d}{\d s} , \frac{\d}{\d u} \right )  du,  {P_{\bw}}(J,[0,t],s,x)\t J^*( \bm) \left ( \frac{\d}{\d x} , \frac{\d}{\d t} \right )\right] du \, dt\\
 &+ \int_0^1\int_0^u \left [ {P_{\bw}}(J,[0,u],s,x)\t J^*( \bm) \left ( \frac{\d}{\d s} , \frac{\d}{\d u} \right )  du,  {P_{\bw}}(J,[0,t],s,x)\t J^*( \bm) \left ( \frac{\d}{\d x} , \frac{\d}{\d t} \right )\right] dt \, du.
\end{align*}
\normalsize
Since $\d(\overline{m})=\bw$,  from the above and \eqref{ppps}, we can see that the curvature  $da+[a,a]$ of $a$ satisfies:
\begin{align*}\allowdisplaybreaks
 &\left(da+[a,a]  \right)\left(\dx,\ds \right) =-\int_0^1 {P_{\bw}}(J,[0,u],s,x)\t \left (J^*(\M)\left (\dx,\ds,\du\right )\right )du\\
&-\int_0^1 \int_0^t \left\langle {P_{\bw}}(J,[0,u],s,x)  \t  J^*(\bm) \left ( \frac{\d}{\d x} , \frac{\d}{\d u} \right  )  ,  {P_{\bw}}(J,[0,u],s,x)\t\,J^*(\bm) \left ( \ds,\dt\right)\right\rangle du dt\\
&+\int_0^1 \int_0^t \left\langle {P_{\bw}}(J,[0,u],s,x)   \t J^*(\bm) \left ( \frac{\d}{\d s} , \frac{\d}{\d u} \right  ) ,  {P_{\bw}}(J,[0,u],s,x)\t\,J^*(\bm) \left ( \dx,\dt\right)\right\rangle du dt.
\end{align*}
By using {Lemma} \ref{um}, we can see, given that
{$Q_{(\bw,\bm)}(\G,[s,s_0])
=\sum_{n=0}^{+ \infty} h^n \oint_{s_0}^s \underbrace{a(\g'(v))*\dots*a(\g'(v))}_{n \textrm{ terms}} dv
$}:
\begin{Theorem} We have, using the notation of Definition \ref{3path}, thus $J_x(t,s)=\G(t,s,x)$:   
\begin{align*}\allowdisplaybreaks
 &\frac{d}{dx} Q_{(\bw,\bm)} (J_x)=-\int_0^1 \int_0^1 Q_\w(J,[1,s],x)\M\left (\frac{\d}{\d x} J(t,s,x),\frac{\d}{\d s} J(t,s,x) ,\frac{\d}{\d t} J(t,s,x)\right )Q_\w(J,[s,0],x) \, dt\, ds \\&+
 \int_0^1 \int_0^1 \int_0^t   Q_\w(J,[1,s],x) 
 \Big\langle {P_{\bw}}(J,[0,u],s,x)\t   \bm \left ( \frac{\d}{\d s} J(u,s,x),  \frac{\d}{\d u} J(t,s,x)\right  ) , \\ &\quad \quad \quad \quad \quad \quad \quad \quad \quad \quad  \quad \quad \quad    {P_{\bw}}(J,[0,u],s,x)\t\,\bm \left ( \dx J(t,s,x) ,\dt J(t,s,x)\right)\Big\rangle
 Q_\w(J,[s,0],x) \, du\, dt\, ds
 \\&-
 \int_0^1 \int_0^1 \int_0^t   Q_\w(J,[1,s],x) 
 \Big\langle {P_{\bw}}(J,[0,u],s,x)  \t  \bm \left ( \frac{\d}{\d x} J(u,s,x),  \frac{\d}{\d u} J(t,s,x)\right  )  , \\ &\quad \quad \quad \quad \quad \quad \quad \quad \quad \quad  \quad \quad \quad    {P_{\bw}}(J,[0,u],s,x)\t\,\bm \left ( \ds J(t,s,x) ,\dt J(t,s,x)\right)\Big\rangle
 Q_\w(J,[s,0],x) \, du\, dt\, ds\\ &+\left ( \int_0^1 P(J^t,[0,s],x)\,\t \bm \left ( \frac{\d}{\d x} J^t(s,x), \frac{\d}{\d s}  J^t(s,x) \right ) du\right)  Q_{(\bw,\bm)  }(\G_x)
 \\&
 -Q_{(\bw,\bm)  }(\G_x)\left ( \int_0^1 P(J^b,[0,s],x)\,\t \bm \left ( \frac{\d}{\d x} J^b(s,x), \frac{\d}{\d s}  J^b(s,x) \right ) du\right)\,.  
  \end{align*}
  \end{Theorem}
  {The tensor $\langle, \rangle$ vanishes in the primitive space of $I\h$. Therefore, for fully primitive connections the second and third terms of the previous equation are zero. 
If we, furthermore, suppose that $J$ is a thin homotopy then all of the remaining terms disappear. Thence we can see that the {two-dimensional} holonomy of a fully primitive connection is invariant under thin homotopy; c.f.  the comments at the end of Lemma \ref{um}. Thus:}
\begin{Theorem}
 Suppose that {$(\bw,\bm)$} is fully primitive. The assignments ${\Fc^\times_{(\bw,\bm)}}$, defined in \eqref{F1} and \eqref{F2}, descend to a functor from the fundamental thin 2-groupoid of $M$ (see Remark \ref{tth2g}) into the 2-groupoid ${\Cc}^\times({{\H}^*_{\rm gl}})$.
\end{Theorem}
Also:
\begin{Theorem}
Let {$(\w,m_1,m_2)$} be a fully primitive $\H$-valued Hopf 2-connection in a manifold $M$; where $\H=(\d\colon I \to H, \rho)$. Suppose that the 2-curvature 3-form $\M=d \bm - \bw \wedge^\t \bm $ vanishes (Definition \ref{2c3f}), explicitly $dm_1=0$, $dm_2-\w \wedge^\t m_1=0$ and  $\w \wedge^\t m_2=0$. The {two-dimensional} holonomy $\F^\times_{(\bw,\bm)}(\G) \in I^*_{\rm gl}$ of a 2-path $\G\in \P^2(M)$ depends only on the homotopy class of $\Gamma\colon [0,1]^2 \to M$, relative to the boundary.
\end{Theorem}

\subsection{Bare 2-connections and their holonomy}\label{bare}
 Consider a unital pre-crossed module  $\hat{\A}=(\d\colon A \to B, {\tra,\tla})$, of bare algebras. The following parallels the discussion in subsection  \ref{lhopf}. A local bare 2-connection, on a manifold $M$,  yields a {two-dimensional} holonomy, living in the sesquigroupoid $\hat{\Cc}^+({\hat{\A}})$; see \ref{bacmc}. This sesquigroupoid is a 2-groupoid if $\hat{\A}$ is a crossed module. 
 \begin{Definition} [Bare 2-connection]
 A (local) bare 2-connection is given by a triple $(\w,m_1,m_2)$ of differential forms in $M$, namely
$\w \in \Omega^1(M,B)$ and $ m_1,m_2 \in \Omega^2(M,A),$
such that the vanishing of the fake curvature condition holds:
$\d (m_1)= d \w$ and $\d(m_2)= -\frac{1}{2} \w \wedge^{[,]} \w. $
Given a formal parameter $h$ we can build the pair $(\bw,\bm)$ of graded forms:
$\bw=h\w$ and $\bm=h m_1 +h^2 m_2. $  Therefore we have $\d(\bm)=\F_{\bw} $
where $\F_{\bw} \in \Omega^2(M,B\h)$ is the graded curvature of $\bw$, thus:
$\F_{\bw}=d \bw-\frac{1}{2} \bw \wedge^{[,]} \bw=h d \w-h^2\frac{1}{2} \w \wedge^{[,]} \w $.
 \end{Definition}
 
\begin{Definition}[Curvature of a  bare 2-connection]
 The bare 2-curvature 3-form of the bare 2-connection  {$(\w,m_1,m_2)$} is the graded 3-form $\M \in  \Omega^3(M,A\h)$, such that:
 $$\M=d\overline{m} - \bw \wedge^\t \bm, $$
where $b \t a\doteq b\tra a-a \tla b$, for $a \in A$ and $b \in B$. If $X,Y,Z$ are vector fields in $M$, then for example:
\begin{align*}
(\w \wedge^\t m_1) (X,Y,Z)&= \w(X) \tra m(Y,Z) + \w(Y) \tra m(Z,X) + \w(Z) \tra m(X,Y)\\ & -  m(Y,Z)\tla \w(X)  +  m(Z,X)\tla\w(Y)   +  m(X,Y)\tla \w(Z) .
\end{align*}
 \end{Definition}

\subsubsection{The two dimensional holonomy of a bare 2-connection}\label{2wer}
  We resume the notation of \ref{ch}. Let $M$ be a manifold, $\hat{\A}=(\d\colon A \to B, {\tra,\tla})$ a pre-crossed module of bare algebras, and   $(\w,m_1,m_2)$ a bare 2-connection, seen as a pair $(\bw,\bm)$ of graded forms. 
  Looking at \eqref{2dd}, and considering a 2-path   {$\morts{x}{\g_0}{\G}{\g_1}{y}$} in $M$, we  put: 
 \begin{align}\label{dr}
 R_{(\bw,\bm)}(\G,[s,s_0])=-\int_{s_0}^s \int_0^1 {P_{\bw}}(\G,[0,u],v) \tra \bm \left ( \frac{\d}{\d v} \G(u,v), \frac{\d}{\d u} \G(u,v) \right )\tla {P_{\bw}}(\G,[u,1],v) \,du\,dv \in A.
\end{align}
Also put $R_{(\bw,\bm)}(\G)\doteq R_{(\bw,\bm)}(\G,[1,0])$.
\begin{Theorem}
 The    assignment:
 \begin{align}\allowdisplaybreaks\label{Fp1}
&{\textrm{ for a 1-path } \g \textrm{ in  } M  :} & &&&& &{\Fc^+_{(\bw,\bm)}}(x \ra{\g} y)=*\ra{P_\bw(\g)^{-1}} *,&&& & \end{align}
 which is a morphism of the sesquigroupoid    $\hat{\Cc}^+({\hat{\A}})$, and, for a 2-path in $M$,
\begin{equation}\label{Fp2}
 {\Fc^+_{(\bw,\bm)}}\left(
{\hskip-1cm
\xymatrix
{ 
& {x} \ar@/^1pc/[rr]^{{\g_2} }\ar@/_1pc/[rr]_{\g_1} & \Uparrow {\G} & {y}
}
} \right)=\hskip-1cm\morplus{\left(P_\bw(\g_1)\right)^{-1}}{ R_{(\bw,\bm)}(\G)}{\left(P_\bw(\g_2)\right)^{-1}}  ,
\end{equation}
 which is a 2-morphism of the sesquigroupoid   $\hat{\Cc}^+({\hat{\A}})$,
preserves composition of {1-path}s,  vertical compositions of 2-paths, whiskerings of 2-paths by 1-paths, and vertical reverses or 2-paths. If furthermore $\hat{\A}$ is a crossed module of bare algebras, then  ${\Fc^+_{(\bw,\bm)}}$ also preserves horizontal compositions and horizontal reverses of 2-paths.
\end{Theorem}
\begin{Proof}
 That the image of ${\Fc^+_{(\bw,\bm)}}$ is in the sesquigroupoid  $\hat{\Cc}_{\hat{\A}}^+$ follows from \eqref{dr}, \eqref{2dd}, \eqref{compr} and the {definition} of crossed modules of bare algebras; see \ref{dcmba}. These imply that, for each $s,s_0 \in [0,1]$, and each 2-path $\G$:
 \begin{align}\label{relp}
  P_\bw(\G,[0,u],s) = P_\bw(\G,[0,u],s_0)+\d\left(R_{(\bw,\bm)}(\G,[s,s_0])\right).
 \end{align}

 The most difficult bit is the preservation of horizontal compositions and reverses in the crossed module case. We deal only with horizontal compositions, since horizontal reverses are dealt with similarly. Consider two 2-paths $\G_1$ and $\G_2$ such that $\G=\G_1\G_2$, their horizontal composition, is well defined. We have:
  \begin{align*}
 R_{(\bw,\bm)}&(\G,[s,s_0])=-\int_{s_0}^s \int_0^1 {P_{\bw}}(\G,[0,u],v) \tra \bm \left ( \frac{\d}{\d v} \G(u,v), \frac{\d}{\d u} \G(u,v) \right )\tla {P_{\bw}}(\G,[u,1],v) \,du\,dv \\
 &= \quad-\int_{s_0}^s \int_0^1 {P_{\bw}}(\G_1,[0,1],v) {P_{\bw}}(\G_2,[0,u],v) \tra \bm \left ( \frac{\d}{\d v} \G_2(u,v), \frac{\d}{\d u} \G_2(u,v) \right )\tla {P_{\bw}}(\G_2,[u,1],v) \,du\,dv \\
 &-\int_{s_0}^s \int_0^1 {P_{\bw}}(\G_1,[0,u],v) \tra \bm \left ( \frac{\d}{\d v} \G_2(u,v), \frac{\d}{\d u} \G_2(u,v) \right )\tla {P_{\bw}}(\G_1,[u,1],v){P_{\bw}}(\G_2,[0,1],v) \,du\,dv\, .
\end{align*}
Now use \eqref{relp} in 
${P_{\bw}}(\G_2,[0,1],v)$ (end of the bottom line) and in ${P_{\bw}}(\G_1,[0,1],v)$ (beginning of the top line):
 \begin{align*}
 R_{(\bw,\bm)}&(\G,[s,s_0])=R_{(\bw,\bm)}(\G_1,[s,s_0]) \tla {P_{\bw}}(\G_2,[0,1],s_0)+ {P_{\bw}}(\G_1,[0,u],s_0) \tra R_{(\bw,\bm)}(\G_2,[s,s_0])+A,
 \end{align*}
where, using \eqref{dr}, and ulteriorly the second Peiffer relation for crossed modules of bare algebras:
\begin{align*}\allowdisplaybreaks
 A&= \int_{s_0}^s \int_0^1\int_{s_0}^v \int_0^1  {P_{\bw}}(\G_1,[0,u],v) \tra \bm \left ( \frac{\d}{\d v} \G_1(u,v), \frac{\d}{\d u} \G_1(u,v) \right )\tla {P_{\bw}}(\G_1,[u,1],v)\\ & \quad \quad \quad \d\left( {P_{\bw}}(\G_2,[0,u'],v') \tra \bm \left ( \frac{\d}{\d v'} \G_2(u',v'), \frac{\d}{\d u'} \G_2(u',v') \right )\tla {P_{\bw}}(\G_2,[u',1],v')\right) \,du'\,dv'  \,du\,dv \\
 &+\int_{s_0}^s \int_0^1 \int_{s_0}^v \int_0^1  \d \left ({P_{\bw}}(\G_1,[0,u'],v') \tra \bm \left ( \frac{\d}{\d v'} \G_1(u',v'), \frac{\d}{\d u'} \G_1(u',v') \right )\tla {P_{\bw}}(\G_1,[u',1],v')\right)\\ & \quad \quad\quad \quad {P_{\bw}}(\G_2,[0,u],v) \tra \bm \left ( \frac{\d}{\d v} \G_2(u,v), \frac{\d}{\d u} \G_2(u,v) \right )\tla {P_{\bw}}(\G_2,[u,1],v) \,du'\,dv' \, du\, dv\\
 &= \int_{s_0}^s \int_0^1\int_{s_0}^v \int_0^1  {P_{\bw}}(\G_1,[0,u],v) \tra \bm \left ( \frac{\d}{\d v} \G_1(u,v), \frac{\d}{\d u} \G_1(u,v) \right )\tla {P_{\bw}}(\G_1,[u,1],v)\\ & \quad \quad \quad  {P_{\bw}}(\G_2,[0,u'],v') \tra \bm \left ( \frac{\d}{\d v'} \G_2(u',v'), \frac{\d}{\d u'} \G_2(u',v') \right )\tla {P_{\bw}}(\G_2,[u',1],v') \,du'\,dv'  \,du\,dv \\
 &+\int_{s_0}^s \int_0^1 \int_{s_0}^v \int_0^1  {P_{\bw}}(\G_1,[0,u'],v') \tra \bm \left ( \frac{\d}{\d v'} \G_1(u',v'), \frac{\d}{\d u'} \G_1(u',v') \right )\tla {P_{\bw}}(\G_1,[u',1],v')\\ & \quad \quad\quad \quad {P_{\bw}}(\G_2,[0,u],v) \tra \bm \left ( \frac{\d}{\d v} \G_2(u,v), \frac{\d}{\d u} \G_2(u,v) \right )\tla {P_{\bw}}(\G_2,[u,1],v) \,du'\,dv' \, du\, dv\\
 &=\int_{s_0}^s \int_0^1 \int_{s_0}^s \int_0^1  {P_{\bw}}(\G_1,[0,u'],v') \tra \bm \left ( \frac{\d}{\d v'} \G_1(u',v'), \frac{\d}{\d u'} \G_1(u',v') \right )\tla {P_{\bw}}(\G_1,[u',1],v')\\ & \quad \quad\quad \quad {P_{\bw}}(\G_2,[0,u],v) \tra \bm \left ( \frac{\d}{\d v} \G_2(u,v), \frac{\d}{\d u} \G_2(u,v) \right )\tla {P_{\bw}}(\G_2,[u,1],v) \,du\,dv \, du'\, dv'\\
 &=R_{(\bw,\bm)}(\G_1,[s,s_0])\,\, R_{(\bw,\bm)}(\G_2,[s,s_0]).
\end{align*}
Therefore the following holds in the crossed module case (but not-necessarily in the pre-crossed module case): 
 \begin{multline}\label{gcomp}
 R_{(\bw,\bm)}\big((\G_1\G_2),[s,s_0]\big)=R_{(\bw,\bm)}(\G_1,[s,s_0]) \tla {P_{\bw}}(\G_2,[0,1],s_0)\\+ {P_{\bw}}(\G_1,[0,u],s_0) \tra R_{(\bw,\bm)}(\G_2,[s,s_0])+R_{(\bw,\bm)}(\G_1,[s,s_0]) R_{(\bw,\bm)}(\G_2,[s,s_0]).
\end{multline}
Putting $s_0=0, s_1=1$ and looking at \eqref{phc},  we see that ${\Fc^+_{(\bw,\bm)}}$ preserves  horizontal composition of 2-paths. Horizontal reverses are dealt with similarly.  \end{Proof}

\subsubsection{Bare 2-curvature and two dimensional holonomy}\label{bcbh}
Recall {Definition} \eqref{3path}. Consider a 3-path $J$ connecting the 2-paths $\G_0$ and $\G_1$. Consider a unital pre-crossed module  $\hat{\A}=(\d\colon A \to B, {\tra,\tla})$, of bare algebras. Consider a bare 2-connection  $(\w,m_1,m_2)$ in a manifold $M$. 
The difference between the bare {two-dimensional} holonomies $R_{(\bw,\bm)}(\G_0)$ and $R_{(\bw,\bm)}(\G_1)$ is determined by the 2-curvature 3-form $\M$ of ${(\bw,\bm)}$, {Definition} \ref{2c3f}. In the chain-complex case this appeared in \cite{CFM2}, and we use essentially the same argument. 

 Given $x \in [0,1]$, we have a 2-path $\G_x$, where $\G_x(t,s)=J(t,s,x)$. For $t,s,x \in [0,1]$, put ${P_{\bw}}(J,[t,1],s,x)={P_{\bw}}(J_{(s,x)},[t,1])$, where $J_{(s,x)}(t)=J(t,s,x)$. We have (where we use Leibniz law and equations \eqref{dsi2}, \eqref{dsi3}):
\begin{align*}\allowdisplaybreaks
 \frac{d }{ dx} & R_{(\bw,\bm)}(\G_x) =- \frac{d }{ dx}\int_{0}^1 \int_0^1 {P_{\bw}}(\G_x,[0,u],v) \tra \bm \left ( \frac{\d}{\d v} \G_x(u,v), \frac{\d}{\d u} \G_x(u,v) \right )\tla {P_{\bw}}(\G_x,[u,1],v) \,du\,dv \\
  =&- \int_{0}^1 \left( \frac{d }{ dx} \int_0^1 {P_{\bw}}(\G_x,[0,u],v)\right) \tra \bm \left ( \frac{\d}{\d v} \G_x(u,v), \frac{\d}{\d u} \G_x(u,v) \right )\tla {P_{\bw}}(\G_x,[u,1],v) \,du\,dv 
  \\
 &- \int_{0}^1 \int_0^1 {P_{\bw}}(\G_x,[0,u],v) \tra J^*(\bm) \left ( \frac{\d}{\d v} , \frac{\d}{\d u}  \right )\tla \left (\frac{d }{ dx}  {P_{\bw}}(\G_x,[u,1],v)\right) \,du\,dv \\
  &- \int_{0}^1 \int_0^1 {P_{\bw}}(\G_x,[0,u],v) \tra\left ( \frac{d }{ dx}\bm \left ( \frac{\d}{\d v} \G_x(u,v), \frac{\d}{\d u} \G_x(u,v) \right )\right) \tla {P_{\bw}}(\G_x,[u,1],v) \,du\,dv \\
 &= \int_{0}^1 \int_0^1   \int_0^u {P_{\bw}}(\G_x,[0,u'],v) \, \F_{\bw} \left ( \frac{\d}{\d x} \G_x(u',v), \frac{\d}{\d u'} \G_x(u',v) \right )\,\\ & \quad \quad \quad\quad \quad \quad \quad  {P_{\bw}}(\G_x,[u',u],v)\tra J^*(\bm) \left ( \frac{\d}{\d v} , \frac{\d}{\d u}  \right )\tla {P_{\bw}}(\G_x,[u,1],v)\, du' \,du\,dv
 \\& +  \int_{0}^1 \int_0^1   {P_{\bw}}(\G,[0,u],v) \bw \left(\frac{\d}{\d x} \G_x(u,v)\right)  \tra J^*(\bm) \left ( \frac{\d}{\d v} , \frac{\d}{\d u}  \right )\tla {P_{\bw}}(\G_x,[u,1],v) \,du\,dv\\ &
 + \int_{0}^1 \int_0^1 \int_u^1 {P_{\bw}}(\G_x,[0,u],v) \tra J^*(\bm) \left ( \frac{\d}{\d v} , \frac{\d}{\d u}  \right )\tla  {P_{\bw}}(\G,[u,u'],v) \\ & \quad \quad \quad \quad \quad \quad \, \F_{\bw} \left ( \frac{\d}{\d x} \G_x(u,v), \frac{\d}{\d u'} \G_x(u,v) \right )\, {P_{\bw}}(\G_x,[u',1],v)\, du' \,du\,dv\\
 &- \int_{0}^1 \int_0^1 {P_{\bw}}(\G_x,[0,u],v) \tra J^*(\bm) \left ( \frac{\d}{\d v} , \frac{\d}{\d u}  \right )\tla  \bw \left(\frac{\d}{\d x} \G_x(u,v)\right)  {P_{\bw}}(\G,[u,1],v). \,du\,dv \\
 & - \int_{0}^1 \int_0^1 {P_{\bw}}(\G_x,[0,u],v) \tra  J^*( d \bm) \left (\dx, \frac{\d}{\d v}, \frac{\d}{\d u}  \right )\tla {P_{\bw}}(\G_x,[u,1],v) \,du\,dv \\
 & + \int_{0}^1 \int_0^1 {P_{\bw}}(\G_x,[0,u],v) \tra  \left (\frac{\d}{\d v}  J^*(\bm)\left ( \frac{\d}{\d u},\dx\right )+ \frac{\d}{\d u}  J^*(\bm)\left (\dx,  \frac{\d}{\d v}\right )  \right)\tla {P_{\bw}}(\G_x,[u,1],v) \,du\,dv .
\end{align*} 
Now look at the last term. Integrating by parts, and using   \eqref{dsi2} and  \eqref{dsi3},  we have:
\begin{align*}\allowdisplaybreaks
 & \int_{0}^1 \int_0^1 {P_{\bw}}(\G_x,[0,u],v) \tra  \frac{\d}{\d v}  J^*(\bm)\left ( \frac{\d}{\d u},\dx\right )  \tla {P_{\bw}}(\G_x,[u,1],v) \,du\,dv &\\&=
 - \int_{0}^1 \int_0^1 \left( \frac{\d}{\d v}  {P_{\bw}}(\G_x,[0,u],v)\right) \tra   J^*(\bm)\left ( \frac{\d}{\d u},\dx\right )  \tla {P_{\bw}}(\G_x,[u,1],v) \,du\,dv 
 \\ &\quad -\int_{0}^1 \int_0^1  {P_{\bw}}(\G_x,[0,u],v) \tra   J^*(\bm)\left ( \frac{\d}{\d u},\dx\right )  \tla  \left( \frac{\d}{\d v}  {P_{\bw}}(\G_x,[u,1],v)\right) \,du\,dv \\&\quad +
  \int_{0}^1  {P_{\bw}}(\G_x,[0,u],1) \tra  J^*(\bm)\left ( \frac{\d}{\d u},\dx\right )  \tla {P_{\bw}}(\G_x,[u,1],1)\,du\\&\quad -
  \int_{0}^1  {P_{\bw}}(\G_x,[0,u],0) \tra  J^*(\bm)\left ( \frac{\d}{\d u},\dx\right )  \tla {P_{\bw}}(\G_x,[u,1],0) \,du\\
 &=\int_{0}^1 \int_0^1   \int_0^u {P_{\bw}}(\G_x,[0,u'],v) \, \F_{\bw} \left ( \frac{\d}{\d v} \G_x(u',v), \frac{\d}{\d u'} \G_x(u',v) \right )\,\\ & \quad \quad \quad\quad \quad \quad \quad  {P_{\bw}}(\G,[u',u],v)\tra J^*(\bm) \left ( \frac{\d}{\d u} , \frac{\d}{\d x}  \right )\tla {P_{\bw}}(\G_x,[u,1],v)\, du' \,du\,dv
 \\&\quad +  \int_{0}^1 \int_0^1   {P_{\bw}}(\G,[0,u],v) \bw \left(\frac{\d}{\d v} \G_x(u,v)\right)  \tra J^*(\bm) \left ( \frac{\d}{\d u} , \frac{\d}{\d x}  \right )\tla {P_{\bw}}(\G_x,[u,1],v) \,du\,dv\\ &\quad
 + \int_{0}^1 \int_0^1 \int_u^1 {P_{\bw}}(\G_x,[0,u],v) \tra J^*(\bm) \left ( \frac{\d}{\d u} , \frac{\d}{\d x}  \right )\tla  {P_{\bw}}(\G,[u,u'],v) \\ & \quad \quad \quad \quad \quad \quad \, \F_{\bw} \left ( \frac{\d}{\d v} \G_x(u,v), \frac{\d}{\d u'} \G_x(u,v) \right )\, {P_{\bw}}(\G_x,[u',1],v)\, du' \,du\,dv\\
 &\quad- \int_{0}^1 \int_0^1 {P_{\bw}}(\G_x,[0,u],v) \tra J^*(\bm) \left ( \frac{\d}{\d u} , \frac{\d}{\d x}  \right )\tla  \bw \left(\frac{\d}{\d v} \G_x(u,v)\right)  {P_{\bw}}(\G,[u,1],v). \,du\,dv \\&\quad +
  \int_{0}^1  {P_{\bw}}(\G_x,[0,u],1) \tra  J^*(\bm)\left ( \frac{\d}{\d u},\dx\right )  \tla {P_{\bw}}(\G_x,[u,1],1)\,du\\&\quad -
  \int_{0}^1  {P_{\bw}}(\G_x,[0,u],0) \tra  J^*(\bm)\left ( \frac{\d}{\d u},\dx\right )  \tla {P_{\bw}}(\G_x,[u,1],0) \,du.
\end{align*}
Analogously, integrating by parts, using \eqref{h1d}, and noting that $J(0,s,x)$ and $J(1,s,x)$ each are constant:
\begin{align*}\allowdisplaybreaks
 \int_{0}^1 \int_0^1 & {P_{\bw}}(\G_x,[0,u],v) \tra \frac{\d}{\d u}  J^*(\bm)\left (\dx,  \frac{\d}{\d v}\right )  \tla {P_{\bw}}(\G_x,[u,1],v) \,du\,dv\\ &=  \int_{0}^1 \int_0^1 {P_{\bw}}(\G_x,[0,u],v) \bw \left(\frac{\d}{\d u} \G_x(u,v)\right)\tra   J^*(\bm)\left (\dx,  \frac{\d}{\d v}\right )  \tla {P_{\bw}}(\G_x,[u,1],v) \,du\,dv
 \\ &\quad - \int_{0}^1 \int_0^1 {P_{\bw}}(\G_x,[0,u],v) \tra  J^*(\bm)\left (\dx,  \frac{\d}{\d v}\right )  \tla \bw \left(\frac{\d}{\d u} \G_x(u,v)\right) {P_{\bw}}(\G_x,[u,1],v) \,du\,dv.
\end{align*}
Putting everything together we see:
\begin{align*}\allowdisplaybreaks
\frac{d }{ dx}  R_{(\bw,\bm)}(\G_x) &=A- \int_{0}^1 \int_0^1  {P_{\bw}}(\G_x,[0,u],v) \tra   J^*( \M) \left (\dx, \frac{\d}{\d v}, \frac{\d}{\d u}  \right )  \tla {P_{\bw}}(\G_x,[u,1],v) \,du\,dv  \\&\quad\quad +
  \int_{0}^1  {P_{\bw}}(\G_x,[0,u],1) \tra  J^*(\bm)\left ( \frac{\d}{\d u},\dx\right )  \tla {P_{\bw}}(\G_x,[u,1],1)\,du\\&\quad\quad -
  \int_{0}^1  {P_{\bw}}(\G_x,[0,u],0) \tra  J^*(\bm)\left ( \frac{\d}{\d u},\dx\right )  \tla {P_{\bw}}(\G_x,[u,1],0) \,du,
\end{align*}
where: 
\begin{align*}\allowdisplaybreaks
 A&= \int_{0}^1 \int_0^1   \int_0^u {P_{\bw}}(\G_x,[0,u'],v) \, \F_{\bw} \left ( \frac{\d}{\d x} \G_x(u',v), \frac{\d}{\d u'} \G_x(u',v) \right )\,\\ & \quad \quad \quad\quad \quad \quad \quad  {P_{\bw}}(\G_x,[u',u],v)\tra J^*(\bm) \left ( \frac{\d}{\d v} , \frac{\d}{\d u}  \right )\tla {P_{\bw}}(\G_x,[u,1],v)\, du' \,du\,dv\\
&+ \int_{0}^1 \int_0^1 \int_u^1 {P_{\bw}}(\G_x,[0,u],v) \tra J^*(\bm) \left ( \frac{\d}{\d v} , \frac{\d}{\d u}  \right )\tla  {P_{\bw}}(\G,[u,u'],v) \\ & \quad \quad \quad \quad \quad \quad \, \F_{\bw} \left ( \frac{\d}{\d x} \G_x(u,v), \frac{\d}{\d u'} \G_x(u,v) \right )\, {P_{\bw}}(\G_x,[u',1],v)\, du' \,du\,dv\\
 &+\int_{0}^1 \int_0^1   \int_0^u {P_{\bw}}(\G_x,[0,u'],v) \, \F_{\bw} \left ( \frac{\d}{\d v} \G_x(u',v), \frac{\d}{\d u'} \G_x(u',v) \right )\,\\ & \quad \quad \quad\quad \quad \quad \quad  {P_{\bw}}(\G,[u',u],v)\tra J^*(\bm) \left ( \frac{\d}{\d u} , \frac{\d}{\d x}  \right )\tla {P_{\bw}}(\G_x,[u,1],v)\, du' \,du\,dv\\ &
 + \int_{0}^1 \int_0^1 \int_u^1 {P_{\bw}}(\G_x,[0,u],v) \tra J^*(\bm) \left ( \frac{\d}{\d u} , \frac{\d}{\d x}  \right )\tla  {P_{\bw}}(\G,[u,u'],v) \\ & \quad \quad \quad \quad \quad \quad \, \F_{\bw} \left ( \frac{\d}{\d v} \G_x(u,v), \frac{\d}{\d u'} \G_x(u,v) \right )\, {P_{\bw}}(\G_x,[u',1],v)\, du' \,du\,dv=B+C+D+E.
\end{align*}
Swapping orders of integration we have:
\begin{align*}\allowdisplaybreaks
 B&=\int_{0}^1 \int_0^1   \int_{u'}^1 {P_{\bw}}(\G_x,[0,u'],v) \, \F_{\bw} \left ( \frac{\d}{\d x} \G_x(u',v), \frac{\d}{\d u'} \G_x(u',v) \right )\,\\ & \quad \quad \quad\quad \quad \quad \quad  {P_{\bw}}(\G_x,[u',u],v)\tra J^*(\bm) \left ( \frac{\d}{\d v} , \frac{\d}{\d u}  \right )\tla {P_{\bw}}(\G_x,[u,1],v)\, du \,du'\,dv.
\end{align*}
Let us now use the fact that $\d(\bm)=\bw$, together with equation \eqref{ppba}, to see that:
\begin{align*}
 B+E&=\int_{0}^1 \int_0^1   \int_{u'}^1 \Big \{ {P_{\bw}}(\G_x,[0,u'],v) \, \bm \left ( \frac{\d}{\d x} \G_x(u',v), \frac{\d}{\d u'} \G_x(u',v) \right ),\,\\ & \quad \quad \quad\quad \quad \quad \quad  {P_{\bw}}(\G_x,[u',u],v)\tra J^*(\bm) \left ( \frac{\d}{\d v} , \frac{\d}{\d u}  \right )\tla {P_{\bw}}(\G_x,[u,1],v)\Big\}\, du \,du'\,dv.
\end{align*}
We analogously reach:
\begin{align*}
 C+D&=\int_{0}^1 \int_0^1   \int_0^u \Big \{ {P_{\bw}}(\G_x,[0,u'],v) \, \bm \left ( \frac{\d}{\d v} \G_x(u',v), \frac{\d}{\d u'} \G_x(u',v) \right ), \,\\ & \quad \quad \quad\quad \quad \quad \quad  {P_{\bw}}(\G,[u',u],v)\tra J^*(\bm) \left ( \frac{\d}{\d u} , \frac{\d}{\d x}  \right )\tla {P_{\bw}}(\G_x,[u,1],v)\Big\}\, du' \,du\,dv.
\end{align*}

\begin{Theorem}
 In the conditions stated in the beginning of \ref{bcbh} we have:
 \begin{align*}\allowdisplaybreaks
 \frac{d }{ dx}& R_{(\bw,\bm)}(\G_x)\\ &=- \int_{0}^1 \int_0^1  {P_{\bw}}(\G_x,[0,u],v) \tra    \M \left (\dx J(u,v,x), \frac{\d}{\d v} J(u,v,x), \frac{\d}{\d u} J(u,v,x) \right )  \tla {P_{\bw}}(\G_x,[u,1],v) \,du\,dv  \\&\quad +
  \int_{0}^1  {P_{\bw}}(\G_x,[0,u],1) \tra  \bm\left ( \frac{\d}{\d u} J(u,1,x),\dx J(u,1,x)\right )  \tla {P_{\bw}}(\G_x,[u,1],1)\,du\\&\quad -
  \int_{0}^1  {P_{\bw}}(\G_x,[0,u],0) \tra  \bm\left ( \frac{\d}{\d u} J(u,0,x) ,\dx   J(u,0,x) \right )  \tla {P_{\bw}}(\G_x,[u,1],0) \,du\\
  &\quad +\int_{0}^1 \int_0^1   \int_{u'}^1 \Big \{ {P_{\bw}}(\G_x,[0,u'],v) \, \bm \left ( \frac{\d}{\d x} \G_x(u',v), \frac{\d}{\d u'} \G_x(u',v) \right ),\,\\ & \quad \quad \quad\quad \quad \quad \quad  {P_{\bw}}(\G_x,[u',u],v)\tra J^*(\bm) \left ( \frac{\d}{\d v} , \frac{\d}{\d u}  \right )\tla {P_{\bw}}(\G_x,[u,1],v)\Big\}\, du \,du'\,dv\\
  &\quad +\int_{0}^1 \int_0^1   \int_0^u \Big \{{P_{\bw}}(\G_x,[0,u'],v) \, \bm \left ( \frac{\d}{\d v} \G_x(u',v), \frac{\d}{\d u'} \G_x(u',v) \right ), \,\\ & \quad \quad \quad\quad \quad \quad \quad  {P_{\bw}}(\G,[u',u],v)\tra J^*(\bm) \left ( \frac{\d}{\d u} , \frac{\d}{\d x}  \right )\tla {P_{\bw}}(\G_x,[u,1],v)\Big\}\, du' \,du\,dv .
\end{align*}
\end{Theorem}
 \begin{Remark}
{ In a crossed module of bare algebras the tensor $\{,\}$ of \eqref{ppba} vanishes. If we furthermore suppose that $J$ is a thin homotopy then all of the remaining terms of the previous equation vanish, thence we can see that the {two-dimensional} holonomy of a bare 2-connection with values in a crossed module of bare algebras is  invariant under thin homotopy. Also in the crossed module case, the {two-dimensional} holonomy $\F^+_{(\bw,\bm)}(\G)$ of a bare 2-connection with zero 2-curvature 3-tensor depends only on the homotopy class of the 2-path $\G\colon [0,1]^2 \to M$, relative to the boundary. }
\end{Remark}
From this remark it follows that:
\begin{Corollary} Suppose that $\hat{\A}=(\d\colon A \to B, {\tra,\tla})$ is a crossed module of bare algebras. 
 The assignment ${\Fc^+_{(\bw,\bm)}}$ defined in \eqref{Fp1} and \eqref{Fp2} descends to a 2-functor from the fundamental thin 2-groupoid of $M$ (see remark \ref{tth2g}) into the sesquigroupoid $\hat{\Cc}^+({\hat{\A}})$, which in this case a 2-groupoid, since $\hat{\A}$ is a crossed module. 
\end{Corollary}

\subsection{The blur and fuzzy holonomies of a Hopf 2-connection}\label{fuzzy}

 Let $\H=(\d \colon I \to H, \rho)$ be a crossed module of Hopf algebras. As explained in \ref{qwer}, the notation and nomenclature of which we freely use,  consider the bare algebra $I^0\subset I$, where $I^0={\rm ker}(\e\colon I \to \C)$, which is {closed} under the action $\rho$ of $H$, since $I$ is a coalgebra module. Also recall the  pre-crossed module of bare algebras
$\widehat{\BA}(\H)=\big(t \colon I^0 \otimes_\rho H \to H, \tra, \tla\big).$
We will also consider its associated bare algebra crossed module ${\BA}(\H)={\cal R}(\widehat{\BA}(\H))=\big(t \colon 
\underline{I^0\otimes_\rho H} \to H, \tra, \tla\big);$ see \ref{refl}.
 Recall that we have  algebra maps:
 \begin{align*}
i_I\colon u \in  I & \mapsto u \tn 1_H \in I \tn_\rho H\\ i_H\colon x \in  H & \mapsto 1_I \tn x \in  I \tn_\rho H.
\end{align*}
Clearly, $i_I$ restricts to a map $i_H \colon \Prim(I) \to  I^0\otimes_\rho H $.

One has a sesquigroupoid $\hat{\Cc}^+(\widehat{\BA}(\H))$ and a 2-groupoid   ${\Cc}^+({\BA(\H)})$; \ref{bacmc}. Recall Definition \ref{fully}. The blur holonomy of an $\H$-valued fully-primitive Hopf 2-connection $(\w,m_1,m_2)$  takes values in $\hat{\Cc}^+({\widehat{\BA}(\H)})$, whereas its fuzzy holonomy takes values in ${\Cc}^+({\BA(\H)})$.
Our starting point is equations  \eqref{2dd} and \eqref{2mb}.
\begin{Definition}[Blur Holonomy and Fuzzy Holonomy]
{Let $M$ be a manifold, provided with a Hopf 2-connection ${(\w,m_1,m_2)}$, whose graded form is $(\bw,\bm)$, taken to be fully primitive (it  suffices that both $m_1$ and $m_2$ take values in $I^0 \subset I$, however being essential that $\w$ be primitive.)  We can see $\big(i_H(\w),i_I(m_1),i_H(m_2)\big)$, whose graded form we denote by  $i(\bw,\bm))$, as being a bare 2-connection with values in the pre-crossed module $\widehat{\BA}(\H)$ or, projecting onto $\underline{I^0\otimes_\rho H}$, as a bare 2-connection with values in the quotient crossed module $\BA(\H)$. The {two-dimensional} holonomy (as a bare 2-connection) of the former will be called the blur holonomy ${\hat{\Fz}^+_{(\bw,\bm)}}\colon \P(M) \to {\hat\Cc}^+(\widehat{{\cal BA}}(\H))$, of ${(\bw,\bm)}$. We call the {two-dimensional} holonomy of the latter the fuzzy holonomy ${{\Fz}^+_{(\bw,\bm)}}\colon \P(M) \to \Cc^+({\cal BA}(\H))$, of $(\w,m_1,m_2)$. }
\end{Definition}
Therefore, looking at \ref{2wer}, the blur holonomy of a 2-path  \smash{$\morts{x}{\g_0}{\G}{\g_1}{y}$} has the form below:
\begin{equation}
 {\hat{\Fz}^+_{(\bw,\bm)}}\left (  \morts{x}{\g_0}{\G}{\g_1}{y}\right ) =\hskip-1cm\morplus{\left(P_\bw(\g_1)\right)^{-1}}{ R_{(\bw,\bm)}(\G)}{\left(P_\bw(\g_2)\right)^{-1}}
\end{equation}
where, using the notation of Theorem \ref{majmain0}, and given that one-dimensional holonomies are group-like:
\begin{equation}\label{lkjh}
\begin{split}
&R_{(\bw,\bm)}(\G)= -\int_0^1\int_0^1 i_H\left ({P_{\bw}}(\G,[0,u],s)\right)  \,i_I\left ( \bm \left ( \frac{\d}{\d s} \G(u,s), \frac{\d}{\d u} \G(u,s) \right )\right ) i\,_H\left ( {P_{\bw}}(\G,[u,1],s)\right ) \, du \,ds\\
                &=-\int_0^1\int_0^1 \left ({P_{\bw}}(\G,[0,u],s) \t   \bm \left ( \frac{\d}{\d s} \G(u,s), \frac{\d}{\d u} \G(u,s) \right ) \right )\tn   {P_{\bw}}(\G,[0,1],s) \, du \,ds \in (I^0 \otimes_\rho H)[[h]]. 
\end{split}
\end{equation}

Recall that if $\hat{\A}=(\d\colon A \to B,\tra,\tla)$ is a pre-crossed module of bare algebras, then the associated crossed module of Hopf algebras is denoted by ${\cal R}(A)= (\d\colon \underline{A} \to B,\tra,\tla)$; see \ref{refl}. The algebra projection $A \to \underline{A}$ is denoted by $a \mapsto \underline{a}$. 
The fuzzy holonomy $\Fz^+_{(\bw,\bm)}$, of the Hopf 2-connection ${(\w,m_1,m_2)}$,  which is a projection of the seemingly much finer, blur holonomy $\widehat{\Fz^+}_{(\bw,\bm)}$, has the form:
\begin{equation}
 {{\Fz}^+_{(\bw,\bm)}}\left (  \morts{x}{\g_0}{\G}{\g_1}{y}\right ) =\hskip-1cm\morplus{\left(P_\bw(\g_1)\right)^{-1}}{ \underline{R_{(\bw,\bm)}}(\G)}{\left(P_\bw(\g_2)\right)^{-1}}
\end{equation}
We have put:
\begin{equation}\label{lkjhl}
\underline{R_{(\bw,\bm)}}(\G)=-\int_0^1\int_0^1 \underline{\left ({P_{\bw}}(\G,[0,u],s) \t   \bm \left ( \frac{\d}{\d s} \G(u,s), \frac{\d}{\d u} \G(u,s) \right ) \right )\tn   {P_{\bw}}(\G,[0,1],s) }\, du \,ds .
\end{equation}
By construction, the fuzzy holonomy ${{\Fz}^+_{(\bw,\bm)}}\colon \P(M) \to \Cc^+({\cal BA}(\H))$ of $(\bw,\bm)$ preserves all compositions in $\P(M)$ and in   $\Cc^+({\cal BA}(\H))$. The blur holonomy may  not preserve the horizontal compositions of 2-morphisms. 
\subsubsection{The relation between the exact and fuzzy holonomies}
There is no evident general relation between the blur  $\widehat{\Fz}_{(\bw,\bm))}$ and  exact $\Fc_{(\bw,\bm)}$ holonomies of a fully primitive Hopf 2-connection $(\w,m_1,m_2)$. Nevertheless, the fuzzy holonomy can be derived from the exact holonomy.  

{Let $\H=(\d \colon I \to H, \rho)$ be a crossed module of Hopf algebras. 
Consider the crossed module inclusion of \ref{cminc}, a map 
${\rm Inc}_\H \colon  (\H)^*_{{\rm gl}} \to  {\cal BA}(\H)^*_\bullet$  of crossed modules of groups, therefore inducing a map of 2-groupoids:}
$$\Cc^\times\big(  (\H)^*_{{\rm gl}} \big) \ra{\Cc^\times\big({\rm Inc}_\H\big) } \Cc^\times \big( {\cal BA}(\H)^*_\bullet\big).  $$
The following theorem is one of the main results of this paper. Recall the notation of Theorem \ref{rel}:
\begin{Theorem}Let $(\w,m_1,m_2)$ be a fully primitive Hopf 2-connection. The  diagram below commutes:
 \begin{equation}\label{comp0}
  \hskip-2.5cm \xymatrix{
   &&\Cc^\times\Big(\big (\H\big)^*_{\rm gl}\Big) 
   \ar[d]^{\Cc^\times\left ({\rm Inc}_{\H}\right)}&&&\P(M)\ar[lll]_{ {\cal F}^\times_{(\bw,\bm)}  }   \ar[d]^{ \Fz^+_{(\bw,\bm)}  }   \\
   && \Cc^\times\Big({\cal BA}\left(\H\right)^*_\bullet \big) \ar[rrr]_{ {\cal T}_{ {\cal BA}\left (\H\right)}}  & & &\Cc^+\Big({\cal BA}\left(\H\right) \Big)
   } 
   \end{equation}
\end{Theorem}
\begin{Proof}    At the level of 1-paths this is tautological. Recall the construction of \ref{exa}. 
Let $\G=\morts{x}{\g_0}{\G}{\g_1}{y}$ be a 2-path, and as usual put $\g_s(t)=\G(t,s)$. For each $s\in [0,1]$, consider the 2-morphisms below,  of $\Cc^+ ({\BA(\H)})$:
 $$\hskip-1cm\morplush{{P_{\bw}}(\g_0)^{-1}}{  \underline{\left (\left(Q_{(\bw,\bm)}(\G,[0,s])\right) -1 \right)\tn {{P_{\bw}}(\g_0)^{-1}}  }}{{P_{\bw}}(\g_s)^{-1}} \an \hskip-1cm\morplush{{P_{\bw}}(\g_0)^{-1}}{\underline{R_{(\bw,\bm)}(\G,[0,s]) }} {{P_{\bw}}(\g_s)^{-1}} .$$
 These coincide with $ {\cal T}_{ {\cal BA}\left (\H\right)} \circ \Cc^\times\left ({\rm Inc}_{\H}\right) \circ {\cal F}^\times_{(\bw,\bm)}(\G)$ and  $\Fz^+_{(\bw,\bm)}(\G) $, respectively, 
  for $s=1$; see equations \eqref{relf} and \eqref{F2}. 
We  prove that, for each $s \in [0,1]$ (and thus for $s=1$), we have: 
$$ \underline{ \left (Q_{(\bw,\bm)}(\G,[0,s])) -1\right) \tn {({P_{\bw}}(\g_0)^{-1}) }} =\underline{R_{(\bw,\bm)}(\G,[0,s])}.$$
By definition, and using the second Peiffer law for crossed module of bare algebras in the fourth equality:
 \begin{align*}
  &\frac{d}{ ds }
  \underline{\left(Q_{(\bw,\bm)}(\G,[0,s]) -1 \right)\tn {{P_{\bw}}(\g_0)^{-1}}}
  = \underline{\left(\frac{d}{ ds }  Q_{(\bw,\bm)}(\G,[0,s])\right)\tn {{P_{\bw}}(\g_0)^{-1}}} \\
  & =-\underline{\left ( \Big(\int_0^1 {P_{\bw}}(\G,[0,u],s)\t \bm \left ( \frac{\d}{\d s} \G(u,s), \frac{\d}{\d u} \G(u,s) \right )  du\Big) Q_{(\bw,\bm)  }(\G,[s,s_0])\right) \tn {{P_{\bw}}(\g_0)^{-1}}}\\ 
  & =-\underline{\left ( \int_0^1 {P_{\bw}}(\G,[0,u],s)\t \bm \left ( \frac{\d}{\d s} \G(u,s), \frac{\d}{\d u} \G(u,s) \right )  du\right) \left (Q_{(\bw,\bm)  }(\G,[s,s_0])\right) \tn  {{P_{\bw}}(\g_0)^{-1}}}\\
  &=-\underline{\left ( \int_0^1 {P_{\bw}}(\G,[0,u],s)\t \bm \left ( \frac{\d}{\d s} \G(u,s), \frac{\d}{\d u} \G(u,s) \right )  du\right)\tn \d\left( \left (Q_{(\bw,\bm)  }(\G,[s,s_0])\right)\right) {({P_{\bw}}(\g_0)^{-1})}}\\
  &=-\underline{\left ( \int_0^1 {P_{\bw}}(\G,[0,u],s)\t \bm \left ( \frac{\d}{\d s} \G(u,s), \frac{\d}{\d u} \G(u,s) \right )  du\right)\tn \d\left( \left (Q_{(\bw,\bm)  }(\G,[s,s_0])\right) {P_{\bw}}(\g_0)^{-1} \right )}\\
   &=-\underline{\left ( \int_0^1 {P_{\bw}}(\G,[0,u],s)\t \bm \left ( \frac{\d}{\d s} \G(u,s), \frac{\d}{\d u} \G(u,s) \right )  du\right)\tn {{P_{\bw}}(\g_s)^{-1}} },
 \end{align*}
since $\bm$ is primitive and using \eqref{tytyt}. Whereas, by \eqref{lkjhl}, noting that  ${P_{\bw}}(\G,[0,1],s)={P_{\bw}}(\g_s)$:
  \begin{align*}
  &\frac{d}{ ds }\underline{\left(R_{(\bw,\bm)}(\G,[0,s]) \right)}=-\underline{\left(\int_0^1 \left ({P_{\bw}}(\G,[0,u],s) \t \bm \left ( \frac{\d}{\d s} \G(u,s), \frac{\d}{\d u} \G(u,s) \right )du \right) \tn {P_{\bw}}(\G,[0,1],s)\right)}\\
  &\doteq \underline{\left(\int_0^1 {P_{\bw}}(\G,[0,u],v) \t \bm \left ( \frac{\d}{\d s} \G(u,s), \frac{\d}{\d u} \G(u,s) \right )du\right) \tn  {P_{\bw}}(\g_s)^{-1}}.
   \end{align*}
  { Thence $\underline{\left(Q_{(\bw,\bm)}(\G,[0,s]) -1 \right)\tn {{P_{\bw}}(\g_0)^{-1}}}=\underline{R_{(\bw,\bm)}(\G,[0,s]) }$, for all $s \in [0,1]$.}
\end{Proof}

\subsubsection{The relation between the exact and bare holonomies}

Now consider a unital  crossed module $\A=(\d\colon A \to B, {\tra,\tla})$, of bare algebras. From $\A$ we can define the group crossed module $\A_\bullet^*=\big(\d \colon  \Ad^* \to B^*, \t \big)$, as explained in \ref{bag}. We can also form (see \ref{bdga}) the differential crossed module $\Lie(\A)=(\d \colon \Lie(A) \to \Lie(B), \t)$, and, therefore (see \ref{bcvcv})  the crossed module of Hopf algebras $\U(\Lie(\A))=(\d \colon \U(\Lie(A)) \to \U(\Lie(B)), \rho)$. From $\U(\Lie(\A))$, we define (see \ref{hag}) the crossed module of groups $\U(\Lie(\A))^*_{\rm gl}= (\d\colon \U(\Lie(A))^*_{\rm gl} \to \U(\Lie(B))^*_{\rm gl},\t)$. By Theorem \ref{zxzx}, we have a crossed module map ${\rm Proj}_\A=(\xi^\flat_A,\xi_B) \colon \U(\Lie(\A))^*_{\rm gl} \to \A_\bullet^*,$ called the crossed module counit.  
We  have (see \ref{gcmc}) 2-groupoids $\Cc^\times  ( {\A_\bullet^*})$ and $\Cc^\times  ( { \U(\Lie(\A))^*_{\rm gl} })$. The crossed module counit ${\rm Proj}_\A \colon \U(\Lie(\A))^*_{\rm gl} \to \A_\bullet^*$ induces a 2-functor {$ \Cc^\times \left ( { \U(\Lie(\A))^*_{\rm gl} } \right)\ra{\Cc^\times( {\rm Proj}_\A)} \Cc^\times \left ( {\A_\bullet^*}\right) $}.

Recall  \ref{qwer} that from  the Hopf algebra crossed module $\U(\Lie(\A))$ we can  define a crossed module of bare algebras $\BA(\U(\Lie(\A)))=\big(t\colon \underline{{\U^0(\Lie(\A))}\tn_\rho \U(\Lie(B))} \to \U(\Lie(B)), \tra,\tla\big)$.  {By using the construction in  \ref{bacmc}, we have 2-groupoids}  $\Cc^+ ({\BA(\U(\Lie(\A))})$ and $\Cc ^+({\A})$.  Theorem \ref{fgfgfg} gives us a map of  crossed modules of bare algebras:  $K_A\doteq (\kappa_\A, \xi_B)\colon \BA(\U(\Lie(\A)))\to \A$, the bare counit, thus a 2-groupoid map $\Cc^+\left ({\BA(\U(\Lie(\A))}\right) \ra{\Cc^+(K_\A)} \Cc^+\left ({\A}\right)$.

If {$(\bw,\bm)$} is a bare 2-connection in $\A$, seen as a pair of graded forms, we can also look at it as being a fully primitive Hopf 2-connection in  $\U(\Lie(\A))$.  We will thus consider the exact and fuzzy holonomies $\Fc_{(\bw,\bm)}^\times$ and ${\Fz^+_{(\bw,\bm)}}$, of  {$(\bw,\bm)$}.  These define maps $\F^\times_{(\bw,\bm)}\colon \P(M) \to  \Cc^\times \left ( { \U(\Lie(\A))^*_{\rm gl} } \right)$ and ${\Fz^+_{(\bw,\bm)}}\colon \P(M) \to\Cc\left ({\BA(\U(\Lie(\A))}\right)$, preserving all compositions, boundaries and inverses. Of course, we can also consider the bare holonomy $\F^+_{(\bw,\bm)}\colon \P(M) \to \Cc^+(\A)$ of $(\bw,\bm)$.

Now for one of the main theorems of this paper. With a different formulation, and in the case of crossed modules derived from chain-complexes of vector spaces, this was mentioned in \cite{CFM2} and in \cite{AF2} (the latter reference giving a detailed proof). 
\begin{Theorem}
 The following diagram commutes (once again recall  Theorem \ref{rel}):
 $$ {\xymatrix{  
 & & \P(M)\ar[dl]_{\F^\times_{(\bw,\bm)}} \ar[dr]^{ {\Fz^+_{(\bw,\bm)}}}  \ar@/_2pc/[ddr]_{\F^+_{(\bw,\bm)}}&    
 \\
 & \Cc^\times \left ( { \U(\Lie(\A))^*_{\rm gl} } \ar[d]_{\Cc^\times( {\rm Proj}_\A)}\right) &  &   \Cc^+\Big ({\BA(U(\Lie(\A))}\Big)\ar[d]^{\Cc^+(K_\A)}  \\
 & \Cc^\times \left ( {\A_\bullet^*} \right) \ar[rr]_{\Tc_\A} & &  \Cc^+\left ({\A}\right)
 }}$$
\end{Theorem}
\begin{Proof}
This follows directly from diagram \eqref{Comp2}, in the Introduction. The most difficult part of it  (the commutativity of the square  \eqref{comp0}) was already proven. All of the remaining polygons of the diagram \eqref{Comp2}, commute either by naturality, definition  or  by using Lemma  \ref{opopo}.
\end{Proof}
\bibliographystyle{alpha}
\bibliography{HAXM.bib}

\begin{thebibliography}{{Wag}06}

\bibitem[AC12]{AC}
Camilo~Arias {Abad} and Marius {Crainic}.
\newblock {Representations up to homotopy of Lie algebroids.}
\newblock {\em {J. Reine Angew. Math.}}, 663:91--126, 2012.

\bibitem[AP96]{PA}
Z.~{Arvasi} and T.~{Porter}.
\newblock {Simplicial and crossed resolutions of commutative algebras.}
\newblock {\em {J. Algebra}}, 181(2):426--448, art. no. 0128, 1996.

\bibitem[{Arv}04]{A}
Zekeriya {Arvasi}.
\newblock {Crossed modules of algebras.}
\newblock {\em {Math. Comput. Appl.}}, 9(2):173--182, 2004.

\bibitem[AS14a]{AF2}
Camilo~Arias Abad and Florian Schaetz.
\newblock Higher holonomies: comparing two constructions; arxiv:1404.0729.
\newblock {\em arXiv preprint arXiv:1404.0729}, 2014.

\bibitem[AS14b]{AF1}
Camilo {Arias Abad} and Florian {Sch\"atz}.
\newblock {Holonomies for connections with values in $L_\infty$-algebras.}
\newblock {\em {Homology Homotopy Appl.}}, 16(1):89--118, 2014.

\bibitem[{Bar}95]{BN}
Dror {Bar-Natan}.
\newblock {On the Vassiliev knot invariants.}
\newblock {\em {Topology}}, 34(2):423--472, 1995.

\bibitem[{Bau}91]{Baues}
Hans~Joachim {Baues}.
\newblock {\em {Combinatorial homotopy and 4-dimensional complexes.}}
\newblock Berlin etc.: Walter de Gruyter, 1991.

\bibitem[BC04]{BC}
John~C. {Baez} and Alissa~S. {Crans}.
\newblock {Higher-dimensional algebra. VI: Lie 2-algebras.}
\newblock {\em {Theory Appl. Categ.}}, 12:492--538, 2004.

\bibitem[BHS11]{BHS}
Ronald {Brown}, Philip~J. {Higgins}, and Rafael {Sivera}.
\newblock {\em {Nonabelian algebraic topology. Filtered spaces, crossed
  complexes, cubical homotopy groupoids. With contributions by Christopher D.
  Wensley and Sergei V. Soloviev.}}
\newblock Z\"urich: European Mathematical Society (EMS), 2011.

\bibitem[BL04]{BL}
John~C. {Baez} and Aaron~D. {Lauda}.
\newblock {Higher-dimensional algebra. V: 2-Groups.}
\newblock {\em {Theory Appl. Categ.}}, 12:423--491, 2004.

\bibitem[BM05]{BreenMessing}
Lawrence {Breen} and William {Messing}.
\newblock {Differential geometry of gerbes.}
\newblock {\em {Adv. Math.}}, 198(2):732--846, 2005.

\bibitem[BM06]{BaMa}
John~W. {Barrett} and Marco {Mackaay}.
\newblock {Categorical representations of categorical groups.}
\newblock {\em {Theory Appl. Categ.}}, 16:529--557, 2006.

\bibitem[Bro81]{BHiggins}
Higgins~Philip~J. Brown, Ronald.
\newblock The equivalence of $\infty $-groupoids and crossed complexes.
\newblock {\em Cahiers de Topologie et Géométrie Différentielle
  Catégoriques}, 22(4):371--386, 1981.

\bibitem[{Bro}82]{KB}
Kenneth~S. {Brown}.
\newblock {Cohomology of groups.}
\newblock {Graduate Texts in Mathematics, 87. New York-Heidelberg-Berlin:
  Springer-Verlag. X, 306 p., 4 figs. DM 74.00 {\$} 29.60 (1982).}, 1982.

\bibitem[{Bro}99]{Br}
Ronald {Brown}.
\newblock {Groupoids and crossed objects in algebraic topology.}
\newblock {\em {Homology Homotopy Appl.}}, 1:1--78, 1999.

\bibitem[BS76]{BrS}
Ronald {Brown} and Christopher~B. {Spencer}.
\newblock {Double groupoids and crossed modules.}
\newblock {\em {Cah. Topologie G\'eom. Diff\'er. Cat\'egoriques}}, 17:343--362,
  1976.

\bibitem[BS05]{BS}
John~C. Baez and Urs Schreiber.
\newblock {Higher gauge theory; arxiv math/0511710}.
\newblock 2005.

\bibitem[BSCS07]{BSAS}
John~C. Baez, Danny Stevenson, Alissa~S. Crans, and Urs Schreiber.
\newblock From loop groups to 2-groups.
\newblock {\em Homology Homotopy Appl.}, 9(2):101--135, 2007.

\bibitem[CF15]{CFM3}
Lucio~Simone {Cirio} and Jo{\~{a}}o {Faria Martins}.
\newblock {Infinitesimal 2-braidings and differential crossed modules.}
\newblock {\em {Adv. Math.}}, 277:426--491, 2015.

\bibitem[CFM12a]{CFM2}
Lucio~Simone Cirio and Jo\~{a}o Faria~Martins.
\newblock {Categorifying the $sl(2,C)$ Knizhnik-Zamolodchikov Connection via an
  Infinitesimal 2-Yang-Baxter Operator in the String Lie-2-Algebra; arxiv
  1207.1132}.
\newblock 2012.

\bibitem[CFM12b]{CFM1}
Lucio~Simone Cirio and Jo{\~a}o Faria~Martins.
\newblock Categorifying the knizhnik-zamolodchikov connection.
\newblock {\em {{Differ. Geom. Appl.}}}, {30}({3}):{238--261}, {2012}.

\bibitem[{Che}73]{Chen}
Kuo-Tsai {Chen}.
\newblock {Iterated integrals of differential forms and loop space homology.}
\newblock {\em {Ann. Math. (2)}}, 97:217--246, 1973.

\bibitem[CIKL14]{CIKL}
J.M. {Casas}, N.~{Inassaridze}, E.~{Khmaladze}, and M.~{Ladra}.
\newblock Adjunction between crossed modules of groups and algebras.
\newblock {\em Journal of Homotopy and Related Structures}, 9(1):223--237,
  2014.

\bibitem[CP94]{CP}
A.~{Caetano} and R.F. {Picken}.
\newblock {An axiomatic definition of holonomy.}
\newblock {\em {Int. J. Math.}}, 5(6):835--848, 1994.

\bibitem[DIKL12]{DIKL}
Guram {Donadze}, Nick {Inassaridze}, Emzar {Khmaladze}, and Manuel {Ladra}.
\newblock {Cyclic homologies of crossed modules of algebras.}
\newblock {\em {J. Noncommut. Geom.}}, 6(4):749--771, 2012.

\bibitem[Ell93]{El}
Graham~J. Ellis.
\newblock Homotopical aspects of lie algebras.
\newblock {\em Journal of the Australian Mathematical Society (Series A)},
  54:393--419, 6 1993.

\bibitem[FLV07]{FLV}
J.M. {Fern\'andez Vilaboa}, M.P. {L\'opez L\'opez}, and E.~{Villanueva Novoa}.
\newblock {Cat$^{1}$-Hopf algebras and crossed modules.}
\newblock {\em {Commun. Algebra}}, 35(1):181--191, 2007.

\bibitem[FMM11]{FMM}
Jo\~{a}o Faria~Martins and Aleksandar Mikovic.
\newblock {Lie crossed modules and gauge-invariant actions for 2-BF theories}.
\newblock {\em Adv.Theor.Math.Phys.}, 15:1059--1084, 2011.

\bibitem[FP10]{FMP1}
Jo{\~{a}}o {Faria Martins} and Roger {Picken}.
\newblock {On two-dimensional holonomy.}
\newblock {\em {Trans. Am. Math. Soc.}}, 362(11):5657--5695, 2010.

\bibitem[FP11a]{FMP2}
Jo{\~{a}}o {Faria Martins} and Roger {Picken}.
\newblock {Surface holonomy for non-Abelian 2-bundles via double groupoids.}
\newblock {\em {Adv. Math.}}, 226(4):3309--3366, 2011.

\bibitem[FP11b]{FMP3}
Jo{\~{a}}o {Faria Martins} and Roger {Picken}.
\newblock {The fundamental Gray 3-groupoid of a smooth manifold and local
  3-dimensional holonomy based on a 2-crossed module.}
\newblock {\em {Differ. Geom. Appl.}}, 29(2):179--206, 2011.

\bibitem[FW11]{FW}
Ya\"el {Fr\'egier} and Friedrich {Wagemann}.
\newblock {On Hopf 2-algebras.}
\newblock {\em {Int. Math. Res. Not.}}, 2011(15):3471--3501, 2011.

\bibitem[{Ger}66]{Ger}
M.~{Gerstenhaber}.
\newblock {On the deformation of rings and algebras. II.}
\newblock {\em {Ann. Math. (2)}}, 84:1--19, 1966.

\bibitem[HKK00]{HKK}
K.A. {Hardie}, K.H. {Kamps}, and R.W. {Kieboom}.
\newblock {A homotopy 2-groupoid of a Hausdorff space.}
\newblock {\em {Appl. Categ. Struct.}}, 8(1-2):209--234, 2000.

\bibitem[{Kas}95]{Kassel}
Christian {Kassel}.
\newblock {\em {Quantum groups.}}
\newblock New York, NY: Springer-Verlag, 1995.

\bibitem[{Kon}93]{Ko}
Maxim {Kontsevich}.
\newblock {Vassiliev's knot invariants.}
\newblock In {\em {I. M. Gelfand seminar. Part 2: Papers of the Gelfand seminar
  in functional analysis held at Moscow University, Russia, September 1993}},
  pages 137--150. Providence, RI: American Mathematical Society, 1993.

\bibitem[KP02]{KP}
K.H. {Kamps} and T.~{Porter}.
\newblock {2-groupoid enrichments in homotopy theory and algebra.}
\newblock {\em {$K$-Theory}}, 25(4):373--409, 2002.

\bibitem[{Mac}98]{McLane}
Saunders {Mac Lane}.
\newblock {\em {Categories for the working mathematician. 2nd ed.}}
\newblock New York, NY: Springer, 2nd ed edition, 1998.

\bibitem[{Maj}95]{MajidQG}
Shahn {Majid}.
\newblock {\em {Foundations of quantum group theory.}}
\newblock Cambridge: Cambridge Univ. Press, 1995.

\bibitem[Maj12]{MajidHA}
Shahn Majid.
\newblock {Strict quantum 2-groups; arxiv 1208.6265}.
\newblock 2012.

\bibitem[{Mur}88]{Mu}
M.K. {Murray}.
\newblock {Another construction of the central extension of the loop group.}
\newblock {\em {Commun. Math. Phys.}}, 116(1):73--80, 1988.

\bibitem[MW50]{MLW}
Saunders {MacLane} and J.H.C. {Whitehead}.
\newblock {On the 3-type of a complex.}
\newblock {\em {Proc. Natl. Acad. Sci. USA}}, 36:41--48, 1950.

\bibitem[PS86]{PS}
Andrew {Pressley} and Graeme {Segal}.
\newblock {Loop groups.}
\newblock {Oxford Mathematical Monographs. Oxford: Clarendon Press. VIII, 318
  p., (1986).}, 1986.

\bibitem[Str96]{St}
Ross Street.
\newblock {Categorical structures.}
\newblock {Hazewinkel, M. (ed.), Handbook of algebra. Volume 1. Amsterdam:
  North-Holland. 529-577 (1996).}, 1996.

\bibitem[SW11]{SW1}
Urs {Schreiber} and Konrad {Waldorf}.
\newblock {Smooth functors vs. differential forms.}
\newblock {\em {Homology Homotopy Appl.}}, 13(1):143--203, 2011.

\bibitem[SW13]{SW2}
Urs {Schreiber} and Konrad {Waldorf}.
\newblock {Connections on non-abelian gerbes and their holonomy.}
\newblock {\em {Theory Appl. Categ.}}, 28:476--540, 2013.

\bibitem[Var84]{Var}
V.S. Varadarajan.
\newblock {\em {Lie groups, Lie algebras, and their representations. Reprint of
  the 1974 edition.}}
\newblock {Graduate Texts in Mathematics. 102. New York, NY: Springer. xiii,
  430 p.}, 1984.

\bibitem[{Wag}06]{W}
Friedrich {Wagemann}.
\newblock {On Lie algebra crossed modules.}
\newblock {\em {Commun. Algebra}}, 34(5):1699--1722, 2006.

\bibitem[Whi41]{WJHC}
J.H.C. Whitehead.
\newblock {On adding relations to homotopy groups.}
\newblock {\em Ann. of Math. (2)}, 1941.

\bibitem[WW]{WoW}
Friedrich Wagemann and Christoph Wockel.
\newblock A cocycle model for topological and lie group cohomology.
\newblock {\em Transactions of the American Mathematical Society (to appear).}

\end{thebibliography}

\end{document}